\documentclass[11pt, reqno]{amsart}
\usepackage{a4wide}
\usepackage{wrapfig}
\usepackage{amsmath}
\usepackage{amsthm} 
\usepackage{mathrsfs}
\usepackage{xcolor}
\usepackage{graphicx}
\usepackage{amssymb}
\usepackage{ dsfont }
\usepackage{url}	
\usepackage{multicol}
\usepackage{tikz}
\usepackage{amsmath}
\usepackage{fancyhdr} 
\usetikzlibrary{matrix}
\usetikzlibrary{arrows}
\usepackage{subfig}
\usepackage{hyperref}
\usepackage{listings}

\usepackage{esint}

\usepackage{varwidth}

\allowdisplaybreaks

\usepackage{esint,amssymb}
\usepackage{enumerate}
\usepackage{mathrsfs}

\makeatletter
\@namedef{subjclassname@2020}{\textup{\it 2020} Mathematics Subject Classification}
\makeatother

\usepackage[latin1]{inputenc}
\usepackage[T1]{fontenc} 
\usepackage[english]{babel}
\usepackage{verbatim}
\usepackage{calligra}
\usepackage{dsfont}
\usepackage{amsfonts}
\usepackage{amsmath}
\usepackage{amsthm}
\usepackage{amssymb}
\usepackage{mathtools}
\usepackage{mathrsfs}
\usepackage{color}
\usepackage{graphicx}
\usepackage{bm}
\usepackage{srcltx}

\usepackage{tikz}
\usepackage[retainorgcmds]{IEEEtrantools}

\usepackage{graphicx}		
\usepackage[hypcap=false]{caption}

\renewcommand{\div}{\operatorname{div}}

\newcommand{\tr}{\operatorname{tr}}

\theoremstyle{definition}
\newtheorem{definition}{Definition}[section]

\newtheorem{rmk}[definition]{Remark}

\theoremstyle{plane}
\newtheorem{thm}[definition]{Theorem}
\newtheorem{prop}[definition]{Proposition}
\newtheorem{cor}[definition]{Corollary}
\newtheorem{lemma}[definition]{Lemma}

\newcommand{\tsc}{\textsc}
\newcommand{\tsl}{\textsl}

\newcommand{\mbb}{\mathbb}

\newcommand{\mc}{\mathcal}
\newcommand{\mf}{\mathfrak}

\newcommand{\veps}{\varepsilon}

\newcommand{\wtilde}{\widetilde}
\newcommand{\vphi}{\varphi}
\newcommand{\oline}{\overline}

\newcommand{\ra}{\rightarrow}

\newcommand{\g}{\gamma}
\renewcommand{\k}{\kappa}
\newcommand{\s}{\sigma}
\renewcommand{\t}{\tau}
\newcommand{\z}{\zeta}

\newcommand{\de}{\delta}

\newcommand{\bull}{\bullet}

\newcommand{\R}{\mathbb{R}}
\newcommand{\Q}{\mathbb{Q}}

\newcommand{\N}{\mathbb{N}}
\newcommand{\Z}{\mathbb{Z}}

\renewcommand{\P}{\mathbb{P}}
\newcommand{\D}{\mathbb{D}}

\renewcommand{\div}{{\rm div}\,}

\newcommand{\Id}{{\rm Id}\,}

\newcommand{\adj}{{\rm adj}\,}

\newcommand{\dt}{ \, {\rm d} t}

\newcommand{\loc}{{\rm loc}}

\allowdisplaybreaks

\def\d{\partial}
\def\div{{\rm div}}
\def\divergence{{\rm div}}
\newcommand{\dd}{{\rm d}}

\newcommand{\be}{\begin{equation}}
\newcommand{\ee}{\end{equation}}
\newcommand{\bp}{\begin{proof}}
\newcommand{\ep}{\end{proof}}
\newcommand{\bel}{\begin{equation}\label}
\newcommand{\eeq}{\end{equation}}
\newcommand{\bea}{\begin{eqnarray}}
\newcommand{\eea}{\end{eqnarray}}
\newcommand{\bee}{\begin{eqnarray*}}
\newcommand{\eee}{\end{eqnarray*}}
\newcommand{\ben}{\begin{enumerate}}
\newcommand{\een}{\end{enumerate}}

\newcommand{\la}{\lambda}

\definecolor{codegreen}{rgb}{0,0.6,0}
\definecolor{codegray}{rgb}{0.5,0.5,0.5}
\definecolor{codepurple}{rgb}{0.58,0,0.82}
\definecolor{backcolour}{rgb}{0.95,0.95,0.92}

\lstdefinestyle{mystyle}{
	backgroundcolor=\color{backcolour},
	commentstyle=\color{codegreen},
	keywordstyle=\color{magenta},
	numberstyle=\tiny\color{codegray},
	stringstyle=\color{codepurple},
	basicstyle=\footnotesize,
	breakatwhitespace=false,
	breaklines=true,
	captionpos=b,
	keepspaces=true,
	numbers=left,
	numbersep=5pt,
	showspaces=false,
	showstringspaces=false,
	showtabs=false,
	tabsize=2
}
\numberwithin{equation}{section}

\usepackage{pgfplots}
\pgfplotsset{compat=newest}


\theoremstyle{definition}

\numberwithin{ej}{section}

\author[F. Fanelli]{Francesco Fanelli}
\thanks{The work of F.F. has been partially supported by the project CRISIS (ANR-20-CE40-0020-01), operated by the French National Research Agency (ANR),
and by the project GHOST (PID2025-169177NB-I00), operated by the Spanish State Research Agency (AEI).
F.F. also aknowledges the support of the Basque Government through the BERC 2022-2025 program, of the Spanish State Research Agency (AEI)
through the BCAM Severo Ochoa excellence accreditation CEX2021-001142, and of the European Union through the COFUND program [HORIZON-MSCA-2022-COFUND-101126600-SmartBRAIN3].}
\address[F. Fanelli]{Basque Center for Applied Mathematics, Alameda de Mazarredo 14, E-48009 Bilbao, Basque Country, SPAIN \newline
\& Ikerbasque -- Basque Foundation for Science, Plaza Euskadi 5, E-48009 Bilbao, Basque Country, SPAIN} 
\email{ffanelli@bcamath.org}

\author[M. E. Mart\'inez-Martini]{M. Eugenia Mart\'inez-Martini}
\thanks{M. E. M. M. was supported by Centro de Modelamiento Matem\'atico (CMM) BASAL fund FB210005 for center of excellence from ANID-Chile. Part of this work was done during M. E. M. M. postdoctoral position in the framework of the project "HiCE", funded by the Universit\'e Claude Bernard Lyon 1 within the programme SENS 2022.}
\address[M. E. Mart\'inez-Martini]{Cente of Mathematical Modeling, Universidad de Chile, Casilla 170, Correo 3, Beauchef 851, Santiago, CHILE}
\email{maria.martinez.m@uchile.cl}

\begin{document}


\title[Navier-Stokes-Fourier system with discontinuous densities]{Local well-posedness for viscous heat-conducting fluids with discontinuous densities}

\keywords{Navier-Stokes-Fourier system; discontinous density; local well-posedness; maximal regularity; striated regularity. \vspace{.1cm}}

\date{\today}

\begin{abstract}
In this paper, we consider the Navier-Stokes-Fourier system in space dimension $d\geq2$. We prove local well-posedness
of solutions having density only in $L^p\cap L^\infty$, for suitable $p\in\,]1,+\infty[\,$, with no further regularity requirements.
In particular, our result yields existence and uniqueness of solutions in the case of discontinuous densities presenting a jump across a smooth interface.

The existence proof is based on the use of the Hoff effective viscous flux and on a combination of both maximal regularity and Chemin's striated regularity
arguments. Uniqueness, under the assumption of slightly more regular initial velocity fields and temperatures, relies on a Lagrangian approach,
supplemented with an original bootstrap argument which shows improved time integrability properties of the solution.

\end{abstract}

\maketitle


\medskip
\section{Introduction and main result}

The goal of this work is to study the \emph{full compressible Navier-Stokes} system in space dimension $d\geq2$, in a low regularity functional framework
which is able to capture possible singularities (\tsl{i.e.} discontinuities) in the density function.

\subsection{The system of equations} \label{ss:system}

The full compressible Navier-Stokes system describes the motion of a viscous, compressible and heat-conducting fluid.
The state variables of the system are $\rho(t,x)\geq0$ the density of the fluid, $u(t,x)\in \R^d$ its velocity field and $\theta(t,x)\geq0$ its temperature.

Here, $d\in\N\setminus\{0\}$ denotes the space dimension, $t\geq0$ the time variable and $x\in\Omega\subset\R^d$ space variable. 
In this work, we decide to focus on the physically relevant cases $d=2,3$, but actually our arguments work for any $d\geq2$.
Throughout this paper we assume
\[
(t,x)\,\in\,\R_+\times\R^d\,,\qquad\quad \mbox{ with }\quad d\geq2\,.
\]

In order to write down the equations of motion, let us introduce
the scalar fields $P = P(\rho,\theta)$ and $e = e(\rho,\theta)$ to represent, respectively, the pressure and the internal energy per unit mass
of the fluid. These are functions of the density and temperature, and are determined through appropriate equations of state which will be specified below.
Finally, we also introduce the total energy per unit mass $E(t,x)\geq0$ as the scalar field
\[
 E\,:=\,\frac{1}{2}\,|u|^2\,+\,e\,.
\]
Fixed an initial datum $\big(\rho, \rho u, \rho E\big)_{|t=0}\,=\,\big(\rho_0, m_0, \Pi_0\big)$, the system of equations reads
\begin{equation}\label{FNS}
\left\{ 
	\begin{aligned}
	&\partial_t \rho\, +\, \divergence (\rho\, u)\,=\,0  
	\\[1ex]
	&\partial_t (\rho\, u)\, +\, \divergence(\rho\, u \otimes u)\,+\,\nabla P\, -\, \divergence \tau\,=\,0 
	\\[1ex]
	&\partial_t \left( \rho\, E\right)\, +\, \divergence \big( \left( \rho \,E\,+\, P\right)\,u \big)\,+\,\div q\,=\, \divergence(\tau \cdot u )\,,
\end{aligned}	
\right.
\end{equation}
where $\t$ is the viscous stress tensor and $q$ represents the heat flux.

System \eqref{FNS} is complemented with the following constitutive relations.
First of all, throughout this work we shall focus on the case of Newtownian fluids: we assume the viscous stress tensor $\tau$ to be linear in the velocity
gradients and given by
\begin{equation} \label{def:tau}
	\tau\,:= \,2\,\mu\,\D(u)\,+\,\lambda\,\div u\,\Id\,,
\end{equation} 
where $\Id$ denotes the $d\times d$ identity matrix and $\D(u)$ denotes the symmetric part of the Jacobian matrix of $u$, defined as 
\[
	\D(u)\,:=\, \frac12\,  \big( Du + \nabla u \big)\,, \qquad\qquad 
\left( Du\right)_{ij}\,=\, \partial_ju^i \quad \mbox{ and } \quad \left( \nabla u \right)_{ij}\,=\,\partial_i u^j\,.	
\]
The viscosity coefficients $\lambda$ and $\mu$ are assumed to be constant and to satisfy 
\[
\mu\,>\,0\qquad  \mbox{ and }\qquad \nu\,:=\,\lambda\, +\, 2\,\mu\,>\,0\,,	
\]
conditions which ensure the ellipticity of Lam\'e operator, namely the second order operator
\begin{equation}\label{DefinitionL}
	\mathcal L\,:=\, -\,\mu\,\Delta\, -\, \left(\lambda+\mu\right)\,\nabla \divergence\,.
\end{equation}
Observe that one has the equality
\[
-\,\div\tau\,=\,\mc Lu\,.
\]

Next, we assume the heat flux $q$ to be given by the Fourier law
\[
q\,:=\,-\, \kappa\,\nabla \theta\,,
\]
where we suppose $\kappa \,>\,0$ to be constant. In fact, we believe to be possible, within our approach, to consider the case of heat conductivity
coefficients $\k=\k(\rho)$ depending on the density. However, for simplicity, we stick to the case
$\k\,\equiv\,{\rm cst}$ constant here.

Finally, we specify the structural assumptions on the pressure $P(\rho,\theta)$ and internal energy $e(\rho,\theta)$.
Following \cite{Chik-D}, here we restrict our attention to the case where the pressure is given by
\begin{equation}\label{DefinitionP}
	P( \rho, \theta ): = \pi_0 ( \rho ) + \theta \pi_1 ( \rho ),
\end{equation}
where $\pi_0$ and $\pi_1$ are given smooth functions of the density only.
At the same time, we assume that the internal energy $e$ depends linearly on the temperature $\theta$. Then, taking into account the Maxwell relation
\[
	\rho^2\, \partial_\rho e\, =\, P\, -\, \theta\, \partial_\theta P\,,	
\]
we find the explicit expression
\begin{equation}\label{definitione}
	e(\rho,\theta)\,=\, c_v\, \theta\, +\, h(\rho)\,, 	
\end{equation}
where $c_v>0$ is the specific heat at constant volume and the smooth function $h$ satisfies the relation
\begin{equation}\label{RelationBetweenPih}
	\rho^2\, h'(\rho)\,=\, \pi_0 (\rho)\,. 	
\end{equation}
For simplicity of exposition, from now on we will assume $c_v=1$.

\subsection{Overview of the literature} \label{ss:overview}

The mathematical study of compressible viscous and heat-conducting fluids has been the subject of extensive research over the past decades.
Local well-posedness of classical (say, $\mc C^\infty$) solutions was first established by Nash \cite{Nash}; later, Itaya \cite{It_1971, It_1976} considered
the framework of H\"older continuous solutions.
The first global well-posedness result for strong solutions dates back to Matsumura and Nishida \cite{Mat-Ni}, who solved the system near the trivial equilibrium
in Sobolev spaces $H^s$ of high enough regularity.
Subsequently, Danchin \cite{D-2001-loc, D-2001_glob} developed a well-posedness theory in critical functional spaces of Besov type.
This approach was then extended by Chikami and Danchin \cite{Chik-D}, who used a Lagrangian formulation approach to improve the theory (especially concerning
uniqueness results) of the full compressible Navier-Stokes system in critical Besov spaces.

In comparison, results concerning existence weak solutions, which are usually based on a physically relevant energy inequality and the notion of finite energy,
are much more recent. First results are due to Hoff \cite{Hoff_1997} (which we will comment more in detail below)
and Feireisl \cite{Feir_book, Feir_IUMJ}. Preliminary studies, where one suitably modifies the constitutive laws for large values
of the density and temperature, can be found in P.-L. Lions's book \cite{Lions}. We refer to the book \cite{F-N_book} of Feireisl and Novotn\'y
for a general introduction on the topic, precise
statements and proofs, and additional references. It must be noted that all these results are global in time, but they only concern existence of solutions,
whereas they lack, in general, of a uniqueness statement, despite partial weak-strong uniqueness principles have been obtained (see \tsl{e.g.} \cite{F-N_2012}).

More recent studies on the full compressible Navier-Stokes system \eqref{FNS} concern the case of open systems, namely of non-homogeneous boundary conditions,
especially in the context of the celebrated Rayleigh-B\'enard problem.
They focus on both weak solutions theory \cite{F-N_open, Chaud-F} and their asymptotic behaviour for large times \cite{F-SG, Ch-F},
and on strong solutions theory \cite{A-B-C-F}.

\subsection{Goal of the study} \label{ss:goal}
Inspired by recent advances in the case of isentropic flows (which we will quickly review in Subsection \ref{ss:barot-NS}), in this paper we address
the question of the well-posedness of equations \eqref{FNS} in presence of \emph{density discontinuities}, like for instance density
jumps across a smooth interface. More precisely, we look for a theory able to ensure both \emph{uniqueness of solutions} and a
\emph{qualitative description} of the evolution of those discontinuities.

It must be noted that this question cannot be solved in a strong solutions framework, which is not able to capture density discontinuities,
owing to the high regularity demanded to solve the problem in that setting.
At the same time, all the results mentioned above concerning the weak solution theory, despite being able to consider discontinuous density functions,
do not provide enough qualitative information to track the dynamics of the singularities (namely, of the region of discontinuity of the density), nor do they
provide enough regularity to establish uniqueness (as already mentioned above).

In the present work, we are able to build up a well-posedness theory, albeit only local in time, which reach the above mentioned goals:
we obtain existence and uniqueness of solutions, together with a description of the evolution of the density discontinuities.

In order to better motivate and explain our result, which we will rigorously state in Subsection \ref{ss:NSF-temp},
let us make an \tsl{excursus} on the case of isentropic flows, in which one assumes that the entropy of the system is constant.
In particular, by Gibb's law, this implies that the internal energy $e$ of the fluid only depends on the density $\rho$, thus
the same holds for the pressure function $P$, owing to \eqref{definitione}. At the mathematical level, one simply neglects the total energy balance,
namely the last equation appearing in system \eqref{FNS}, and only considers the first two equations appearing therein.

\subsection{The case of isentropic flows} \label{ss:barot-NS}
Of course, the same questions mentioned above for system \eqref{FNS} arise also for the isentropic case, in which the equations simplify to the
so-called barotropic Navier-Stokes system
\begin{equation} \label{eq:bar-NS}
\left\{  \begin{array}{l}
          \d_t\rho\,+\,\div\big(\rho\,u\big)\,=\,0 \\[1ex]
          \d_t\big(\rho\, u\big)\, +\, \div(\rho\, u \otimes u)\,+\,\nabla P(\rho)\, +\,\mc Lu \,=\,0\,.
         \end{array}
\right.
\end{equation}

For this system, besides the Lions-Feireisl theory of finite energy weak solutions \cite{Lions, Feir_2001}, 
another theory of finite energy weak solutions was developed by Hoff \cite{Hoff_1991, Hoff_JDE}, under the assumption of initial data having
small initial energy, so being sufficiently close to the equilibrium $\big(\rho_{\rm eq}, u_{\rm eq}\big)\,=\,(\oline\rho,0)$,
for some constant $\oline\rho>0$.
Building on previous observations due to Serre \cite{Serre}, the key idea of Hoff theory was to use the (nowadays celebrated) \emph{effective viscous flux},
together with the small energy assumption and time-weighted energy estimates, to improve the qualitative information available for the constructed class of solutions.
In passing, we mention that the use of the effective viscous flux represents a crucial ingredient also of the Lions-Feireisl theory of weak solutions
(which, however, may have arbitrarily large initial data, not necessarily close to equilibrium).

Despite those additional properties, Hoff were only able to prove the existence of solutions within his framework, at least in the general
case $d\geq2$ (the case $d=1$ is a bit special and was dealt with in \cite{Hoff_1991}). Later on, some improvements
were obtained: for planar flows $d=2$, and under a linear assumption on the pressure function $P(\rho)\approx \rho$,
he established \cite{Hoff_2002} propagation of density discontinuities, globally in time, along $\mc C^{1,\alpha}$ curves. Uniqueness of such solutions was handled
only afterwards, in \cite{Hoff_2006}, but required stringent assumptions on the velocity field
(namely, to be Lipschitz continuous, a fact which was only guaranteed in $2$-D \cite{Hoff_2002}) and either on the pressure (which was required to be linear)
or on the density (roughly speaking, a condition $\nabla\rho\in L^p$ was needed, which is of course incompatible with jump discontinuities).

A general well-posedness result for system \eqref{eq:bar-NS}, ensuring uniqueness and propagation of singularities (namely, density discontinuities)
in any space dimension $d\geq1$ and for general pressure laws, was first given in \cite{DFP}, but only for short times. The analysis of \cite{DFP}
mixed arguments from maximal regularity, used to propagate smoothness on the effective viscous flux, and striated regularity \tsl{\`a la Chemin} \cite{Ch1991, Ch1993},
employed to gain Lipschitz continuity of the velocity field despite the singularities appearing in the density function.
As a matter of fact, the Lipschitz continuity of $u$ was crucial both for uniqueness and for propagation of density jumps.
Later on, the result of \cite{DFP} was made global in time by Liao and Zodji \cite{L-Z}, who resorted to Hoff's original approach
(energy estimates with time weights) and combined it with the striated regularity argument. Zodji further extended this analysis
to two-fluid systems \cite{Z-2025} and to the case of system \eqref{eq:bar-NS} in presence of density-dependent viscosity coefficients \cite{Z-glob},
that is $\mu=\mu(\rho)$ and $\lambda=\lambda(\rho)$ in definition \eqref{def:tau} of the viscous stress tensor $\t$.

\medbreak
Going back to the non-isentropic system \eqref{FNS}, Hoff's work \cite{Hoff_1997} represented an extension of his theory from \cite{Hoff_JDE}
to the case of entropy variations of the fluid. However, despite a series of better qualitative properties of solutions
was established in \cite{Hoff_1997}, to the best of our knowledge uniqueness still remained an elusive question within his theory for \eqref{FNS},
as well as the question of propagation of density jumps across interfaces.

\subsection{Reformulation of the system and main result} \label{ss:NSF-temp}

As already mentioned, in the present work we prove local in time well-posedness for system \eqref{FNS} in a functional framework
able to ensure uniqueness of solutions, as well as propagation of singularities (jump discontinuities) in the density functions.
Before stating our main result, some preparation is in order.

With the goal of getting full gain from the presence of the heat flux term, it is better to write the total energy
balance 
in terms of the temperature $\theta$. The precise computations, which are rigorously justified in our setting, are detailed in Appendix \ref{app:equiv}.
In turn, we can recast equations \eqref{FNS} in terms of the set of unkowns $\big(\rho, u, \theta\big)$ as follows:
\begin{equation}\label{eq:NSF}
	\left\{ 
		\begin{aligned}
		&\partial_t \rho + \div(\rho\, u)=0 
		\\[1ex]
		&\partial_t(\rho\,u)\, +\,\div(\rho\,u\otimes u)\,+\,\mc Lu\, +\, \nabla P\,=\,0		\\[1ex]
		&\partial_t (\rho\,\theta)\, +\,\div(\rho\,\theta\,u)\, +\,\pi_1(\rho)\,\theta\, \divergence u \,-\,\k\,\Delta\theta\,=\,\mc B(\nabla u)\,,
	\end{aligned}	
	\right.
\end{equation}	
where $\mc L$ has been introduced in \eqref{DefinitionL} and where we have defined the bilinear operator $\mc B$ as
\begin{equation} \label{def:B_u}
 \mc B(\nabla u)\,:=\,\lambda\, |\divergence u|^2\, +\, \mu\, \vert \nabla u \vert^2\, +\, \mu\, \nabla u : \nabla u\,,
\end{equation}
with, for given matrices $A$ and $B$, the notation $A:B := \sum_{j,k} A_{jk} B_{kj}\,=\,\tr(AB)$.
We dub system \eqref{eq:NSF}-\eqref{def:B_u} the \emph{Navier-Stokes-Fourier system} and we are going to use this denomination throughout this work.

Since we will consider the case of small density variations close to a constant $\oline\rho>0$, in particular the initial density will have a positive lower bound,
which we will be able to propagate for later times. Thus,
we can safely supplement system \eqref{eq:NSF} with the initial datum
\begin{equation*}
		\big(\rho, u, \theta\big)_{|t=0}\,=\, \big(\rho_0, u_0, \theta_0\big)\,,
\end{equation*}
for which suitable assumptions will be better specified below.

Let us also explicitly state some conditions and conventions.
First of all, we append to system \eqref{eq:NSF} the folowing boundary conditions at spatial infinity:
\[
u \to 0\,, \qquad \theta \to 0 \qquad \mbox{ and } \qquad \rho\to \oline \rho\qquad\qquad \mbox{ for }\qquad |x|\to+\infty\,,
\] 
where $\oline\rho>0$ is the (constant) density reference state. Without loss of generality, from now on we simply choose $\oline \rho\,=\,1$. 
Next, again with no loss of any generality, we assume that $\pi_0(1)=0$, so that $\pi_0(\rho)\approx \varpi(\rho)\,(\rho-1)$
and $P(\rho,\theta)$ satisfies $P(1,0)=0$.

Finally, in stating our main result, we need to make use of some notation and concepts properly defined later, especially in Section \ref{s:tools}.
In order to make the reading of the statement clearer, we point out the following:
\begin{itemize}
\item $B^{s}_{p,r}$ stands for the \emph{non-homogeneous} Besov space, $\dot B^s_{p,r}$ for the \emph{homogeneous} Besov space;
both are given in Definition \ref{d:B};
\item $\mbb L^{\infty,p}$ denotes the set of vector fields $X\in L^\infty$ such that $\nabla X\in L^p$, see \eqref{d:SetVectorField};
\item the notion of \emph{non-degenerate family} $\mc X\,=\,\bigl(X_{\la}\bigr)_{1\leq\la\leq m}$ of vector-fields is given in Definition \ref{d:NonDeg};
\item the directional derivative $\d_{X}$ along a vector field $X$ is introduced in Definition \ref{d:DDField};
\item given a family of vector fields $\mc X$, the space $\mbb L^{p}_{\mc X}$ is formed by functions $f\in L^\infty$ which are more regular
(namely, $\d_Xf \in L^p$) along all the directions $X\in\mc X$;
the associated norm is denoted by $\Vert \cdot \Vert_{\mbb L^{p}_{\mc X}}$; this is made precise in Definition \ref{d:LpdirecX}.
\end{itemize}

We are now ready to state the main result of this paper.

\begin{thm} \label{th:well-p}
Let the space dimension $d\geq 2$.
Let the indices $p$, $r$ and $s$ in $\,]1,+\infty[\,$ and $p_0\in\,]1,+\infty]$ be such that 
\begin{align}
\label{HypForPs}
& p>d\,, \qquad p_0> p\,, \qquad  1\,-\,\dfrac d{2p}\,<\,\frac1r\,<\,1\,-\,\frac d{4p}\,, \\
\label{HypForS}
&\quad 
\qquad\qquad\qquad
\dfrac2r\,-\,1 \,<\, \dfrac 1s\, <\, \dfrac1r
\quad \mbox{ and }\quad
\frac{1}{r}\,+\,\left(1\,+\,\frac{d}{p}\right)\,\frac{1}{2\,s}\,<\,\frac32\,.
\end{align}
Moreover, define the indices $r_0$ and $r_1$ by the formulas
\begin{equation}\label{DefOfRs}
\frac1{r_0}\,=\, \frac1r\,+\,\frac d{2p}\,-\,\frac d{2p_0}\,-\,1 \qquad \mbox{ and }\qquad \frac1{r_1}\,=\,\frac1r\,-\,\frac12\,,
\end{equation}
where we implicitly assume that $p_0$ in \eqref{HypForPs} is so large that $1/r_0>0$.

\paragraph*{\bf \tsc{Existence.}}
There exists $\veps_0>0$ such that the following property holds true.
Take an initial datum $\big(\rho_0,u_0,\theta_0\big)$ such that the following conditions are satisfied:
\begin{itemize}
 \item $\z_0\,:=\,\rho_0\,-\,1\;\in\,\big(L^p\cap L^\infty\big)(\R^d)$;
 \item $\theta_0\,\in\, \big(B^{2-2/s}_{p/2,s}\cap L^p\big)(\R^d)$;
 \item $w_0\,:=\,u_0\,-\,v_0\; \in\, \dot B^{2-2/r}_{p,r}(\R^d)$, where we have set
$v_0\,:=\,-\,\frac{1}{\nu}\,\nabla(\Id-\Delta)^{-1}P(\rho_0,\theta_0)$, with $\nu:=\lambda+2\mu$.
\end{itemize}
Assume that, for some $\veps\in\,]0,\veps_0]$, one has
\begin{equation}\label{th_HypothesisForZetaPropositionEstimatesForModifiedSolution}
\Vert \zeta_0 \Vert_{L^\infty}\,\leq\,\veps\,.	
\end{equation}
Finally, assume that  there exists  a non-degenerate family   $\mc X_0\,=\,\bigl(X_{0,\la}\bigr)_{1\leq\la\leq m}$ of vector-fields in
$\mbb L^{\infty,p}$ such that $\rho_0$ belongs to  $\mbb L^{p}_{\mc X_0}$.

Then, there exists a time $T>0$ and a weak solution $\big(\rho, u, \theta\big)$ to system \eqref{eq:NSF} defined on $[0, T]\times\R^d$, such that 
\begin{equation}\label{th_TesisLInftyZeta}
\Vert \rho-1\Vert_{L_T^\infty\left( L^\infty\right)}\,\leq\,4\,\varepsilon\,
\end{equation}
and such that the following regularity properties are satisfied:
\begin{enumerate}
\item $\zeta\,:=\,\rho\,-\,1\;\in L^\infty\big([0,T]\times\R^d\big)\,\cap\,\mc C\big([0,T]; L^{q}(\R^d)\big)$ for any $q\in[p,+\infty[\,$;
\item $\theta \in \mc C\big([0,T]; B^{2-2/s}_{p/2,s}\cap L^{p}(\R^d)\big)$, with $\partial_t \theta, \nabla^2\theta\; \in L^s\big([0,T]; L^{p/2}(\R^d)\big)$;
	\item $u\,=\,w\,+\,v$, with $v\,:=\,-\,\frac{1}{\nu}\,\nabla(\Id-\Delta)^{-1}P(\rho,\theta)$ belonging
to $\mc C\big([0,T]; W^{1,q}(\R^d)\big)$ for any $q\in[p,+\infty[\,$ and where one has
	$w \in \mc C\big([0,T];\dot B^{2-2/r}_{p,r}(\R^d)\big)\cap L^{r_0}\big([0,T];L^{p_0}(\R^d)\big)$, with in addition
$\nabla w \in L^{r_1}\big([0,T]; L^{p}(\R^d)\big)$ and $\partial_t w, \nabla^2 w\;\in L^{r}\big([0,T]; L^{p}(\R^d)\big)$.
\end{enumerate}
Furthermore, one has  $\nabla u\,\in\,L^1\bigl([0,T];L^\infty(\R^d)\bigr)$.
In particular, $u$ admits a unique flow $\psi_u$ with bounded gradient, defined as the unique solution of
\begin{equation} \label{eq:flow}
\psi_u(t,x)\,=\,x\,+\,\int_0^t u\bigl(\tau,\psi_u(\tau,x)\bigr)\,d\tau \qquad\qquad \mbox{ for all }\qquad (t,x)\in[0,T]\times\R^d\,.
\end{equation}
Finally, if we define $X_{t,\lambda}$ by the formula 
$X_{t,\lambda}(x)\,:=\,\d_{X_{0,\lambda}}\psi_u\bigl(t,\psi_u^{-1}(t,x)\bigr)$
then, for all time $t\in[0,T]$, the family  $\mc X_t:=\bigl(X_{t,\lambda}\bigr)_{1\leq\lambda\leq m}$ remains  non-degenerate  and belongs to the space $\mbb L^{\infty,p}$,
and one has the property $\rho(t)\,:=\,1+\z(t)\,\in\,\mbb L^{p}_{\mc X_t}$.

\paragraph*{\bf \tsc{Uniqueness}}
Assume that $p\in\,\big]2d\,,4d\big]$ and reinforce condition \eqref{HypForPs} for $r$ with the condition
$7/8\,<\, 1/r\, <\, 1-d/(4p) $. Define the indices $R$ and $S$ such that
\[
 \frac{1}{R}\,:=\,\frac{4}{r}\,-\,3\qquad \mbox{ and }\qquad 
 \frac{1}{S}\,:=\,\frac{1}{2}\,\left(\frac{1}{R}\,-\,\frac{1}{2}\right)\,=\,\frac{2}{r}\,-\,\frac{7}{4}\,.
\]
and assume that $w_0\,\in\,B^{2-2/R}_{p,R}$ and $\theta_0\,\in\,B^{2-2/S}_{p/2,S}$.
Then, the solution $\big(\rho,u,\theta\big)$ is unique.
\end{thm}

Some remarks about the previous statement are in order. Here, we limit ourselves to focus on the gap between the assumptions required for existence and those
imposed to get uniqueness. Specific comments on the former set of hypotheses will be formulated in Section \ref{s:a-priori}, after
the statement of Proposition \ref{PropositionEstimatesForModifiedSolution}.

\begin{rmk} \label{r:hyp-u_1}
Observe that, by definition \eqref{DefOfRs}, one has
\[
 \frac{1}{S}\,<\,\frac{2}{r}\,-\,1\,<\,\frac{1}{s}\,,
\]
so one has the embedding $B^{2-2/S}_{p/2,S}\,\hookrightarrow\,B^{2-2/s}_{p/2,s}$. The additional assumption over $\theta_0$ for having uniqueness
is then a higher regularity condition.

At the same time, we notice that $2-2/R>0$, because $1/R<1/r<1$. So, one has that $B^{2-2/R}_{p,R}\,\equiv\,L^p\,\cap\,\dot B^{2-2/R}_{p,R}$,
see also Remarks \ref{r:Besov} and \ref{r:hyp-theta} below for more details. Thus, this assumption not only contains 
a higher regularity requirement (as $R>r$), but also  a low frequency control on $w_0$.
\end{rmk}

\begin{rmk} \label{r:hyp-d}
The assumption $2d\,<\,p\,\leq\,4d $ in the uniqueness part is imposed mostly for technical reasons. It helps a lot in some computations of Section \ref{s:uniqueness}
(see Proposition \ref{p:G-integr} and its proof, carried out in Subsection \ref{ss:proof-prop-u}) and allows us to prove uniqueness by mean of
a unified and simple argument. 
\end{rmk}

Finally, as anticipated in Subsection \ref{ss:system}, let us comment on the assumption $d\geq2$.

\begin{rmk} \label{r:dim}
The one-dimensional setting is typically easier to handle (for instance, compare the statements of Theorems 1.1 and 1.2 of \cite{DFP}).
Hence, in the case $d=1$, we expect similar (or even better) results as ours to hold.

However, some of our arguments work well for $d\geq2$, whereas they would require a special treatement when $d=1$.
For the sake of conciseness, we prefer to avoid to deal with those modifications here and simply focus on the case $d\geq2$.
\end{rmk}

\subsection{Strategy of the proof} \label{ss:strategy}

In this last part of the introduction, we briefly comment on the general strategy we adopt to prove Theorem \ref{th:well-p}, and on the related difficulties.

In order to prove our main result, we resort to the same strategy of \cite{DFP}, here employed in the context of the full compressible system \eqref{FNS}.
Precisely, this is based on the following main ingredients: a reformulation of the system by use of the effective viscous flux;
maximal regularity estimates for the related equation; use of a striated regularity approach to obtain a control on the
$L^1_T(L^\infty)$ norm of the gradient $\nabla u$ of the velocity field.

However, compared with \cite{DFP}, we will face major differences in the analysis, which also represent main technical difficulties of our study.
These differences arise mostly due to the strong coupling of the momentum equation with not only
the mass conservation equation, but also the temperature equation.
As a matter of fact, the bilinear term $\mc B(\nabla u)$ appearing in the last equation of \eqref{eq:NSF}
is a source of troubles, essentially because its quadratic dependence on $\nabla u$ is very expensive in terms of time integrability properties. At the same time,
that bad term, together with $\Delta\theta$, will give a feedback (through the pressure function) also at the level of the equation for $w$ itself,
making the maximal regularity argument much more delicate for closing the estimates.
Connected with this problem, we point out that, here, we will need to use the control $\nabla u\in L^1_T(L^\infty)$
(actually, the $L^1$ condition must be upgraded to $L^m$, for some $m>1$)
not only for the proof of uniqueness, as it was the case in \cite{DFP}, but also for closing the estimates and thus establishing existence of solutions.

Finally, we comment on the uniqueness proof. As it is classically done in previous studies related to the isentropic case \eqref{eq:bar-NS}, 
uniqueness is based on $L^2$ stability estimates carried out on the Lagrangian formulation of the equations. Indeed,
the density is (almost) constant along the flow, so the use of Lagrangian coordinates avoids the problem of estimating the difference of the densities,
which would result in a loss of derivatives: owing to the hyperbolicity of the mass conservation equation,
performing stability estimates for $\rho$ would require a control on $\nabla\rho$, which is of course not compatible with the presence
of jump discontinuities in $\rho$.
Thus, also in our framework we will pass to the Lagrangian formulation and perform $L^2$ stability estimates for the difference of two
solutions; then, one closes the stability estimates by use of the Gr\"onwall lemma.
In order to do so, we will employ an original argument, consisting in making a bootstrap of the time integrability
of the solutions; in turn, this will yield that the terms appearing in the exponent of the Gr\"onwall estimate are $L^1_T$,
thus implying uniqueness in the case of same initial datum. It is precisely to perform this bootstrap of the time integrability
that the higher regularity assumptions on the initial $w_0$ and $\theta_0$ are required.

\subsection*{Outline of the paper}
Before moving on, let us give a general overview of the rest of the paper.
In Section \ref{s:tools}, we introduce the main analytical tools needed in our study. There, one also finds all the definitions needed to fully
understand the main statement, Theorem \ref{th:well-p}.
In Section \ref{s:a-priori} we perform \tsl{a priori} estimates on equations \eqref{eq:NSF}, which we previously reformulate in terms of the effective viscous
flux $w$. In Section \ref{s:existence} we complete the proof of the existence of a solution; Section \ref{s:uniqueness} is instead devoted to the Lagrangian
formulation of system \eqref{eq:NSF} and the proof of uniqueness. Finally, in Appendix \ref{app:equiv} we detail the computations leading
from equations \eqref{FNS} to equations \eqref{eq:NSF}, showing their equivalence within our functional framework.

\section{Functional and harmonic analysis toolbox} \label{s:tools}

This section is devoted to the presentation of several tools from functional and harmonic analysis which are needed in our analysis.
In Subsection \ref{ss:LP}, we recall the basics of Littlewood-Paley decomposition and Besov spaces; there, we also state some useful
interpolation inequalities and study some properties of the operator $(\Id-\Delta)^{-1}$.
Subsection \ref{ss:MR} contains classical maximal regularity estimates and their extension to the Lam\'e operator. Finally, in Subsection \ref{ss:TR}
we give an introduction to tangential regularity \tsl{\`a la Chemin} \cite{Ch1991, Ch1993, Ch1995}.

\subsection{Littlewood-Paley theory and Besov spaces} \label{ss:LP}

We  here briefly present some elements of the  Littlewood-Paley theory in $\R^d$, 
as it will come into play  for proving our main result. We refer e.g. to Chapter 2 of \cite{BahouriCheminDanchin} for more details.

We start by introducing a dyadic partition of unity with respect to the Fourier variable.
For this, fix a smooth radial function $\chi$ supported in the ball $B(0,2)$, equal to $1$ in a neighborhood of $B(0,1)$
and decreasing (actually, non-increasing) along radial directions. Set
$\varphi\left(\xi\right)=\chi\left(\xi\right)-\chi\left(2\xi\right)$ and
$\vphi_j(\xi):=\vphi(2^{-j}\xi)$ for all $j\geq0$.
The non-homogeneous dyadic blocks $(\Delta_j)_{j\in\Z}$ are defined by\footnote{Here and throughout this work, we agree  that  $f(D)$ stands for 
the pseudo-differential operator operator associated to the function $f(\xi)$, that is the operator (formally) defined by $u\mapsto\mc{F}^{-1}(f\,\mc{F}u)$.} 
$$
\Delta_j\,:=\,0\quad\mbox{ if }\; j\leq-2,\qquad\Delta_{-1}\,:=\,\chi(D)\qquad\mbox{ and }\qquad
\Delta_j\,:=\,\varphi(2^{-j}D)\quad \mbox{ if }\;  j\geq0\,.
$$
We  also introduce the following low frequency cut-off operators:
\begin{equation} \label{eq:S_j}
S_j\,:=\,\chi(2^{-j}D)\,=\,\sum_{k\leq j-1}\Delta_{k}\quad \mbox{ for }\quad j\geq0,\qquad\qquad
 S_j=0\ \mbox{ for }\ j<0.
\end{equation}
It is well known that for any $u\in\mc{S}'(\R^d)$,  one has the equality 
$$
u\,=\,\sum_{j\geq-1}\Delta_ju\qquad\qquad \mbox{ in the sense of }\qquad \mc{S}'(\R^d)\,.
$$
This is the Littlewood-Paley decomposition of the tempered distribution $u$.

In some parts of our analysis, see especially maximal regularity estimates from Subsection \ref{ss:MR}, we need to use \emph{homogeneous}
spectral cut-off operators. To this end, we define the operators $\dot\Delta_j$ and $\dot S_j$ by the formulas
\begin{equation} \label{eq:S_j-dot}
\dot\Delta_j\,:=\,\varphi(2^{-j}D)\qquad \mbox{ and }\qquad \dot S_j\,=\,\chi(2^{-j}D)\qquad\qquad \mbox{ for all }\quad  j\in\Z\,.
\end{equation}
Note  that the homogeneous decomposition
\begin{equation}
\label{eq:LPh} u\,=\,\sum_{j\in\Z}\dot\Delta_ju
\end{equation}
holds true only modulo  polynomials, which makes formula \eqref{eq:LPh} unwieldy. A way to have equality in \eqref{eq:LPh} in the sense of tempered distributions 
is to restrict the attention to elements of  the set $\mc S_h'(\R^d)$ of tempered distributions $u$ such that
\begin{equation} \label{def:S_h}
\lim_{j\ra-\infty}\dot S_ju\,=\,0\qquad\qquad \mbox{ in }\qquad \mc S'(\R^d)\,.
\end{equation}
We refer to \cite{Cobb} and references therein for more insights on the spaces $\mc S'_h(\R^d)$.

After this preparation, we can now introduce the definition of non-homogeneous and homogeneous Besov spaces.
\begin{definition} \label{d:B}
Let $s\in\R$ and $1\leq p,r\leq+\infty$.
\begin{enumerate}[(i)]
\item The \emph{non-homogeneous Besov space}
$B^{s}_{p,r}$ is the  set of tempered distributions $u$ for which
$$
\|u\|_{B^{s}_{p,r}}\,:=\,\left\|\left(2^{js}\,\|\Delta_ju\|_{L^p}\right)_{j\geq-1}\right\|_{\ell^r}\,<\,\infty\,.
$$
\item The \emph{homogeneous Besov space} $\dot B^s_{p,r}$ is the subset of homogeneous tempered distributions $u\in \mc S'_h(\R^d)$ such that
$$
\|u\|_{\dot B^{s}_{p,r}}\,:=\,\left\|\left(2^{js}\,\|\dot\Delta_ju\|_{L^p}\right)_{j\in\Z}\right\|_{\ell^r}\,<\,\infty\,.
$$
\end{enumerate}
\end{definition}

It is well known that, for any $s\in\R$, one has $B^s_{2,2}\equiv H^s$ with equivalence of norms,
and that non-homogeneous (respectively, homogeneous) Besov spaces are interpolation spaces between Sobolev spaces $W^{k,p}$ (respectively, $\dot W^{k,p}$). 
Furthermore, for all $p\in\,]1,+\infty[\,$, one has the following continuous embeddings (see the proof in \cite[Chap. 2]{BahouriCheminDanchin}): 
$$
\dot B^0_{p,\min(p,2)}\hookrightarrow L^p\hookrightarrow \dot B^0_{p,\max(p,2)}\qquad \mbox{ and }\qquad 
B^0_{p,\min(p,2)}\hookrightarrow L^p\hookrightarrow B^0_{p,\max(p,2)}\,.
$$

Next, we point out a useful fact about the relation between non-homogeneous and homogeneous Besov spaces.
\begin{rmk} \label{r:Besov}
When $s>0$, it is possible to see that one has $B^s_{p,r}\,\equiv\,L^p\,\cap\,\dot B^{s}_{p,r}$, with equivalence of norms.
\end{rmk}

We shall often use also the embedding properties that are stated in the following proposition. 
\begin{prop}\label{p:embed}
Let $1\leq p_1\leq p_2\leq+\infty$. 
The space $B^{s_1}_{p_1,r_1}(\R^d)$ is continuously embedded into the space $B^{s_2}_{p_2,r_2}(\R^d)$ if
$$
s_2\,<\,s_1-d\left(\frac{1}{p_1}-\frac{1}{p_2}\right)\qquad\mbox{ or }\qquad
s_2\,=\,s_1-d\left(\frac{1}{p_1}-\frac{1}{p_2}\right)\;\;\mbox{ and }\;\;r_1\,\leq\,r_2\,. 
$$
The space $\dot B^{s_1}_{p_1,r_1}(\R^d)$ is continuously embedded in the space $\dot B^{s_2}_{p_2,r_2}(\R^d)$ if
$$
s_2\,=\,s_1-d\left(\frac{1}{p_1}-\frac{1}{p_2}\right)\;\;\mbox{ and }\;\;r_1\,\leq\,r_2\,. 
$$
\end{prop}

Next, we make use of Littlewood-Paley decomposition and Besov spaces to derive some useful interpolation inequalities,
of Gagliardo-Nirenberg type that we will use regularly in the proof of existence.
Their proof relies on the use of Littlewood-Paley decomposition and follows a standard strategy: cutting functions into low and high frequencies and, then,
estimating the low frequencies by means of the low regularity norm and the high frequencies by use of the high regularity norm.
For conciseness, we omit the details here.

\begin{lemma} \label{l:GN-ineq}
For any $p\in[1,+\infty]$, there exists a universal constant $C>0$ such that
\[
 \left\Vert\nabla f \right\Vert_{L^{p}}\,\le\, C\, \left\Vert f \right\Vert_{L^{p}}^{1/2}\,\left\Vert \Delta f \right\Vert_{L^{p}}^{1/2}\,.
\]
In addition, for any $p\in\,]d,+\infty[\,$, one has the inequality
\[
 \left\Vert\nabla f \right\Vert_{L^{p}}\,\le\, C\, \left\Vert f \right\Vert_{L^{p/2}}^{1-\alpha}\, \left\Vert \Delta f \right\Vert_{L^{p/2}}^{\alpha}\,,
\qquad\qquad \mbox{ with }\qquad \alpha\, =\,\dfrac12\,+\, \dfrac d{2p}\,.
\]
\end{lemma}

\begin{lemma}\label{LemmaUsefulInequalityGeneralVersion}
Let the triplet $(q,q_0,q_1)\in[1,+\infty]$ be such that $\max\big\{q_0,q_1\big\}\leq q$ and that
$1/q_1\,-\,1/q\,\leq\,1/d$. Define the number $\alpha\in\,]0,1[\,$ by the formula
\[
 \alpha\,:=\,\frac{1+d/q- d/q_1}{1+d/q_0-d/q_1}\,. 
\]

Then, there exists a universal constant $C>0$ such that
	\begin{equation}\label{LemmaUsefulInequalityGeneralVersioneq1}
\Vert f \Vert_{L^q}\,\le\, \Vert f \Vert_{L^{q_0}}^{\alpha}\,\Vert \nabla f \Vert_{L^{q_1}}^{1-\alpha}\,.
	\end{equation}
\end{lemma}

We conclude this part by studying some continuity properties of the operator $(\Id-\Delta)^{-1}$,
which will be useful in estimating the functional norms of the vector field $v$ introduced in Theorem \ref{th:well-p}.
We start by recalling Lemma 2.3 and Corollary 2.4 of \cite{DFP}.

\begin{lemma}\label{l:D-1}
For any $p\in[1,+\infty]$, the operator $\Delta(\Id-\Delta)^{-1}$ maps continuously $L^p$ into itself. In addition, there
exists a constant $C>0$, depending only on $d$, such that, for any $p\in[1,+\infty]$ and any $f\in L^p$, one has
$$
\|\Delta(\Id-\Delta)^{-1}f\|_{L^p}\,\leq\, C\,\|f\|_{L^p}\,.
$$
\end{lemma}

\begin{cor} \label{c:Delta-Id}
Let $u$ solve the elliptic equation $(\Id-\Delta)u\,=\,f$ on $\R^d$, with $f\in L^p$ for some $p\in[1,+\infty]$.
Then one has $u\in W^{1,p}$, with $\Delta u\in L^p$, together with the estimate
$$
\|u\|_{W^{1,p}}\,+\,\|\Delta u\|_{L^p}\,\leq\,C\,\|f\|_{L^p}\,,
$$
for some suitable constant $C>0$ depending just on $d$.

In the case  $1<p<+\infty$, one further deduces  $u\in W^{2,p}$, with the estimate
$$
\|u\|_{W^{2,p}}\,\leq\,C\,\|f\|_{L^p}\,.
$$
\end{cor}

Combining the above Corollary \ref{c:Delta-Id} with Lemma \ref{LemmaUsefulInequalityGeneralVersion}, we infer that,
for any $q\in\,]1,+\infty[\,$ fixed, given some $f\in L^q$, one has
$\nabla(\Id-\Delta)^{-1}f\in W^{1,q}$. By Sobolev embeddings, this implies that
$\nabla(\Id-\Delta)^{-1}f\in L^m$, with $1/m\,:=\,1/q\,-\,1/d$ if $q<d$, any $m\in[q,+\infty[\,$ if $q=d$ and
$m=+\infty$ in the case $q>d$.
Keeping $m=p$ fixed and varying the index $q$, we deduce the following property, which will turn out to be useful in the computations of Section \ref{s:a-priori}.

\begin{cor}\label{ClaimEstimationForOperatorPGeneralp}
Let $d\geq2$ and $p\in\,\left]\frac{d}{(d-1)}\,,\,+\infty\right[\,$. Define $\overline p\,:=\,\frac {p\,d}{d+p}\,\in\,\left]1,\min\big\{p,d\big\}\right[\,$.
There exists a universal multiplicative constant
$C>0$ such that, for any $q\in[\oline p,p]$ and any $f\in L^q$, one has the following estimate:
	\begin{equation}\label{EstimationForOperatorPGeneralp}
		\left\Vert\nabla \left( \Id-\Delta\right)^{-1}f\right\Vert _{L^p}\,\le\, C\,\Vert f\Vert_{L^{q}}\,.
	\end{equation}
\end{cor}

\subsection{Maximal regularity} \label{ss:MR}

In this section, we revisit some known results about maximal regularity estimates for the heat kernel and their extension to treat the Lam\'e operator $\mc L$.
We refer to Section 2B of \cite{DFP} and references therein for the proofs of these estimates. However, we recall here most of the statements, in order to
clarify the numerology of Theorem \ref{th:well-p} and to fix the notation to be used in the computations of Sections \ref{s:existence} and \ref{s:uniqueness}.

\subsubsection{Maximal regularity for the heat kernel}
Maximal regularity estimates for the heat operator $\d_t-\Delta$ are derived from writing the Duhamel formula associated to the heat equation,
differentiating twice with respect to the space variable and finally estimating the initial datum and the forcing term separately.

Concerning the initial datum, one simply uses the characterisation of the homogeneous Besov spaces through the heat kernel,
as given in the following statement (see Theorem 2.34 of \cite{BahouriCheminDanchin} for details and a proof).

\begin{prop}\label{PropositionBahouri}
Let $s>0$ and $(p,r)\in[1,+\infty]^2$. There exists a constant $C>0$ such that
\[
C^{-1}\,\Vert f \Vert _{\dot B^{-s}_{p,r}}\,\le\,\left\Vert \Vert t^{s/2}e^{t\Delta}f\Vert_{L^p}\right\Vert_{L^r\left( \R_+; dt/t\right)}\,\le\,
C\, \Vert f\Vert_{\dot B^{-s}_{p,r}}\,.
\]
\end{prop}

Because of the non-linearity appearing in the momentum equation in system \eqref{eq:NSF}, it is convenient to estimate not only the second order derivatives,
but also the first order ones and the term of zero order. Here below, we use the couple of indices $\big(r_j,p_j\big)$ related to $\nabla^jf$,
where $j\in\big\{0,1,2\big\}$, with the convention that $r_j$ is related to time integrability and $p_j$ denotes the integrability with respect to the space variable.
In addition, for $j\in\big\{0,1,2\big\}$, we define the operator
\[
\mathcal A_j\,:\; f\,\mapsto\, \int_0^t \nabla^j e^{(t-\tau)\Delta}f(\tau, \cdot)\,d\tau\,. 
\]

The following conditions of integrability indices with respect to time and space will then be well-suited
for treating both the initial datum and the forcing term:
for finite integrability indices $1<p_2\leq p_1$ and $p_0\in[p_2,+\infty]$, we choose
\begin{align}
\label{MaximalRegularityCondition1}
&	\frac{2}{r_2}\,+\,\frac{d}{p_2}\,-\,\frac{d}{p_1}\,>\,1\,, \qquad\qquad 
	\frac{2}{r_2}\,+\, \frac{d}{p_2}\,=\, 1\,+\, \frac2{r_1}\,+\, \frac d{p_1}\,, \\
	\label{MaximalRegularityCondition3}
&\qquad	2\,-\,\frac{2}{r_2}\,-\, \frac{d}{p_1}\,<\,0\,, \qquad\qquad
	\frac{2}{r_2}\,+\,\frac{d}{p_2}\,=\, 2\,+\,\frac{2}{r_0}\,+\, \frac d{p_0}\,.
\end{align}	

The next lemma treats the case $j=2$, corresponding to classical maximal regularity results. We refer to Lemma 7.3 of \cite{LemR} for details.

\begin{lemma}\label{LemmaMaximalRegularitySecondDerivative}
Let $T>0$ and $1<p_2, r_2<+\infty$.
Then $\mathcal A_2$ is a bounded operator from $L^{r_2}\big(\,]0,T[\,; L^{p_2}(\R^d)\big)$ into itself. Moreover, one has the estimate
\[
	\Vert \mathcal A_2 f \Vert_{L_T^{r_2}(L^{p_2})}\,\le\, C\, \Vert f \Vert_{L_T^{r_2}(L^{p_2})}\,.
\]
\end{lemma}

The next results, related to the cases $j=1$ and $j=0$, correspond respectively to Lemma 2.7 and Lemma 2.8 of \cite{DFP}.

\begin{lemma}\label{LemmaMaximalRegularityFirstDerivative}
Let $T>0$ and $1<p_2, r_2<+\infty$. Fix Lebesgue exponents $p_1$ and $r_1$ satisfying the condition $0\leq 1/p_2-1/p_1<1/d$ together with \eqref{MaximalRegularityCondition1}.

Then $\mathcal A_1$ is a bounded operator from $L^{r_2}\big(\,]0,T[\,; L^{p_2}(\R^d)\big)$ into $L^{r_1}\big(\,]0,T[\,; L^{p_1}(\R^d)\big)$.
In addition, the following estimates hold true:
\[
	\Vert \mathcal A_1 f \Vert_{L_T^{r_1}(L^{p_1})}\,\le\, C\, \Vert f \Vert_{L_T^{r_2}(L^{p_2})}\,.
\]

\end{lemma}	

\begin{lemma}\label{LemmaMaximalRegularityZeroDerivative}
Let $T>0$ be a given time and $1<p_2, r_2<+\infty$. Take Lebesgue exponents $p_0$ and $r_0$ such that the conditions $0\leq 1/p_2-1/p_0<2/d$
and \eqref{MaximalRegularityCondition3} are satisfied. 

Then $\mathcal A_0$ is a linear bounded operator from the space $L^{r_2}\big(\, ]0,T[\,; L^{p_2}(\R^d)\big)$ into the space
$L^{r_0}\big( \,]0,T[\,; L^{p_0}(\R^d)\big)\cap \mathcal C\big(\, ]0, T[\,; \dot B^{2-2/r_2}_{p_2, r_2}\big)$. In addition,
one has the estimate
\[
\Vert \mathcal A_0 f \Vert_{L_T^{\infty}\left(\dot B^{2-2/r_2}_{p_2, r_2}\right)}\,+\, 
	\Vert \mathcal A_0 f \Vert_{L_T^{r_0}(L^{p_0})}\,\le\, C\, \Vert f \Vert_{L_T^{r_2}(L^{p_2})}\,.
\]
\end{lemma}	

\subsubsection{Maximal Regularity for $\mathcal L$}
The previous results can be applied to deduce maximal regularity estimates for the Lam\'e operator
$\mc L$ defined in \eqref{DefinitionL}.
More precisely, consider the equation 
\begin{equation}\label{EllipticEquation}
	\begin{cases}
	\partial_t U\, +\, \mathcal L U\, =\,F
	\\[1ex]
	U_{|t=0}\,=\,U_0\,. 
	\end{cases}	
\end{equation}
Using the Leray-Helmholtz decomposition $\Id\,=\,\P\,+\,\Q$, where $\P$ denotes the $L^2$-projector onto the divergence-free
vector fields and $\Q$ its $L^2$-orthogonal operator, one can prove the following result
(see Proposition 2.9 of \cite{DFP} for its proof).

\begin{lemma}\label{LemmaEstimationSolutionElllipticEquation}
Let the set of indices $(p_j, r_j)_{j=0,1,2}$ satisfy 
\[
	1<p_2, r_2<+\infty, \qquad r_2<r_0, \qquad r_2<r_1, \qquad p_0\ge p_2, \qquad p_1\ge p_2, 	
\]
and relations \eqref{MaximalRegularityCondition1} and \eqref{MaximalRegularityCondition3}. Let $U_0 \in \dot B_{p_2, r_2}^{s_2}$, with $s_2:= 2-2/r_2$,
and $F \in L_{\loc}^{r_2}\big( \R_+; L^{p_2}( \R^d)\big)$.

Then, for all $T>0$, system \eqref{EllipticEquation} has a unique solution
$U \in \mathcal C\left( [0,T]; \dot B^{s_2}_{p_2, r_2}(\R^d)\right)\cap L^{r_0}\left( [0,T]; L^{p_0}(\R^d)\right)$, with $\nabla U\in L^{r_1}\left([0,T];L^{p_1}(\R^d)\right)$ and
$\partial_t U$, $\nabla^2U\; \in L^{r_2}\left( [0,T]; L^{p_2}(\R^d)\right)$.
Moreover, there exists a constant $C>0$ (depending on $\mu$, $\lambda$, $d$, $p_0$, $p_1$, $p_2$ and $r_2$) such that the following estimate holds: 
\begin{align}\label{LemmaEstimationSolutionElllipticEquation_Eq1}
	&\Vert U \Vert_{L_T^\infty\left( \dot B_{p_2, s_2}^{r_2}\right)} + \Vert U \Vert_{L_T^{r_0}\left(L^{p_0}\right)} + \Vert \nabla U \Vert_{L_T^{r_1}\left(L^{p_1}\right)}	\\ 
\nonumber	& 
	\qquad\qquad\qquad  + \Vert \left( \partial_t U, \nabla^2 U\right)\Vert_{L_T^{r_2}\left(L^{p_2}\right)}\, 
	\le\, C\left( \Vert U_0 \Vert_{\dot B^{s_2}_{p_2, r_2}} + \Vert F\Vert_{L_T^{r_2}\left( L^{p_2}\right)}\right)\,. 
\end{align}	
\end{lemma}	

\subsection{Introduction to striated regularity} \label{ss:TR}
Now, we recall some definition and basic properties related to striated (also dubbed tangential) regularity,
in the spirit of Chemin's works \cite{Ch1991, Ch1993, Ch1995}.
We will mainly follow the presentation and notation of \cite{DFP}.

Denote by $C_b(\R^d)$ the set of continuous and bounded functions over $\R^d$. We use the same notation also for vector fields.
For $p\in\,]d,+\infty[\,$, define the set of vector fields
\begin{equation}\label{d:SetVectorField}
\mathbb L^{\infty, p}\,:=\, \left\{Y \in C_b(\R^d) \; \Big| \quad \Vert \nabla Y \Vert_{L^p}\,<\,+\,\infty \right\}\,,
\end{equation}
endowed with the norm
\begin{equation}\label{DefLInftyPNorm}
\Vert Y \Vert_{\mathbb L^{\infty, p}}\,:=\, \Vert Y\Vert_{L^\infty}\, +\, \Vert \nabla Y\Vert_{L^p}	\,.
\end{equation}

We now present the definition of \emph{tangential regularity} of functions along a vector field $Y$. In order to enlarge the class of functions
for which the definition makes sense, we will define it in the weak sense.
\begin{definition}\label{d:DDField}
Let $p\in\,]d,+\infty[\,$ and $Y \in \mathbb L^{\infty, p}$ be a given vector field.
Given $f\in L^{p'}_\loc$, with $1/p\,+\,1/p'\,=\,1$, we define
\[
\partial_Y f\, :=\,\div\big(f\,Y\big)\,-\,f\,\div Y\,. 
\]
\end{definition}

Observe that, if in addition $f\in W^{1,1}_\loc$, then this quantity coincides with the usual definition
\[
 \d_Yf\,\equiv\,Y\cdot\nabla f\,=\,\sum_{j=1}^d Y^j \partial_j f\,.
\]
It goes without saying that the definition above extends to vector-valued and matrix-valued functions $f$, by simply acting componentwise.

Next, we consider a family $\mathcal X= \left( X_{\lambda}\right)_{1\le \lambda\le m}$ of vector fields.
We introduce the following convenient definition.

\begin{definition}\label{d:NonDeg}
Given a set $\mathcal X= \left( X_{\lambda}\right)_{1\le \lambda\le m}$ of vector fields belonging to $\mathbb L^{\infty, p}$,
we say that $\mc X$ is a \emph{non-degenerate family} of vector fields if 
\[
I(\mathcal X)\,:=\, \inf_{x \in \R^d}\,\sup_{\Lambda \in \Lambda^m_{d-1}}\left\vert \overset{d-1}{\wedge}X_\Lambda(x)\right\vert^{1/(d-1)}\;>\,0\,, 	
\]
where $\Lambda \in \Lambda_{d-1}^m$ means $\Lambda= \left( \lambda_1, \dots, \lambda_{d-1}\right)$, $\lambda_i\in \left\{1, \dots, m\right\}$,
with $\lambda_i \le \lambda_j$ for $i <j$, and $\overset{d-1}{\wedge}X_\Lambda$ is the unique element of $\R^d$ such that 
\[
\forall Y \in \R^d, \qquad \left( \overset{d-1}\wedge X_\Lambda\right)\cdot Y\, =\, \det\left( X_{\lambda_1}, \dots, X_{\lambda_{d-1}}, Y\right)\,.
\]
For such a family, we also define
\[
||| \mathcal X |||_{\mathbb L^{\infty, p}}\,:=\, \sup_{\lambda \in \Lambda}\Vert X_{\lambda}\Vert_{\mathbb L^{\infty, p}}\,.
\]
More in general, for any normed space $E$, we shall denote 
\[
||| \mathcal X|||_{E}\,:=\, \sup_{\lambda \in \Lambda}\Vert X_\lambda\Vert_{E}\,. 
\]
\end{definition}

We are now ready to ``measure'' the regularity of some function $f$ along given directions.

\begin{definition}\label{d:LpdirecX}
Let $p\in\,]d,+\infty[\,$ and $Y \in \mathbb L^{\infty, p}$. For any $f \in L^\infty$, we write $f \in \mathbb L^p_Y$ if
$\divergence \left( fY\right)\in L^p$, or equivalently if $\d_Yf\in L^p$.
Moreover, if $\mathcal X= \left(X_{\lambda}\right)_{1\le \lambda\le m}$ is a non-degenerate family of vector fields in $\mathbb L^{\infty, p}$, then we set
	\[
	\mathbb L_{\mathcal X}^p\,:=\, \bigcap\limits_{1\le \lambda\le m}\mathbb L_{X_{\lambda}}^p\,,
\quad \mbox{ with }\quad  \Vert f \Vert_{\mathbb L_{\mathcal X}^p}\,:=\,\frac{1}{I(\mathcal X)}\,
\Big( \Vert f \Vert_{L^\infty} ||| \mathcal X |||_{\mathbb L^{\infty, p}} + ||| \divergence(f\mathcal X)|||_{L^p}\Big)\,. 
	\]
\end{definition}

\section{A priori estimates} \label{s:a-priori}

This section collects suitable \tsl{a priori} estimates for smooth solutions to the Navier-Stokes-Fourier system, necessary for our proof of
existence and uniqueness. Therefore, we are going to assume to dispose of a smooth solution $\big(\rho,u,\theta\big)$ to system
\eqref{eq:NSF}-\eqref{def:B_u}, related to the given initial datum $\big(\rho_0,u_0,\theta_0\big)$,
defined on the whole $\R_+\times\R^d$ and decaying sufficiently fast at infinity. We perform formal computations.

The \tsl{a priori} estimates derived here are quite involved and will be performed in several steps. First of all, in Subsection \ref{ss:reform}
we reformulate the momentum equation as a parabolic equation for the effective viscous flux $w$, coupled with the temperature equation for $\theta$
and a continuity equation for the densityh perturbation function $\z\,:=\,\rho-1$. In Section \ref{ss:apriori-stat} we state an existence
result for the reformulated system, which immediately implies the existence statement of Theorem \ref{th:well-p}, and explain
the general strategy of the proof. The proof itself, only pertaining to the derivation of the necessary \tsl{a priori} estimates
(see Section \ref{s:existence} for the actual proof of existence),
is carried out in Subsections \ref{ss:static}, \ref{SubsectionAPrioriEstimatesForDensity},
\ref{ss:apriori-temp}, \ref{SubsectionMaximalRegularityForVelocity}, \ref{SecTanReg} and \ref{ss:close}.

\subsection{Reformulation of the system} \label{ss:reform}

In order to handle discontinuities in the density function, the pressure is obviously an obstacle. Inspired by
Serre \cite{Serre} and Hoff \cite{Hoff_JDE, Hoff_1997}, we introduce the new variable
\begin{equation}\label{Defv}
	v\,:=\, -\,\frac{1}{\nu}\,\nabla ( \Id - \Delta)^{-1} P(\rho, \theta)\,,\qquad\qquad \nu\,:=\,\lambda\,+\,2\,\mu\,,
\end{equation}
where we follow the idea of \cite{DFP} and use the operator $(\Id-\Delta)^{-1}$ instead of $(-\Delta)^{-1}$, as the former is less singular for low frequencies.
For simplicity of notation, and without loss of any generality, from now on let us set
\[
\nu\,=\ \lambda\,+\,2\,\mu\,=\,1\,.
\]
We also define the \emph{effective viscous flux}, namely the modified velocity field 
\begin{equation}\label{def:EVF}
	w\,:=\, u\,-\,v\, =\, u \,+\,\nabla (\Id- \Delta )^{-1} P(\rho, \theta)\,. 
\end{equation}
Finally, we also introduce the density variation function
\[
	\zeta\,:=\, \rho\,-\,1\,,
\]
measuring the departures of the density from the constant value $1$.

Of course,  we also define the initial data for $\z$ and $w$ according to the above formulas:
\[
 w_{|t=0}\,=\,w_0\,:=\,u_0\,-\,v_0 \qquad\qquad \mbox{ and }\qquad\qquad \zeta_{|t=0}\,=\,\zeta_0\,:=\,\rho_0\,-\,1\,,
\]
where we have set $v_0\,=\,-\,\nabla ( \Id - \Delta)^{-1} P(\rho_0, \theta_0)$.

\medbreak
We are now going to reformulate system \eqref{eq:NSF} in terms of the new set of unknowns $\big(\z,w,\theta\big)$.
For this, we will use the decomposition $\rho\,=\,\zeta + 1$ and $u\,=\, w+v$ everywhere; however, for simplicity of writing, we will keep writing
$u$ and $\rho$ when convenient.

We first focus on the density equation: simple computations yield the transport equation
\begin{align}\label{EquationForZeta}
	\partial_t \zeta\, +\, u \cdot \nabla \zeta\,=&\,-\,(\z+1)\,\div u \\
\nonumber
&=\, -\,(\z+1) \, \divergence w\, +\, (\z+1)\, \Delta (\Id-\Delta)^{-1}P(\rho, \theta)\,.
\end{align}
For later use, we remark that, from system \eqref{eq:NSF}, we can compute an equation for the pressure function $P(\rho,\theta)$.
Indeed, combining the first and third equations in \eqref{eq:NSF}, we find
\begin{equation}\label{EquationPressure}
	\partial_t P\,=\, -\, \divergence\big(P\,u\big)\,+ \,
	 \Big( g_0(\rho) + \theta g_1(\rho)\Big)\,\divergence u\, +\, \frac{\pi_1(\rho)}{\rho}\,G
\end{equation}	
where we have defined the functions
\[
 g_0(\rho)\,:=\, \pi_0(\rho)\,-\,\pi_0'(\rho)\,\rho\qquad\qquad \mbox{ and }\qquad\qquad
g_1(\rho)\,:=\, \pi_1(\rho)\,-\,\pi_1'(\rho)\,\rho
\] 
and, recalling the bilinear operator $\mc B(\nabla u)$ from \eqref{def:B_u}, we have defined the forcing term
\begin{equation}\label{defintionF_1}
	G\,:=\, \mc B(\nabla u)\, - \, \theta\, \pi_1(\rho)\, \divergence u\, +\,\k\,\Delta \theta\,.	
\end{equation}

Finally, to close the system, we need an equation for the effective viscous flux $w$. This is derived from the velocity equation in \eqref{eq:NSF},
by using the expansion $u\,=\,w\,+\,v$. Observe that definition \eqref{Defv} of $v$ (where we have fixed the value $\nu=1$) directly yields that
\[
 \mc Lv\,+\,\nabla P\,=\,\nabla(\Id-\Delta)^{-1}P\,,
\]
which is more regular than $P$, owing to the presence of an operator of order $-1$ acting on it.
In addition, we can use the mass conservation balance, namely the equation for $\rho$, to write
\[
 \d_t\big(\rho\,u\big)\,+\,\div\big(\rho\,u\otimes u\big)\,=\,\rho\,\d_tu\,+\,\rho\,(u\cdot\nabla)u\,.
\]
Thus, also using the pressure equation \eqref{EquationPressure}, we obtain the sought equation for $w$:
\begin{equation}\label{EquationForw}
	\rho\,\partial_t w\, +\, \mathcal L w\,=\,- \rho\, F\,,
\end{equation}
where the source term $F$ is given by 
\begin{align}
\label{DefinitionF}	F\,&:=\,	
	(w\cdot \nabla) w\, +\, (w\cdot \nabla) v\, +\, (v\cdot \nabla) w\,+\,
    (v\cdot \nabla) v \,+\,\frac{1}{\rho}\,\nabla\left( \Id - \Delta\right)^{-1} P 
	\\ 
	\nonumber	&\qquad\qquad 
	+\, \nabla \left( \Id - \Delta\right)^{-1}\left( 
	-\,\divergence \big( P\,u\big)\, +\,  \Big( g_0(\rho)\,+\, \theta\,g_1(\rho)\Big)\,\divergence u\, +\, \frac{\pi_1(\rho)}{\rho}\,G\right)\,. 
\end{align}

All in all, by considering the variables $(\zeta, w, \theta)$ as the new set of unknowns, we obtain the following system of equations,
where the terms $F$ and $\mc B(\nabla u)$ have been defined, respectively, in \eqref{DefinitionF} (see also \eqref{defintionF_1} for the definition of $G$
appearing therein) and \eqref{def:B_u}:
\begin{equation}\label{FNSwithTemperatureModified}
	\left\{ 
		\begin{aligned}
		&\partial_t \zeta\, +\, u \cdot \nabla \zeta\,=\, -\,(\z+1) \, \divergence w\, +\, (\z+1)\, \Delta (\Id-\Delta)^{-1}P(\rho, \theta)
		\\
		&\rho\, \partial_t w\, +\, \mathcal L w\, =\,-\, \rho\, F
		\\
		&\rho\, \partial_t \theta\, +\, \rho\, u\cdot \nabla \theta\, +\, \divergence u \, \pi_1(\rho)\,\theta\, -\,\k\,\Delta \theta
		\,=\, \mc B(\nabla u)\,,
	\end{aligned}	
	\right.
	\end{equation}	
complemented with the initial datum
\[
\big(\zeta, w, \theta\big)_{|t=0}\,=\, \big(\zeta_0, w_0, \theta_0\big)\,.
\]

\subsection{Main statement and strategy of the proof} \label{ss:apriori-stat}

The idea is to take advantage of these new variables to prove existence of $\big(\rho, u, \theta\big)$.
More precisely, we will show that there exists $\big(\zeta, w, \theta\big)$ which solves \eqref{FNSwithTemperatureModified} in the weak sense and for which the density
perturbation function $\zeta$ remains small (in the $L^\infty$ norm), at least for short times.
Then, we will eventually infer appropiate bounds for the original set of unknowns $\big(\rho, u, \theta\big)$.

The next statement is the main result of this section, and is the core of the existence proof in Theorem \ref{th:well-p}.

\begin{prop}\label{PropositionEstimatesForModifiedSolution}
Let $d\geq2$.
Let $p$, $r$ and $s$ in $\,]1,+\infty[\,$ and $p_0\in\,]1,+\infty]$ satisfy conditions \eqref{HypForPs} and \eqref{HypForS}. 
Define the indices $r_0$ and $r_1$ by \eqref{DefOfRs}. 

There exists $\veps_0>0$ such that the following property holds true.
Take initial data $\zeta_0 \in  L^p\cap L^\infty$, $w_0 \in \dot B^{2-2/r}_{p,r}$ and $\theta_0\in B^{2-2/s}_{p/2,s}\cap L^p$.
Assume that, for some $\veps\in\,]0,\veps_0]$, one has
\begin{equation}\label{HypothesisForZetaPropositionEstimatesForModifiedSolution}
\Vert \zeta_0 \Vert_{L^\infty}\,\leq\,\veps\,.	
\end{equation}
Finally, assume that  there exists  a non-degenerate family   $\mc X_0\,=\,\bigl(X_{0,\la}\bigr)_{1\leq\la\leq m}$ of vector-fields in
$\mbb L^{\infty,p}$ such that $\rho_0$ belongs to  $\mbb L^{p}_{\mc X_0}$.

Then, there exists a time $T>0$ and a weak solution $\big(\zeta, w, \theta\big)$ to system \eqref{FNSwithTemperatureModified} in the time interval $[0, T]$, such that 
\begin{equation}\label{TesisLInftyZeta}
\Vert \zeta \Vert_{L_T^\infty\left( L^\infty\right)}\,\leq\,4\,\varepsilon\,
\end{equation}
and such that the following regularity properties are satisfied:
\begin{enumerate}
	\item $w \in \mc C_T \left(\dot B^{2-2/r}_{p,r}\right)\cap L_T^{r_0}\left( L^{p_0}\right)$, with $\nabla w \in L_T^{r_1}\left( L^{p}\right)$ and
$\partial_t w, \nabla^2 w\;\in L_T^{r}\left( L^{p}\right)$;
	\item $\theta \in \mc C_T\left( B^{2-2/s}_{p/2,s}\cap L^{p}\right)$, with $\partial_t \theta, \nabla^2\theta\; \in L_T^s\left( L^{p/2}\right)$;
	\item $\zeta\in L_T^{\infty}\left( L^{p}\cap L^\infty\right)$, with $\z\in \mc C_T(L^q)$ for any $q\in[p,+\infty[\,$.
\end{enumerate}
Furthermore, after setting $u\,:=\,w\,+\,v$, with $v$ defined in \eqref{Defv}, one has  $\nabla u\,\in\,L^1_T\bigl(L^\infty\bigr)$,
thus $u$ admits a flow $\psi_u$ with bounded gradient, as defined in \eqref{eq:flow}.
Finally, if we define $X_{t,\lambda}(x)\,:=\,\d_{X_{0,\lambda}}\psi_u\bigl(t,\psi_u^{-1}(t,x)\bigr)$,
then, for all time $t\in[0,T]$, the family  $\mc X_t:=\bigl(X_{t,\lambda}\bigr)_{1\leq\lambda\leq m}$ remains  non-degenerate  and belongs to the space $\mbb L^{\infty,p}$,
and one has the property $\rho(t)\,:=\,1+\z(t)\,\in\,\mbb L^{p}_{\mc X_t}$.
\end{prop}

Some remarks on the previous assumptions are in order. In particular, we comment on the assumptions on the integrability indices which appear
in \eqref{HypForPs} and \eqref{HypForS}.

\begin{rmk} \label{r:s}
Since $p>d$ according to \eqref{HypForPs}, one has $1\,+\,d/p\,<\,2$.
This immediately yields the inequality
\[
 \frac{2}{r}\,-\,1\,<\,\frac{1}{r}\,. 
\]
Therefore, it is possible to choose the index $s\in\,]r,+\infty[\,$ such that the first relation in \eqref{HypForS} is fulfilled. In addition, using
that $1/s<1/r$ and performing simple computations, we deduce that, for the second condition in \eqref{HypForS} to be satisfied, it is enough that
$r$ verifies
\[
 \frac{1}{r}\,<\,\left(1\,+\,\frac{d}{3p}\right)^{-1}\,.
\]
We point out that, under our assumption $d<p$, one has
\[
 1\,-\,\frac{d}{2p}\,<\,\left(1\,+\,\frac{d}{3p}\right)^{-1}\,<\, 1\,-\,\frac{d}{4p}\,.
\]
Hence, the second condition in \eqref{HypForS} would be fulfilled for free, should we strengthen the requirements appearing in \eqref{HypForPs} for $r$.
Both conditions in \eqref{HypForS} will be very useful in the estimates for the temperature function $\theta$ and its derivatives.
\end{rmk}

\begin{rmk} \label{r:hyp-theta}
 Observe that $s>r>1$, hence $2-2/s>0$. Therefore, keeping Remark \ref{r:Besov} in mind, the assumption $\theta_0\in B^{2-2/s}_{p/2,s}$
in particular implies that
\[
 \theta_0\,\in\,L^{p/2}\qquad \mbox{ and }\qquad \Delta\theta_0\,\in\,\dot B^{-2/s}_{p/2,s}\,,
\]
which are the crucial properties we will use for the \tsl{a priori} estimates, see in particular the computations of Subsection \ref{ss:apriori-temp}.

In addition, if one furthermore assumes that $2-2/s-d/p>0$, which in turn reduces to strengthen condition \eqref{HypForS}
to $1/s\,<\,1-d/(2p)$,
by Besov embeddings of Proposition \ref{p:embed} we immediately get that $B^{2-2/s}_{p/2,s}\,\hookrightarrow\,L^p$, so the assumption
$\theta_0\in L^p$ would be satisfied for free.
\end{rmk}

Let us now give an overview of the strategy of the proof to Proposition \ref{PropositionEstimatesForModifiedSolution}.
In this section, we only focus on \tsl{a priori} estimates. The actual proof of existence is strictly tied to the proof of the existence statement in
Theorem \ref{th:well-p} and will be carried out in Section \ref{s:existence}.

The proof of \tsl{a priori estimates} in Proposition \ref{PropositionEstimatesForModifiedSolution} is based on a bootstrap argument,
which we are going to present in some detail.
First of all, we need to fix a time-integrability index $\s>1$ such that
\begin{equation}\label{ConditionSigma}
s\,<\,\sigma\,<\, 2\,s\,\left(1\, +\, \dfrac d{p}\right)^{-1}\,. 
\end{equation}
Thanks to the first inequality appearing in Remark \ref{r:s}, $\s$ is well-defined. This index of integrability will come into play in the definition of
the norm $\mc N(T)$ of the solution, in order to measure time integrability of the $L^\infty$ norm of $\nabla v\,=\,-\nabla^2(\Id-\Delta)^{-1}P$:
it turns out that this term has a better time integrability (see Lemma \ref{LemmaRegularitytheta} and Subsection \ref{SecTanReg}),
but the index $\sigma$ plays a crucial role in the bootstrap argument, in order to close the estimates on some small time interval $[0,T^*]$.

This having been clarified, let us define the quantity $\mc N$ as the sum of the relevant functional norms to be controlled in the proof of
Proposition \ref{PropositionEstimatesForModifiedSolution}. For any time $T>0$, we set
\begin{align}
\label{DefinitionOfN} 
	\mathcal N(T)\,&:=\,\Vert w \Vert_{L_T^\infty \left(\dot B^{2-2/r}_{p,r}\right)\cap L_T^{r_0}\left( L^{p_0}\right)}\,+\,
	\Vert \nabla w \Vert_{L_T^{r_1}\left( L^{p}\right)}\,+\,  \Vert \left( \partial_t w\,,\, \nabla^2 w\right) \Vert_{L_T^{r}\left( L^{p}\right)} \\ 
 \nonumber 	&\qquad\qquad\qquad\qquad\qquad\quad +\,  \Vert  \zeta \Vert_{L_T^{\infty}\left( L^{p}\cap L^\infty\right)} \, +\,
 \Vert \nabla^2 \left( \Id-\Delta\right)^{-1}P\Vert_{L_T^{\sigma}\left( L^\infty\right)}  \\
\nonumber
&\qquad\qquad\qquad\qquad\qquad\qquad\qquad +\,\Vert \theta \Vert_{L_T^{\infty}\left( B^{2-2/s}_{p/2,s}\cap L^{p}\right)} \, +\,
 \left\Vert \left( \partial_t \theta\,,\, \nabla^2\theta\right)\right\Vert_{L_T^s\left( L^{p/2}\right)}\,.
\end{align}
We naturally set $\mc N_0\,:=\,\mc N(0)$ to be the same quantity computed on the initial datum, namely
\[
 \mc N_0\,:=\,
 \Vert w_0 \Vert_{\dot B^{2-2/r}_{p,r}}\, +\,  \Vert  \zeta_0 \Vert_{L^{p}\cap L^\infty}\, +\,  \Vert \theta_0 \Vert_{B^{2-2/s}_{p/2,s}\cap L^{p}}\,.
\]

\begin{rmk} \label{r:norm-temp}
In light of Remarks \ref{r:Besov} and \ref{r:hyp-theta}, the quantity $\|\theta\|_{L^\infty_T(B^{2-2/s}_{p/2,s})}$
could have be splitted into $\left\|\theta\right\|_{L^\infty_T(L^{p/2})}\,+\,\left\|\theta\right\|_{L^\infty_T(\dot B^{2-2/s}_{p/2,s})}$,
and analogue splitting for the norm of the initial datum.

In fact, in the estimates below, we are precisely going to bound the two quantities, $\left\|\theta\right\|_{L^\infty_T(L^{p/2})}$
and $\left\|\theta\right\|_{L^\infty_T(\dot B^{2-2/s}_{p/2,s})}$, separately, by the use of different techniques. Summing up these two quantities,
we will get the control on the full non-homogeneous Besov norm.
\end{rmk}

It is worth remarking that the quantity $\mc N(T)$ already contains the relevant information necessary to bound the Lipschitz norm of the velocity field $u\,=\,w\,+\,v$,
see also Lemmas \ref{LemmaRegularityv} and \ref{LemmaRegularityw} in this respect.
As we will see, the property $\nabla u\in L^{m_1}_T(L^\infty)$, for a suitable $m_1>1$, will follow by precise maximal regularity estimates for $w$ and
$\theta$, combined with tangential regularity estimates for the part of $v$ associated with the density perturbations $\z\,=\,\rho\,-\,1$.
Contrarily to \cite{DFP}, then, tangential regularity and the information $\nabla u\in L^{m_1}_T(L^\infty)$ are necessary not only for uniqueness,
but also for the existence statement in Proposition \ref{PropositionEstimatesForModifiedSolution}.

\medbreak
As a bootstrap hypothesis for bounding $\mc N$, we shall assume that $\mc N(T)$ is finite on a time interval $[0, T_0]$.
By finding appropriate estimates for each term appearing in its definition, we will conclude that there exist $\g>0$ and $k\in\N$ such that 
\begin{equation}\label{BootstrapTesis}
	\mathcal N(T)\,\lesssim\,\mc N_0\,+\, T^{\g}\,\Big( 1\,+\,\mathcal N(T)\Big)^k\,\mathcal N^2(T)\,.	
\end{equation}
It is then clear that, with \eqref{BootstrapTesis} at hand, we can infer the existence of a time $T^*\in\,]0,T_0]$ such that
$\sup_{T\in[0,T^*]}\mc N(T)\,<\,+\,\infty$, thus yielding the sought \tsl{a priori}
bounds necessary to prove Proposition \ref{PropositionEstimatesForModifiedSolution}, hence in turn Theorem \ref{th:well-p}.

To prove \eqref{BootstrapTesis}, we will first use $L^q$ estimates for the transport equation for $\z$ and the trasport-diffusion equation for $\theta$. Next,
following the approach of \cite{DFP}, we will use maximal regularity estimates to control $w$, $\nabla w$ and $\nabla^2 w$, as well as $\nabla^2\theta$.
The most delicate part will consist in bounding $\nabla v=-\nabla^2 \left( \Id\,-\,\Delta\right)^{-1}P(\rho, \theta)$ in $L_T^\sigma\left( L^\infty\right)$,
for which we will resort to striated regularity \tsl{\`a la Chemin} \cite{Ch1991, Ch1993}.
Specifically, for doing this we will need to introduce also another quantity, which we will call $\mc Q(T)$, precisely taking into account the norms
related to tangential regularity. Then, an estimate of the type \eqref{BootstrapTesis} will be in fact proved not for $\mc N(T)$, but rather for
$\mc E(T)\,:=\,\mc N(T)\,+\,\mc Q(T)$.

\medbreak
Before going on, let us point out the following fact.
The smallness assumption \eqref{HypothesisForZetaPropositionEstimatesForModifiedSolution} is crucial for our approach to work, especially
in order to perform maximal regularity estimates for $w$ and $\theta$.
Taking that condition into account, for any $\de>\varepsilon$ we are entitled to define 
\[
T_{(\de)}\,:=\, \sup\left\{T>0\;\Big|\quad \Vert \zeta(t)\Vert_{L^\infty}\,\le\, \de \quad \forall\, t \in [0, T] \right\}\,.	
\] 
On the other hand, since the functions $\pi_0$ and $\pi_1$ are smooth by assumption (in fact, $C^1(\R_+)$ is enough),
with $\pi_0(1)=0$, for any fixed $\delta>0$ there exists $C_*>0$ such that 
\begin{align}
\label{UniformBoundForPi0(2)}
\forall\, \rho \in [1-\de, 1+\de]\,,\qquad\qquad
\vert \pi_0(\rho)\vert\,+\,\big| \pi_1(\rho)-\pi_1(1)\big|\, &\le\, C_*\,\vert \z\vert\,, 
\end{align}
with $\z\,=\,\rho-1$.
In what follows, we will take $\de\,=\,4\veps$ and define
\begin{equation} \label{def:T_0}
 T_0\,:=\,T_{(4\veps)}\,=\,\sup\left\{T>0\;\Big|\quad \Vert \zeta(t)\Vert_{L^\infty}\,\le\, 4\,\veps \quad \forall\, t \in [0, T] \right\}\,.
\end{equation}

\subsection{Static estimates} \label{ss:static}
We collect in this subsection a series of static estimates which derive from assuming that $\mathcal N(T)$ is finite on $[0,T_0]$.
These estimates will be very helpful in the bootstrap argument.

The first result concerns some integrability properties for the temperature function $\theta$. We postpone an important remark at the end of the proof.

\begin{lemma} \label{LemmaRegularitytheta}
Let $T\in[0,T_0]$ such that $\mathcal N(T)<+\infty$. 
Define the two integrability indices $s_0\,>\,s_1\,>\,s$ by the formulas
\begin{equation*}
 s_0\,:=\,\frac{2\,p}{d}\,s\,\left(1\,+\,\frac{d}{p}\right)^{-1}\qquad \mbox{ and }\qquad
 s_1\,:=\,2\,s\,\left(1\,+\,\frac{d}{p}\right)^{-1}\,.
\end{equation*}

Then, for an implicit universal multiplicative constant $C>0$, one has
\[
\left\Vert \theta \right\Vert_{L_T^{s_0}\left( L^\infty\right)}\,+\,
\left\Vert \nabla \theta \right\Vert_{L_T^{s_1}\left( L^p\right)\cap L_T^{2s}\left( L^{p/2}\right)}\,\lesssim\,\mathcal N(T)\,.
\]
\end{lemma}

\begin{proof}
We use the Gagliardo-Nirenberg inequalities of Lemma \ref{l:GN-ineq}. Setting $\alpha$ as in that statement, we infer that 
\begin{align*}
\left\Vert\nabla \theta \right\Vert_{L_T^{s/\alpha}\left(L^{p}\right)}\,&\lesssim\,
\left\Vert \theta \right\Vert_{L_T^\infty\left(L^{p/2}\right)}^{1-\alpha}\, \left\Vert \Delta \theta \right\Vert_{L_T^{s}\left(L^{p/2}\right)}^{\alpha}\,
\lesssim\,\mathcal N(T)\,,  \\
\left\Vert\nabla \theta \right\Vert_{L_T^{2s}\left(L^{p/2}\right)}\,&\lesssim\, \left\Vert \theta \right\Vert_{L_T^\infty\left(L^{p/2}\right)}^{1/2}\,
\left\Vert \Delta \theta \right\Vert_{L_T^s\left(L^{p/2}\right)}^{1/2}\,\lesssim\, \mathcal N(T)\,.  
\end{align*}
Notice that from the definition of $\alpha$ in Lemma \ref{l:GN-ineq}, we have that
\[
\frac{1}{s_1}\,=\,\dfrac{\alpha}s\, <\, \dfrac{\alpha}r\, \le\, 
\dfrac2{r_1}\,=\,\dfrac{2}{r}\,-\,1\,.
\]

Next, from Lemma \ref{LemmaUsefulInequalityGeneralVersion}, we infer that
\[
\left\Vert \theta \right\Vert_{L^\infty}\,\lesssim\,
\left\Vert \theta \right\Vert^{\alpha_0}_{L^p}\,\left\Vert \nabla \theta \right\Vert_{L^p}^{1-\alpha_0}\,, \qquad \mbox{ with }\quad  \alpha_0\, :=\, 1\,-\,\dfrac dp\,. 
\]
Thus, the time integrability index $s_0$ is given by the following formula: 
\[
\dfrac1{s_0}\,:=\,\dfrac{\alpha\,\left( 1-\alpha_0\right)}{s}\, =\, \dfrac{\left( 1+ d/p\right)d/p}{2\,s}\,, 
\]
as claimed. Notice that $s_0$ satisfies
\[
\dfrac1{s_0}\,<\,\frac{1}{s_1}\,. 
\]
The proof of the lemma is thus completed.
\end{proof}

\begin{rmk} \label{r:s_0-s_1}
We stress the fact that, owing to the condition $p>d$ appearing in \eqref{HypForPs} and keeping \eqref{HypForS} in mind,
we get the series of inequalities $s_0>s_1>s>r$.

Finally, we point out that $\s$, defined as in \eqref{ConditionSigma}, satisfies
\[
s\,<\,\s\,<\,s_1\,. 
\]
\end{rmk}

Next, let us present some useful estimates for the pressure function $P=P(\rho,\theta)$ and for the velocity field
$v\,=\,-\nabla(\Id-\Delta)^{-1}P$ (recall that we have previously set $\nu\,=\,\lambda+2\mu\,=\,1$).
As before, we postpone relevant remarks on the statement at the end of its proof.

\begin{lemma} \label{LemmaRegularityv}
Let $T \in [0, T_0]$ such that $\mathcal N(T)<+\infty$. Let the time integrability index $s_0>s$ be defined as in Lemma \ref{LemmaRegularitytheta}
and let $\s$ satisfy the condition in \eqref{ConditionSigma}.

Then, there exists a (implicit) universal multiplicative constant $C>0$ such that 
\begin{equation}\label{LemmaRegularityveq1}
\left\Vert P(\rho, \theta) \right\Vert_{L_T^\infty\left( L^p\right)}\,\lesssim\,\mathcal N(T)\qquad \mbox{ and }\qquad
\left\Vert P(\rho, \theta) \right\Vert_{L_T^{s_0}\left( L^\infty\right)}\,\lesssim\,1\,+\, \mathcal N(T)\,. 
\end{equation}

Moreover, let the vector field $v$ be defined as in \eqref{Defv}. Then, for a suitable constant $C>0$ as above, it holds
\[
\left\Vert v \right\Vert_{L_T^\infty\left( L^p\cap L^\infty\right)}\, +\,
\left\Vert \divergence \, v \right\Vert_{L_T^{s_0}\left( L^\infty\right)}\, +\,
\left\Vert \nabla v \right\Vert_{L_T^\infty\left( L^{p}\right)\cap L_T^{\sigma}\left( L^\infty\right)}\,\lesssim\,\mathcal N(T)\,.
\]
\end{lemma}

\begin{proof}
Recall that the time $T_0$ has been fixed as the maximum time $t>0$ for which $\big|\rho-1\big|\leq 4\veps$ on $[0,t]\times\R^d$. In particular,
this bounds holds true also at the time $T\in[0,T_0]$ given in the statement.
Hence, from the definition of $P(\rho, \theta)= \pi_0(\rho)+ \pi_1(\rho)\theta$, we can write
\begin{align*}
\left\Vert P \right\Vert_{L^\infty}\,\leq\,\left\|\pi_0(\rho)\right\|_{L^\infty}\,+\,\left\|\pi_1(\rho)\right\|_{L^\infty}\,\left\|\theta\right\|_{L^\infty}\,,
\end{align*}
which immediately yields that
\[
 \left\Vert P \right\Vert_{L^{s_0}_T\left(L^\infty\right)}\,\lesssim\,1\,+\,\mc N(T)\,,
\]
owing to Lemma \ref{LemmaRegularitytheta} and the continuity of
the functions $\pi_0$ and $\pi_1$ on the interval $[1-4\veps,1+4\veps]$.

Similarly, for the $L^p$ norm in space, we can use again the definition of $P$, together with the fact that $\pi_0(1)=0$. Hence,
arguing as in \eqref{UniformBoundForPi0(2)} and recalling that $\z\,=\,\rho-1$, we get
\begin{align*}
\left\Vert P \right\Vert_{L^p}\,\lesssim\,\left\|\z\right\|_{L^p}\,+\,\left\|\pi_1(\rho)\right\|_{L^\infty}\,\left\|\theta\right\|_{L^p}\,,
\end{align*}
whence we deduce that
\[
 \left\Vert P \right\Vert_{L^\infty_T\left(L^p\right)}\,\lesssim\,\mc N(T)\,.
\]

Now, let us recall the definition of $v$ as 
\[
v= -\nabla \left( \Id -\Delta\right)^{-1}P(\rho, \theta)\,=\,-\,\nabla \left( \Id -\Delta\right)^{-1}\Big( \pi_0(\rho) \,+\, \pi_1(\rho)\,\theta\Big)\,.
\]
Then, thanks to the continuity of the operator $\nabla \left( \Id-\Delta\right)^{-1}:L^p\to W^{1,p}$ (see Lemma \ref{l:D-1} and Corollary
\ref{c:Delta-Id} above) and \eqref{LemmaRegularityveq1}, we have that 
\[
\left\Vert \nabla \left( \Id-\Delta\right)^{-1}P \right\Vert_{L_T^\infty\left( W^{1,p}\right)}
\,\lesssim\,\left\Vert P \right\Vert_{L_T^\infty\left( L^p\right)}\,\lesssim\,\mathcal N(T)\,.
\]
Notice that, by definition \eqref{DefinitionOfN} of $\mc N(T)$, one has
\[
\left\Vert \nabla v \right\Vert_{L_T^\sigma\left( L^\infty\right)}\,=\,
\left\Vert \nabla^2\left( \Id-\Delta\right)^{-1} P \right\Vert_{L_T^\sigma\left( L^\infty\right)}\,\le\, \mathcal N(T)\,.
\]
Now, to deal with $\divergence v$, we use that the operator $\left( \Id-\Delta\right)^{-1}:L^\infty\to L^{\infty}$ is continuous
(as it can be deduced from Corollary \ref{c:Delta-Id}): thanks to \eqref{LemmaRegularityveq1}, we find that
\[
\divergence v = -\Delta \left( \Id - \Delta\right)^{-1}P =  P -  \left( \Id - \Delta\right)^{-1}P \in L_T^{s_0}\left( L^\infty\right).
\]
To handle $\Vert \nabla \left( \text{Id}-\Delta\right)^{-1}P\Vert_{L_T^{\infty}\left( L^\infty\right)}$, we use the interpolation inequality given in Lemma \ref{LemmaUsefulInequalityGeneralVersion} and the continuity of the operator $\nabla \left( \Id-\Delta\right)^{-1}:L^p\to L^p$ (see again
Lemma \ref{l:D-1} and Corollary \ref{c:Delta-Id}) to get
\begin{align}
\label{EstWAux3}
\Vert \nabla \left( \Id-\Delta\right)^{-1}P\Vert_{L^\infty}\,&\lesssim\,\Vert \nabla \left( \Id-\Delta\right)^{-1}P\Vert^{1-d/p}_{L^p}\,
\Vert \nabla^2 \left( \Id-\Delta\right)^{-1}P\Vert^{d/p}_{L^p}  \\ 
& \lesssim\, \Vert P\Vert^{1- d/p}_{L^p}\,\Vert \Delta \left( \Id-\Delta\right)^{-1}P\Vert^{ d/p}_{L^p}
 \nonumber \\
& \lesssim\,\Vert P\Vert^{1-d/p}_{L^p}\,\left( \Vert \left( \Id-\Delta\right)^{-1}P\Vert_{L^p}\, +\, \Vert P\Vert_{L^p}\right)^{d/p}\,.\nonumber
\end{align}
Then, thanks to \eqref{LemmaRegularityveq1}, we get also the sought bounds for $v$ and its derivatives. Thus,
the proof of Lemma \ref{LemmaRegularityv} is complete. 
\end{proof}

\begin{rmk} \label{r:divv-Dv}
Notice that the quantity $\div v$ enjoys better time integrability (with values in $L^\infty(\R^d)$) than $\nabla v$.
Indeed, as already pointed out in Remark \ref{r:s_0-s_1}, it follows from the definitions that we have
$s<\s<s_1<s_0$.
\end{rmk}

As a consequence of Lemma \ref{LemmaRegularityv}, we obtain suitable bounds for $\div u$ as well. Indeed, taking into account the decomposition
$u= w + v$, we get
\begin{equation}\label{EstimationP(1)}
\divergence u\, \in\, L_T^{r_1}\left( L^p\right)\,,
\end{equation}
for any $T\in\,]0,T_0]$ such that $\mc N(T)<+\infty$. Here $r_1$, defined in \eqref{DefOfRs}, is the integrability index pertaining to $\nabla w$.

As a matter of fact, by assuming $\mathcal N(T)<+\infty$, there is more information regarding the regularity and integrability of $w$.
For instance, Sobolev embeddings imply that $\nabla w \in L_t^r\left(L^\infty\right)$; however, the integration in time can be improved by an interpolation
argument.
All this is the matter of the following lemma. 

\begin{lemma} \label{LemmaRegularityw}
Let $T \in [0, T_0]$ such that $\mathcal N(T)<+\infty$.
Define the integrability indices $r_{0,\infty}\in\,]r_0,+\infty[\,$ and $r_{1,\infty}\in\,]r_1,+\infty[\,$ by the formulas
\[
\frac1{r_{0, \infty}}\,=\, \frac{1}{r}\,+\,\frac{d}{2p}\,-\,1 \qquad \mbox{ and }\qquad
\frac1{r_{1,\infty}}\,=\, \frac{1}{r}\,-\,\frac12\,+\,\frac d{2p}\,.
\]

Then, there exists an implicit universal multiplicative constant $C>0$ such that 
\begin{gather*}
\Vert w \Vert_{L_T^{r_{0,\infty}}\left( L^\infty\right)}\,+\,\Vert \nabla w \Vert_{L_T^{r_{1, \infty}}\left( L^\infty\right)}\, \lesssim\,
\mathcal N(T)\,. 
\end{gather*}
\end{lemma}

\begin{proof}
First of all, we notice that we can write
\[
\frac1{r_{1,\infty}}\,=\,\frac{1}{r_1}\,+\,\frac{d}{2p}\,=\, \frac{1-d/p}{r_1}\, +\, \frac{d/p}{r}\,. 
\]
Then, using the interpolation inequalities stated in Lemma \ref{LemmaUsefulInequalityGeneralVersion}, we have
\begin{equation*} 
\Vert \nabla w \Vert_{L_T^r\left(L^\infty\right)}\,\leq\,\Vert \nabla w\Vert_{L_T^{r_1}\left( L^p\right)}^{1-\frac dp}\,
\left\Vert \nabla^2 w \right\Vert_{L_T^{r}\left( L^p\right)}^{\frac dp}\,\lesssim\,\mc N(T)\,,
\end{equation*}
where we have applyed the Young inequality in the last step.

Proceeding analogously for $w$, one finds that 
\[
\Vert w \Vert_{L^\infty}\,\leq\,\Vert w \Vert_{L^{p_0}}^\alpha\, \Vert \nabla w \Vert_{L^p}^{1-\alpha}\,, \qquad \mbox{ with }\quad 
\alpha\,:=\,\frac{1-d/p}{1+d/p_0-d/p}\,.
\]
This yields a time integrability $r_{0,\infty}$ for $\|w\|_{L^\infty}$ which can be computed as follows:
\begin{align*}
\frac{1}{r_{0,\infty}}\,:=\,\frac{\alpha}{r_0}\,+\,\frac{1-\alpha}{r_1}\,&=\,\frac{1-d/p}{1+d/p_0-d/p}\,\left( \frac d{2p}- \frac d{2p_0}-\frac12\right) + \frac1{r_1}
\,=\,\frac{1}{r}\,+\,\frac{d}{2p}\,-\,1\,.
\end{align*}
Observe that $r_{0,\infty}$ is independent of $p_0$. This concludes the proof of the lemma.
\end{proof}

As a direct consequence of Lemma \ref{LemmaRegularityw}, by arguing as done for \eqref{EstimationP(1)},
we can give a control of $\div u$ in $L^{\infty}$, which is a crucial property for estimating the Lebesgue norms of $\z\,:=\,\rho-1$ and, in turn,
for defining the time $T_0$ as in \eqref{def:T_0}. However, for later applications, we need to upgrade the $L^\infty$ bound to the full gradient $\nabla u$ of
the velocity field.

\begin{cor}\label{CorRegularityu}
Let $T \in [0, T_0]$ such that $\mathcal N(T)<+\infty$. Assume moreover that the quantity $\nabla v\,=\,-\,\nabla^2 \left( \Id - \Delta\right)^{-1}P$ belongs to
$L^{s_1}_T(L^\infty)$, where the index $s_1$ has been introduced in Lemma \ref{LemmaRegularitytheta}.
Then, for a constant $C>0$ as in Lemma \ref{LemmaRegularityw}, one has
\begin{align*}
\left\|u\right\|_{L^{r_{0,\infty}}_T(L^\infty)}\,+\,
\left\|\nabla u\right\|_{L^{r_1}_T(L^p)}\,+\,\Vert\div \, u \Vert_{L_T^{r_{1, \infty}}\left( L^\infty\right)}
\,+\,\left\|\nabla u\right\|_{L^{m_1}_T(L^\infty)}
\,&\lesssim\,\mathcal N(T)\,, 
\end{align*}
where we have defined $m_1\,:=\,\min\big\{r_{1,\infty}\,,\,s_1\big\}$.
\end{cor}

\begin{proof}
The claimed estimates are a direct consequence of the decomposition $u\,=\,w\,+\,v$, combined with Lemmas \ref{LemmaRegularityv} and \ref{LemmaRegularityw}.
Notice that this would only imply $\nabla u\in L^{m'}_T(L^\infty)$, with $m'\,:=\,\min\big\{r_{1,\infty}\,,\,\s\big\}$; it is precisely
here where we use the additional assumption 
\[
\nabla v\,=\,-\,\nabla^2 \left( \Id - \Delta\right)^{-1}P\;\in\, L_T^{s_1}( L^\infty)\,,
\]
from which we deduce the sought property for $\nabla u$. Notice that this condition
is stronger than the one contained in the assumption $\mc N(T)<+\infty$. A control 
for $\nabla v$ in $L^{s_1}_T(L^\infty)$ will be shown in Subsection \ref{SecTanReg} by use of tangential regularity.

Concerning the bound for the term $\div u$, we have only to verify that the integrability index $r_{1,\infty}$, defined in Lemma \ref{LemmaRegularityw}, satisfies
$r_{1,\infty}\leq s_0$, where $s_0$ has been defined in Lemma \ref{LemmaRegularitytheta}.
Establishing the inequality is a matter of performing some computations. First of all, we notice that
\[
\dfrac12\,-\,\dfrac d{2p}\, -\, \dfrac {d^2}{4p^2} \,+\, \dfrac {d^3}{4p^3}\,=\,\left(1\,-\,\frac{d}{d}\right)\,\left(\frac{1}{2}\,-\,\frac{d^2}{4p^2}\right)
\,>\,0\,.
\]
Equivalently, by adding the quantity $\frac{1}{2}\,-\,\frac{d}{2p}$ to both sides of the inequality, we can write
\[
\left(\dfrac12\, -\, \dfrac d{2p} \right)\,\left(1\,-\, \dfrac d{2p}\, -\, \dfrac {d^2}{2p^2}\right)^{-1}\,<\, 1\,-\,\dfrac d{2p}\,,
\]
which implies, owing to \eqref{HypForPs}, that
\[
\left(\dfrac12\, -\, \dfrac d{2p} \right)\,\left(1\,-\, \dfrac d{2p}\, -\, \dfrac {d^2}{2p^2}\right)^{-1}\,<\, \dfrac1r\,.
\]
In other words, we have
\[
\dfrac1{r}\,\left(1\,-\, \dfrac d{2p} \,-\, \dfrac {d^2}{2p^2}\right)\,>\, \dfrac12\,-\, \dfrac d{2p}\,, 
\]
which is the same as having 
\[
\dfrac1{r}\,\left( \dfrac d{2p}\, +\, \dfrac {d^2}{2p^2}\right)\, <\, \dfrac1r\, -\,\dfrac12 \,+\, \dfrac d{2p}\,.
\]
Now, because $s> r$, using the definition of $s_0$ we find
\begin{align*}
\dfrac1{s_0}\,=\,\dfrac1{2s}\,\left( \dfrac dp \,+\, \dfrac {d^2}{p^2}\right)\,<\,
\dfrac1{r}\,\left( \dfrac d{2p}\, +\, \dfrac {d^2}{2p^2}\right)\,<\, \dfrac1r\, -\,\dfrac12\, +\, \dfrac d{2p}\, =\, \dfrac1{r_{1, \infty}}\,.
\end{align*}
In particular, we deduce that $r_{1,\infty}<s_0$ and this completes the proof of the corollary.
\end{proof}

After having collected the previous static estimates, we finally tackle the proof of Proposition \ref{PropositionEstimatesForModifiedSolution},
which will be the matter of the following subsections.

\subsection{Estimates for the density}\label{SubsectionAPrioriEstimatesForDensity}
We begin by presenting the \tsl{a priori} estimates for the density perturbation function $\z\,=\,\rho-1$.

Owing to the initial jump discontinuity in $\rho_0$, which will be propagated also for later times, we can hope only for controls on Lebesgue norms of $\z$.
Thus, take $q \in [1, +\infty[\,$. We multiply equation \eqref{EquationForZeta} by $\vert \zeta\vert^{q-2}\zeta$ and integrate in space, obtaining
\begin{align*}
\frac1q\, \frac{d}{dt} \Vert \zeta \Vert_{L^{q}}^q\, &=\,-\,\frac1q \int_{\R^d}\div u\,\vert \zeta\vert^{q}\,dx\,+\,
\int_{\R^d} (\z+1)\,\div u\,\vert \zeta \vert^{q-2}\zeta\,dx \\
&\leq\,\left(1\,-\,\frac{1}{q}\right)\,\left\|\div u\right\|_{L^\infty}\,\left\|\z\right\|_{L^q}^q\,+\,\left\|\div u\right\|_{L^q}\,\left\|\z\right\|_{L^{q'}}\,,
\end{align*}
where $q'$ denotes the conjugate exponent of $q$, that is $1/q\,+\,1/q'\,=\,1$.
Then, integrating in time, for any $1\le q<+\infty$ and any time $t\geq0$ we find
\begin{align*}
\Vert \zeta(t) \Vert_{L^{q}}\,&\leq\, \Vert \zeta_0 \Vert_{L^{q}} + \left( 1- \frac1q \right)\,
\int_0^t \Vert \div u (\tau)\Vert_{L^\infty}\,\Vert \zeta(\tau)\Vert_{L^q}\,d\tau \\
&\qquad\qquad\qquad\qquad\qquad +\, \int_0^t \Big(\Vert \divergence w(\tau)\Vert_{L^q}\,+\,\Vert \div v(\tau)\Vert_{L^q}\Big)\,d\tau\,.
\end{align*}

Particularly, taking the limit $q\to +\infty$ in the previous inequality and applying the Gr\"onwall lemma, we infer that
\begin{align}
\label{est:z_inf}
\Vert \zeta(t) \Vert_{L^{\infty}}\,&\lesssim\,\left(\Vert \zeta_0 \Vert_{L^{\infty}}\,+\,
\int^t_0\Big(\Vert \divergence w(\tau)\Vert_{L^\infty}\,+\,\Vert \div v(\tau)\Vert_{L^\infty}\Big)\,d\tau\right) \\
\nonumber
&\qquad\qquad\qquad\qquad\qquad\qquad\qquad \times
\exp\left(\int_0^t \Vert \divergence u (\tau)\Vert_{L^\infty}\,d\tau\right)\,.
\end{align}
Analogously, after setting $C_p\,:=\,1\,-\,1/p$, we also get
\begin{align}
\label{EstLq_Zeta}
\Vert \zeta(t) \Vert_{L^{p}}\,&\lesssim\,\left(\Vert \zeta_0 \Vert_{L^{p}}\,+\,
\int^t_0\Big(\Vert \divergence w(\tau)\Vert_{L^p}\,+\,\Vert \div v(\tau)\Vert_{L^p}\Big)\,d\tau\right) \\
\nonumber
&\qquad\qquad\qquad\qquad\qquad\qquad\qquad \times
\exp\left(C_p\,\int_0^t \Vert \divergence u (\tau)\Vert_{L^\infty}\,d\tau\right)\,.
\end{align}

At this point, we define the time $T_1>0$ such that 
\begin{align}\label{DefinitionOfT}
T_1\,&:=\,\sup \left\{ t\,>\,0\quad \Big|\quad \int_0^t \Vert \div u(\tau)\Vert_{L^\infty}\, d\tau\,\leq\,\log 2 \right. \\ 
\nonumber
& \hspace{3.5cm} \left. \mbox{ and }\quad \int_0^t \Big( \Vert \div w(\tau) \Vert_{L^\infty}\,+\,\Vert \div v(\tau)\Vert_{L^\infty}\Big)\,d\tau\,\leq\,
\left\|\z_0\right\|_{L^\infty}\right\}\,.
\end{align}
Owing to hypothesis \eqref{HypothesisForZetaPropositionEstimatesForModifiedSolution}, we then find that 
\begin{equation}\label{EstZeta2}
\forall\,t\in[0,T_1]\,,\qquad	\Vert \zeta (t)\Vert_{L^\infty}\,\le\,4\,\left\|\z_0\right\|_{L^\infty}\,\leq\, 4\,\varepsilon\,.
\end{equation}
Therefore, by maximality of the time $T_0$ defined in \eqref{def:T_0}, we infer that
\[
0\,<\, T_1\,\leq\,T_0\,,
\]
so that in $[0,T_1]$ we can guarantee the validity of all the estimates stated in Subsection \ref{ss:static}.
In addition, from \eqref{EstZeta2} we deduce the inequalities stated in \eqref{UniformBoundForPi0(2)}, specifically
\begin{equation}\label{UniformBoundForPi_0-1}
\forall\,t\in[0,T_1]\,,\qquad	\big\vert \pi_0(\rho)(t)\big\vert\, +\,\big\vert \pi_1(\rho)(t)-\pi_1(1)\big\vert\,\lesssim\,\vert \zeta(t)\vert\,.
\end{equation}

Observe that, thanks to the uniform bound \eqref{EstZeta2} and the decomposition $\rho\,=\,1+\z$, we can also estimate
\begin{equation} \label{est:1-su-rho}
\forall\,T\in[0,T_1]\,,\qquad \left\|\frac{1}{\rho}\right\|_{L^\infty_T(L^\infty)}\,\leq\,\frac{1}{1\,-\,4\,\veps}\,\leq\,\frac{1}{1\,-\,4\,\veps_0}\,\leq\,2\,,
\end{equation}
where, without loss of generality, we can assume that $\veps_0>0$ in the statement of Proposition \ref{PropositionEstimatesForModifiedSolution}
has been chosen so small that $1-4\veps_0\geq1/2$.

Finally, as direct consequence of inequality \eqref{EstLq_Zeta} and Lemma \ref{LemmaRegularityv},
we infer that
\begin{align*}
\forall\,T\in[0,T_1]\,,\qquad
\Vert \zeta \Vert_{L_{T}^\infty\left( L^p\right)}\,&\lesssim\, \Vert \zeta_0\Vert_{L^p}\,+\,
\left( {T}^{1-1/r_1}\,\Vert \divergence w \Vert_{L^{r_1}_{T}\left( L^p\right)}\, +\, {T}\,\Vert \nabla v \Vert_{L_{T}^\infty\left( L^p\right)}\right)\,.
\end{align*}
In turn, owing  to the definition \eqref{DefinitionOfN} of $\mc N(t)$,
we finally deduce that
\begin{equation}\label{EstLp_Zeta}
\forall\,T\in[0,T_1]\,,\qquad
\Vert \zeta \Vert_{L_{T}^\infty\left( L^p\right)}\,\lesssim\,\Vert \zeta_0\Vert_{L^p}\, +\,{T}^{1-1/r_1}\,\big( 1+ \mathcal N({T})\big)\,.
\end{equation}

\subsection{Estimates for the temperature} \label{ss:apriori-temp}

Here, we prove estimates for the norms of the temperature $\theta$ appearing in \eqref{DefinitionOfN}. First of all, we show
bounds for the Lebesgue norms of $\theta$ by means of transport-diffusion type estimates. Only afterwards we will bound $\Delta \theta$ by resorting
to maximal regularity for the heat equation. 

\subsubsection{Basic estimates for the temperature} \label{sss:Leb-temp}
In order to perform basic estimates on the Lebesgue norms of $\theta$, we rewrite
its equation in \eqref{FNSwithTemperatureModified} as a transport-diffusion equation, namely
\begin{equation} \label{eq:temp_new}
\rho\,\partial_t \theta\, +\, \rho\,u\cdot \nabla \theta\, -\,\k\,\Delta \theta
		\,=\,-\,\pi_1(\rho)\,\theta\,\div u\,+\, \mc B(\nabla u)\,,
\end{equation}
where the bilinear term $\mc B(\nabla u)$ has been defined in \eqref{def:B_u}.

First of all, we observe that $\rho\geq1-4\veps>0$, owing to \eqref{EstZeta2}, and $\mc B(\nabla u)\geq 0$ too. Hence, after dividing the previous equation
by $\rho$, the standard maximum principle for second order parabolic equations
implies that $\theta(t,x)\geq0$ for all $(t,x)\in[0,T_1]\times\R^d$ (we refer to \cite{It_1971} for details, see the end of Section 5 therein).

Next, as $\theta\geq0$, we can multiply both members in equation \eqref{eq:temp_new} by $\theta^{q-1}$, for any finite $q\geq1$. Using
the continuity equation for the density $\rho$ and the fact that the contribution coming from the diffusion term $-\k\Delta\theta$ has positive sign,
and keeping the bounds \eqref{EstZeta2} and \eqref{UniformBoundForPi_0-1} in mind, we get
\begin{equation} \label{est:dt_theta}
\frac{d}{dt}\left\|\rho^{1/q}\,\theta\right\|_{L^q}\,\lesssim\,\left\|\div u\right\|_{L^\infty}\,\left\|\theta\right\|_{L^q}\,+\,\left\|\mc B(\nabla u)\right\|_{L^q}\,.
\end{equation}
It is then a matter of bounding $\mc B(\nabla u)$ in $L^q$, for suitable values of $q\geq1$.
Observe that 
\begin{equation*}
\left\Vert \mc B(\nabla u) \right\Vert_{L^q}\,\lesssim\, \left\Vert \left\vert \nabla u \right\vert^2\right\Vert_{L^q}\,.
\end{equation*}
Hence, we are left to prove that $\big|\nabla u\big|^2\,\in\,L_T^1( L^q)$, for suitable values of $q$. In fact, we are going to show that the previous property
holds true with $q=p$ and with $q=p/2$.

To begin with, we use Corollary \ref{CorRegularityu} and the fact that $r_1/2>1$, a property which follows directly from the definition of $r_1$ in \eqref{DefOfRs}
and the fact that $r<2$, owing to \eqref{HypForPs}.
Hence, we can bound
\begin{equation*}
 \left\|\mc B(\nabla u)\right\|_{L^{r_1/2}_T(L^{p/2})}\,\lesssim\,\left\|\nabla u \right\|_{L^{r_1}_T(L^{p})}^2\,\lesssim\,\mc N^2(T)\,.
\end{equation*}
At the same time, we can take advantage again of Corollary \ref{CorRegularityu} and observe that
\[
\frac{1}{s_2}\,:=\, \frac{1}{m_1}\,+\,\frac{1}{r_1}\,<\,1\,.
\]
Indeed, when $m_1=r_{1,\infty}$, this relation is a direct computation; when $m_1=s$, instead, the validity of this inequality is ensured by the second condition
appearing in \eqref{HypForS}. Therefore, we deduce the estimate
\begin{equation*}
 \left\|\mc B(\nabla u)\right\|_{L^{s_2}_T(L^{p})}\,\lesssim\,\left\|\nabla u \right\|_{L^{r_1}_T(L^{p})}\,
 \left\|\nabla u \right\|_{L^{m_1}_T(L^{\infty})}\,\lesssim\,\mc N^2(T)\,.
\end{equation*}

Gathering the previous bounds together with \eqref{est:dt_theta}, after an integration in time we get
\begin{equation}\label{est:temp_L^p}
\forall\,T\in[0,T_1]\,,\qquad
\Vert \theta \Vert_{L_{T}^\infty\left(L^{p/2}\cap L^p\right)}\,\lesssim\,\Vert \theta_0\Vert_{L^{p/2}\cap L^p}\, +\,{T}^{\g_{r,s}}\,\mathcal N^2({T})\,,
\end{equation}
for a suitable positive exponent $\g_{r,s}>0$, whose precise value depends on $r_1$ and $s_2$, thus on $r$ and $s$ appearing in \eqref{HypForPs} and \eqref{HypForS}.
We observe that, in writing the previous inequality, we have also used the definition of the time $T_1$, given in \eqref{DefinitionOfT} above.


\subsubsection{Maximal regularity for $\Delta \theta$} \label{sss:MR-theta}
Thanks to the bounds established in \eqref{est:temp_L^p}, we can now prove that $\Delta \theta \in L_T^s\left( L^{p/2}\right)$. This property
will be inferred by a maximal regularity argument for equation \eqref{eq:temp_new}, after moving the transport term $\rho u\cdot\nabla\theta$ to the right-hand side.
As a matter of fact, an application of Lemmas \ref{LemmaMaximalRegularitySecondDerivative} and \ref{LemmaMaximalRegularityZeroDerivative} to that equation implies that
\begin{align*}
\left\|\theta\right\|_{L^\infty_T(\dot B^{2-2/s}_{p/2,s})}\,+\,
\left\Vert \left( \partial_t \theta, \Delta \theta\right) \right\Vert_{L_T^s\left( L^{p/2}\right)}
\,&\lesssim\,\left\|\Delta\theta_0\right\|_{\dot B^{-2/s}_{p/2,s}}\,+\,\left\Vert \zeta \partial_t\theta\right\Vert_{L_T^s\left( L^{p/2}\right)} + 
\left\Vert \rho u \cdot \nabla \theta \right\Vert_{L_T^s\left( L^{p/2}\right)} \\ 
& \qquad\quad + \left\Vert \pi_1(\rho)  \theta\divergence u \right\Vert_{L_T^s\left( L^{p/2}\right)} +
\left\Vert \mc B(\nabla u) \right\Vert_{L_T^s\left( L^{p/2}\right)} \,,
\end{align*}
where we recall that we have set $\z\,:=\,\rho-1$ and the bilinear operator $\mc B(\nabla u)$ has been defined in \eqref{def:B_u}.
It is clear that, up to reducing the value of $\veps_0>0$ in the statement of Proposition \ref{PropositionEstimatesForModifiedSolution},
the second term on the right-hand side of the previous inequality can be absorbed into the left-hand side.

The remaining terms on the right-hand side of the previous inequality can be bounded by the use of the following lemma, whose proof is postponed to the end
of this subsection.

\begin{lemma}\label{lem:RHS_MR_theta}
Let $s, r, r_0, r_1$ and $p$ be as in Proposition \ref{PropositionEstimatesForModifiedSolution} and let $\mathcal N(T)$ be the quantity defined in \eqref{DefinitionOfN}.
Then, there exists an implicit multiplicative constant $C>0$, depending only on $\left\|\theta_0\right\|_{L^{p/2}\cap L^p}$,
such that, for all $T\in[0,T_1]$, with $T_1$ defined in \eqref{DefinitionOfT}, one has
\begin{align*}
\left\Vert \Big(\rho u \cdot \nabla \theta\;,\; \pi_1(\rho)\,\theta\, \divergence u\;,\; \mc B(\nabla u)\Big) \right\Vert_{L_T^s\left( L^{p/2}\right)}\,
\lesssim\, T^{\gamma_{\rm temp}}\, \mc N^2(T)\,\Big( 1+ \mathcal N(T)\Big)^2\,,
\end{align*}
where we have defined the exponent $\gamma_{\rm temp}\,:=\,\dfrac{1}{2s}\,-\,
\dfrac1{r_1}\,>\,0$.
\end{lemma}

Notice that $\g_{\rm temp}$ is indeed strictly positive, according to the assumptions in \eqref{HypForS} and the definition of $r_1$ in \eqref{DefOfRs}.
As a direct consequence of Lemma \ref{lem:RHS_MR_theta}, it is easy to conclude that, for any time $T\in [0,T_1]$, one has
\begin{equation}\label{EstimationForTheta}
\left\|\theta\right\|_{L^\infty_T(\dot B^{2-2/s}_{p/2,s})}\,+\,
\left\Vert \left( \partial_t \theta, \nabla^2 \theta \right) \right\Vert_{L_T^s\left( L^{p/2}\right)}\,\lesssim\,\left\|\theta_0\right\|_{\dot B^{2-2/s}_{p/2,s}}\,+\,
T^{\gamma_{\rm temp}}\,\Big( 1+ \mathcal N(T)\Big)^2\,\mathcal N^2(T)\,,
\end{equation}
for a (implicit) universal multiplicative constant $C>0$.

This having been established, let us present the proof of Lemma \ref{lem:RHS_MR_theta}.

\begin{proof}[Proof of Lemma \ref{lem:RHS_MR_theta} ]
We consider each term separately. Here below, we fix any time $T\in [0,T_1]$, where $T_1$ is given in \eqref{DefinitionOfT}.

To begin with, we notice that $\theta \in L_T^\infty\left( L^p\right)$ and $\divergence u\in L_T^{r_1}\left( L^p\right)$, the latter owing to \eqref{EstimationP(1)}.
Moreover, since $\rho \in L^\infty_T(L^\infty)$, see \eqref{EstZeta2}, and the function
$\pi_1$ is continuous, one infers that $\pi_1(\rho )\,\theta\, \divergence u\, \in\, L_T^{r_1}\left( L^{p/2}\right)$.
Observe that
\[
\dfrac2{r_1}\,=\, \dfrac2r\, -\, 1\,<\, \dfrac1s\,<\, \dfrac1r\,,
\]
as it follows from \eqref{HypForS}. Thus $1/s \,-\,1/r_1\,>\,1/s\,-\,2/r_1\,=\,2\,\g_{\rm temp}$. This implies, together with \eqref{est:temp_L^p} and
Corollary \ref{CorRegularityu}, the bound
\begin{align*}
\left\Vert \pi_1(\rho)\, \theta\, \divergence u \right\Vert_{L_T^s\left(L^{p/2}\right) }\,&\lesssim\, T^{1/s\,-\,1/r_1}\,
\left\Vert \pi_1(\rho)\right\Vert_{L_T^\infty\left(L^\infty\right)}\,\left\Vert \divergence u \right\Vert_{L_T^{r_1}\left( L^p\right)} \,
\left\Vert \theta \right\Vert_{L_T^\infty\left( L^p\right)} \\ 
& \lesssim\, T^{\gamma_{\rm temp}}\,\Big(1\,+\,\mc N^2(T)\Big)\,\mathcal N^2(T)\,,
\end{align*}
where the implicit multiplicative constat depends on the norm of $\theta_0$ in the space $L^{p/2}\cap L^p$.
Observe that, here, we are implicitly assuming that $T_1\leq1$. This is not restrictive, of course, up to take a smaller value of $T_1$.

Next, as a direct consequence of Corollary \ref{CorRegularityu}, we get the estimate
\[
\left\Vert \mc B(\nabla u)\right\Vert_{L_T^s\left( L^{p/2}\right)}\,\lesssim\, T^{1/s-2/r_1}\,\left\|\nabla u\right\|^2_{L^{r_1}_T(L^p)}\,
\lesssim\,T^{\gamma_{\rm temp}}\,\mathcal N^2(T)\,.
\]

Finally, for the transport term $\rho\,u\cdot\nabla \theta$, it is enough to use the bounds of Lemma \ref{LemmaRegularitytheta} and Corollary
\ref{CorRegularityu}. Observe that, using the definition of $r_{0,\infty}$ given in Lemma \ref{LemmaRegularitytheta}, one sees that
\begin{align*}
\dfrac1s \,-\, \dfrac 1{r_{0, \infty}}\,-\, \dfrac 1{2s}\,=\,\frac{1}{2s}\,-\,\frac{1}{r}\,-\,\frac{d}{2p}\,+\,1\,>\,
\frac{1}{2s}\,-\,\frac{1}{r_1}\,=\,\g_{\rm temp}\,.
\end{align*}
Indeed, a direct computation shows that $1/r_1\,>\,1\,-\,\frac{d}{2p}\,-\,\frac{1}{r}$ if and only if $1/r\,>\,\frac34\,-\,\frac{d}{4p}$. However,
by assumption \eqref{HypForPs} one has $1/r\,>\,1\,-\,\frac{d}{2p}$, so the previous bound is satisfied for free (as $p>d$).
As a conclusion, since $\rho$ is uniformly bounded over $[0,T_1]\times\R^d$, we can estimate
\begin{align*}
\left\|\rho\,u\cdot\nabla\theta\right\|_{L^s_T(L^{p/2})}\,&\lesssim\,T^{1/(2s)\,-\,1/r_{0,\infty}}\,
\left\|u\right\|_{L^{r_0,\infty}_T(L^\infty)}\,\left\|\nabla\theta\right\|_{L^{2s}_T(L^{p/2})} \\
&\lesssim\,T^{\g_{\rm temp}}\,\mc N^2(T)
\end{align*}
for all $T\in[0,T_1]$, where, once again, we have assumed (with no loss of generality) that the condition $T_1\leq1$ is verified.
The lemma is thus completely proven.
\end{proof}

\subsection{Maximal regularity for the modified velocity}\label{SubsectionMaximalRegularityForVelocity}
In this subsection, we establish suitable \tsl{a priori} estimates for the effective velocity field $w$.
For this, we shall take advantage Lemma \ref{LemmaEstimationSolutionElllipticEquation} and the fact that $w$ solves \eqref{EquationForw}:
treating the term $\rho\,\partial_t w$ as done for $\rho\d_t\theta$ at the beginning of Paragraph \ref{sss:MR-theta}, we find
\begin{align*}	
&\Vert w \Vert_{L_T^\infty\left( \dot B_{p, r}^{2-2/r}\right)\cap L^{r_0}(L^{p_0})}\, +\,
\Vert \nabla w \Vert_{L_T^{r_1}\left(L^{p}\right)}\, +\, \Vert \left( \partial_t w, \nabla^2 w \right)\Vert_{L_T^{r}\left(L^{p}\right)}		\\ 
		& 
		\qquad\qquad\qquad
		\lesssim\,\Vert w_0 \Vert_{\dot B^{2-2/r}_{p, r}}\, +\, \Vert \zeta \Vert_{L^\infty_T \left( L^\infty\right)}\,\Vert \partial_t w\Vert_{L_T^r\left( L^p\right)}\,
    +\, \Vert \rho F \Vert_{L_T^{r}\left( L^{p}\right)}\,,
	\end{align*}
	where $T$ has been defined in \eqref{DefinitionF}.
From this inequality, thanks to \eqref{EstZeta2}, up to taking a smaller $\veps_0$ in the statement of Proposition \ref{PropositionEstimatesForModifiedSolution},
we deduce
\begin{align*}	
&\Vert w \Vert_{L_T^\infty\left( \dot B_{p, r}^{2-2/r}\right)\cap L^{r_0}(L^{p_0})}+
\Vert \nabla w \Vert_{L_T^{r_1}\left(L^{p}\right)}+ \Vert \left( \partial_t w, \nabla^2 w \right)\Vert_{L_T^{r}\left(L^{p}\right)}\,
		\lesssim\,\Vert w_0 \Vert_{\dot B^{2-2/r}_{p, r}}+ \Vert F \Vert_{L_T^{r}\left( L^{p}\right)}\,. 
	\end{align*}
Thus, we need to estimate $\Vert F \Vert_{L_T^{r}\left( L^p\right)}$. It turns out that the terms in $G$, introduced in \eqref{defintionF_1},
are the most delicate to bound and will enforce some of the hypotheses of Proposition \ref{PropositionEstimatesForModifiedSolution}.
For this reason, we start by analysing $G$ before studying the rest of the terms in $F$.

\subsubsection{Bounds for $G$}
Recall the definition of $G$, given in \eqref{defintionF_1}:
\[
G\,:=\, 
\mc B(\nabla u)\,-\, \divergence u \, \theta\, \pi_1(\rho) \,+\, \kappa\,\Delta\theta\,.
\]
The idea to bound it is to apply Corollary \ref{ClaimEstimationForOperatorPGeneralp}, for a suitable value of $q$. Now,
the presence of the term $\Delta \theta$ in the definition of
$G$ and the computations performed in Paragraph \ref{sss:MR-theta} force us to take $q=p/2$. Notice that, since $p>d$, one always has $p/2\,\in\,]\oline{p},p[\,$.
Then, applying Corollary \ref{ClaimEstimationForOperatorPGeneralp} and using the uniform boundedness of $\rho$ on $[0,T_1]\times \R^d$ yield,
for any $T\in[0,T_1]$, the bound
\begin{align*}
&\left\Vert \nabla \left( \Id-\Delta\right)^{-1}\left( \frac{\pi_1(\rho)}{\rho} G\right)\right\Vert _{L_T^r\left(L^p\right)} \\
&\qquad\qquad\qquad\lesssim\,
\left\Vert  \frac{\pi_1(\rho)}{\rho} \right\Vert_{L_T^\infty\left( L^\infty\right)}\,\left\Vert G\right\Vert _{L_T^r\left(L^{p/2}\right)} \\
&\qquad\qquad\qquad\lesssim\,T^{1/r\,-\,1/s}\,\left(\left\Vert \Delta \theta \right\Vert _{L_T^s\left(L^{p/2}\right)}\,+\,
\left\Vert \pi_1(\rho)\,\theta\,\div u \right\Vert _{L_T^s\left(L^{p/2}\right)}\,+\,\left\|\mc B(\nabla u)\right\|_{L^s_T(L^{p/2})}\right)\,. 
\end{align*}
Recalling the definition of $\mc N(T)$ and the bounds proved in Lemma \ref{lem:RHS_MR_theta}, we then obtain
\begin{align}
\label{EstF1}
&\left\Vert \nabla \left(\Id -\Delta\right)^{-1}\left( \dfrac{\pi(\rho)}{\rho}\,G\right)\right\Vert _{L_T^r\left( L^p\right)} \\
\nonumber
&\qquad\qquad\qquad \lesssim\,
T^{1/r\,-\,1/s}\,\mathcal N(T)\, +\,T^{1/r\,-\,1/s+\gamma_{\rm temp}}\,\Big( 1+\mathcal N(T)\Big)^2\,\mathcal N^2(T)
\end{align}
for all times $T\in[0,T_1]$.

\subsubsection{Bounds for the rest of the terms in $F$}
We now pass to the estimate of the other terms appearing in the definition of $F$, see formula \eqref{DefinitionF}.
To begin with, recalling that $-\nabla(\Id-\Delta)^{-1}P(\rho,\theta)\,=\,v$, we can use Lemma \ref{LemmaRegularityv} and inequality \eqref{est:1-su-rho}
to write
\begin{equation}\label{MaximalRegularityForWFirstEstimation}
\left\Vert \frac{1}{\rho}\,v\right\Vert _{L_T^r\left( L^p\right)}\,\lesssim\, T^{1/r}\,\mathcal N(T)\,. 
\end{equation}

Next, using again Lemma \ref{LemmaRegularityv} together with Corollary \ref{CorRegularityu}, from Calder\'on-Zygmund theory we have that 
\begin{align}
\label{MaximalRegularityForWSecondEstimation}
\left\Vert \nabla \left( \Id- \Delta \right) ^{-1}\divergence (P\,u)\right\Vert _{L_T^r\left( L^p\right)}\,&\lesssim\,
\left\Vert P\,u\right\Vert _{L_T^r\left( L^p\right)} \\
\nonumber
&\lesssim\,T^{1/r\,-\,1/{r_{0, \infty}}}\,\|u\|_{L^{r_{0,\infty}}_T(L^\infty)}\,\Vert P\Vert_{L_T^\infty\left( L^p\right)}\\ 
&\lesssim\,T^{1/r\,-\,1/{r_{0, \infty}}}\, \mathcal N^2(T)\,,
\nonumber
\end{align} 
where $r_{0, \infty}$ was introduced in Lemma \ref{LemmaRegularityw} and verifies $r_{0,\infty}>r$.

Moreover, similarly to what done for $G$, we use Corollary \ref{ClaimEstimationForOperatorPGeneralp} twice, once with $q=p$ and the other time
with $q=p/2$, to bound
\begin{align*}
&\left\Vert \nabla \left( \Id-\Delta\right)^{-1}\Big( \big( g_0(\rho) + \theta g_1(\rho)\big)\, \divergence u\Big)\right\Vert_{L_T^r\left( L^p\right)} \\
&\qquad \lesssim\,\left\|g_0(\rho)\right\|_{L^\infty_T(L^\infty)}\,\left\Vert \divergence u \right\Vert_{L_T^{r}\left( L^p\right)}\,+\,
\left\|g_1(\rho)\right\|_{L^\infty_T(L^\infty)}\,\left\|\theta\right\|_{L^\infty_T(L^p)}\,\left\Vert \divergence u \right\Vert_{L_T^{r}\left( L^p\right)}\,.
\end{align*}
From the continuity of the functions $g_0$ and $g_1$ and the bounds provided by Corollary \ref{CorRegularityu}, in turn we infer that
\begin{align}
&\left\Vert \nabla \left( \Id-\Delta\right)^{-1}\Big( \big( g_0(\rho) + \theta g_1(\rho)\big)\, \divergence u\Big)\right\Vert_{L_T^r\left( L^p\right)}
\,\lesssim T^{1/2 }\,\mc N(T)\,\Big(1\,+\,\mc N(T)\Big)\,.
\label{MaximalRegularityForWThirdEstimation}
\end{align}
where we have also used the fact that $1/r\,-\,1/r_1\,=\,1/2$, as it directly follows from the definitions.

Finally, we deal with the remaining terms
\[
(w\cdot \nabla) w\, +\, (w\cdot \nabla) v\, +\, (v\cdot \nabla) w\,+\,(v\cdot \nabla) v\,=\,(u\cdot\nabla)u\,. 
\]
After remarking that
\[
\frac{1}{r}\,-\,\frac{1}{r_{0,\infty}}\,-\,\frac{1}{r_1}\,=\,\frac12\,-\,\frac1{r_{0, \infty}}\,=\,
\frac 12 \,-\, \frac{d}{2p}\,+\,1\,-\,\frac{1}{r}\,>\,0\,,
\]
we can apply the estimates of Corollary \ref{CorRegularityu} and get
\begin{align}
\label{MaximalRegularityForWFourthEstimation}
\Vert (u \cdot \nabla) u \Vert_{L_T^r\left( L^p\right)}\,&\lesssim\,
T^{1/2\,-\,1/r_{0,\infty}}\,\left\|u\right\|_{L^{r_{0,\infty}}_T(L^\infty)}\,\left\|\nabla u\right\|_{L^{r_1}_T(L^p)} \\
&\lesssim\,T^{1/2\,-\,1/r_{0,\infty}}\,\mc N^2(T)\,.
\nonumber
\end{align}

\subsubsection{Final bounds for the effective velocity}
It is now time to gather inequalities \eqref{MaximalRegularityForWFirstEstimation}, \eqref{MaximalRegularityForWSecondEstimation}, \eqref{MaximalRegularityForWThirdEstimation}, \eqref{MaximalRegularityForWFourthEstimation} and \eqref{EstF1}. Assuming, as done above, the additional condition $T_1\leq1$,
for any time $T\in[0,T_1]$ we find
\begin{align*}
 \left\|F\right\|_{L^r_T(L^p)}\,&\lesssim\,\left(T^{\frac{1}{r}-\frac{1}{s}}\,+\,T^{\frac{1}{2}}\right)\,\mc N(T)\,\Big(1\,+\,\mc N(T)\Big) \\
 &\qquad\qquad\qquad +\,
\left(T^{\frac{1}{r}-\frac{1}{s}+\g_{\rm temp}}\,+\,T^{\frac{1}{2}-\frac{1}{r_{0,\infty}}}\right)\,\mc N^2(T)\,\Big(1\,+\,\mc N(T)\Big)^2\,,
\end{align*}
where we have also used the fact that $r<2$. From conditions \eqref{HypForPs} and \eqref{HypForS}, since $p>d$, we see that
\[
 \frac{1}{r}\,-\,\frac{1}{s}\,<\,1\,-\,\frac{1}{r}\,<\,\frac{d}{2p}\,<\,\frac{1}{2}\,.
\]
Next, we remark that we also have
\[
 \frac{1}{r}\,-\,\frac{1}{s}\,<\,\frac{1}{2}\,-\,\frac{1}{r_{0,\infty}}\,=\,\frac{3}{2}\,-\,\frac{1}{r}\,-\,\frac{d}{2p}\,,
\]
since $2/r\,-\,1/s<1<3/2\,-\,d/(2p)$. 
Therefore, after setting
\[
 \g_{\rm vel}\,:=\,\frac{1}{r}\,-\,\frac{1}{s}\,,
\]
from the previous bound we finally obtain, for any $T\in[0,T_1]$, the estimate
\begin{align}	
\nonumber
&\Vert w \Vert_{L_T^\infty\left( \dot B_{p, r}^{2-2/r}\right)\cap L^{r_0}(L^{p_0})}+
\Vert \nabla w \Vert_{L_T^{r_1}\left(L^{p}\right)}+ \Vert \left( \partial_t w, \nabla^2 w \right)\Vert_{L_T^{r}\left(L^{p}\right)} \\
&\qquad \qquad \qquad \qquad \qquad \qquad \qquad 
\lesssim\,\Vert w_0 \Vert_{\dot B^{2-2/r}_{p, r}}\,+\,T^{\g_{\rm vel}}\,\Big(\mc N(T)\,+\,\mc N^4(T)\Big)\,. 
\label{FEstimationForw}
	\end{align}

	\medbreak
Let us observe that, in order to complete the \tsl{a priori} estimates for the quantity $\mc N(T)$ defined in \eqref{DefinitionOfN},
we still need to bound the $L^\s_T(L^\infty)$ norm of $\nabla^2(\Id-\Delta)^{-1}P$, with $\s$ introduced in \eqref{ConditionSigma}. Actually, in the proof
of Corollary \ref{CorRegularityu} we have used the stronger property
\begin{equation}\label{TesisP}
\nabla v=\nabla^2(\Id-\Delta)^{-1}P\,\in\,L^{s_1}_T(L^\infty)\,, \qquad\qquad \mbox{ with }\qquad  \dfrac1{s_1}\,:= \,\dfrac1{2\,s}\,\left(1\,+\, \dfrac d{p}\right)\,.
\end{equation}
Recall that $s_1>\s$, see also Remark \ref{r:s_0-s_1} in this respect. Property \eqref{TesisP} was crucially used in Paragraph \ref{sss:Leb-temp},
when bounding $\mc B(\nabla u)$ in $L^{s_2}_T(L^p)$.

The difficulty in getting \eqref{TesisP} relies on the fact that $P$ is discontinuous, as the density $\rho$ is.
Inspired by the pioneering work of Chemin \cite{Ch1991, Ch1993}, in order to prove it we resort to an
approach based on tangential regularity, which (differently to what happened in the isentropic case \cite{DFP}) is then needed not only for uniqueness,
but also at the level of the proof of existence of solutions.

\subsection{Propagation of tangential regularity}\label{SecTanReg}
Consider the non-degenerate family
\[
\mathcal X_0\,=\, \big( X_{0, \lambda}\big)_{1\le \lambda\le m}	
\]
of vector fields in $\mathbb L^{\infty, p}$ given in the assumptions of Proposition \ref{PropositionEstimatesForModifiedSolution}.
For any $\lambda\in\{1,\ldots m\}$, we define the vector field $X_\lambda$ as the unique solution to the transport problem
\begin{equation}\label{EquationForX}
	\begin{cases}
	\partial_t X_\lambda\, +\, (u\cdot \nabla) X_\lambda\,=\, \partial_{X_\lambda} u	\\[1ex]
	{X_\lambda}_{|t=0}\,=\,X_{0,\lambda}\,. 
	\end{cases}
\end{equation}	
We then define, for all $t\ge0$, the new family 
$\mathcal X(t)\, :=\, \big(X_{\lambda}(t)\big)_{1\le \lambda\le m}$.

For convenience of notation, throughout this subsection we set
\[
U(t)\,:=\,\int_0^t\Vert \nabla u(\tau)\Vert_{L^\infty}\,d\tau\,. 
\]
In addition, we reinforce the definition of $T_1$ given in \eqref{DefinitionOfT}. More precisely, we define the time $T_2>0$ as
\begin{align}
\label{def:T_2}
T_2\,&:=\,\sup \left\{ t\,>\,0\quad \Big|\quad U(t)\,\leq\,\log 2 \right. \\ 
\nonumber
& \hspace{3.5cm} \left. \mbox{ and }\quad \int_0^t \Big( \Vert \div w(\tau) \Vert_{L^\infty}\,+\,\Vert \div v(\tau)\Vert_{L^\infty}\Big)\,d\tau\,\leq\,
\left\|\z_0\right\|_{L^\infty}\right\}\,.
\end{align}
Directly from the respective definitions, one easily sees that $0<T_2\leq T_1$. We are now going to perform estimates on the interval
$[0,T_2]$. We are also going to freely use that $T_2\leq 1$, as we has already done with $T_1$ in the previous computations.
%

\subsubsection{Bounds for the vector fields $X_\lambda$}
For simplicity, let the vector field $X$ be any vector field of the family $\mathcal X(t)$. Particularly, we know that $X$ satisfies \eqref{EquationForX}
with an initial datum $X_0\,\in\,\mc X_0$. Then, since $\d_Xu\,=\,(X\cdot\nabla)u$, one easily finds
\begin{equation}\label{EstimationLInftyNormOfX}
\forall\,t\in[0,T_2]\,,\qquad\qquad	\Vert X(t)\Vert_{L^\infty}\,\leq\, \Vert X_0 \Vert_{L^\infty}\,e^{U(t)}\,.
\end{equation}	
In addition, arguing exactly as in Subsection 4.1 of \cite{DFP} (see also Proposition 4.1 of \cite{D1999} for details), one can prove that
\begin{equation}\label{est:I}
 \forall\,t\in[0,T_2]\,,\qquad\qquad I\big(\mc X(t)\big)\,\geq\,I\big(\mc X_0\big)\,e^{-U(t)}\,.
\end{equation}
This shows that the family $\mc X(t)$ remains non-degenerate for any time $t\geq0$, provided the quantity $U(t)$ is finite.

Next, differentiating \eqref{EquationForX} with respect to $x_j$ yields an equation for $\d_jX$, namely
\begin{equation} \label{eq:d_jX}
\partial_t \partial_j X \,+\, (u \cdot \nabla) \partial_j X\, =\, -\,\partial_j u \cdot \nabla X	\,+\, \partial_j \partial_X u\,. 
\end{equation}
By direct computations, we have that
\[
	\partial_j \partial_Xu \,=\, \partial_X \partial_j u\, +\, (\partial_jX \cdot \nabla) u\,.
\]
Hence, performing $L^p$ estimates for equation \eqref{eq:d_jX} and using the Gr\"onwall inequality, we obtain the following bound:
\begin{equation}\label{estimateForNablaX}
\forall\,t\in[0,T_2]\,,\qquad
\Vert \nabla X(t)\Vert_{L^p}\,\le\, \left( \Vert \nabla X_0\Vert_{L^p} + \int_0^t\Vert \partial_X \nabla u (\tau)\Vert_{L^p}\,d\tau\right)\,e^{C\,U(t)}\,,
\end{equation}
for a universal constant $C>0$, independent of the data and the solution.

\medbreak
At this point, for $T\in[0,T_2]$ we define the quantity
\begin{equation} \label{def:Q}
 \mc Q(T)\,:=\,|||\mc X|||_{L^\infty_T(L^\infty)}\,+\,|||\nabla\mc X |||_{L^\infty_T(L^p)}\,+\,\sup_{t\in[0,T]}\frac{1}{I\big(\mc X(t)\big)}\,.
\end{equation}
We also naturally denote $\mc Q_0\,:=\,|||\mc X_0|||_{\mbb L^{\infty,p}}\,+\,1/I(\mc X_0)$.

Then, it follows from \eqref{EstimationLInftyNormOfX}, \eqref{est:I} and \eqref{estimateForNablaX} that
\begin{align} \label{est:Q_1}
\forall\,T\in[0,T_2]\,,\qquad\qquad \mc Q(T)\,\lesssim\,\mc Q_0\,+\,|||\partial_{\mc X} \nabla u|||_{L^1_T(L^p)}\,,
\end{align}
where we have also used the definition of $T_2$ given in \eqref{def:T_2}.

\subsubsection{Tangential bounds for the pressure}
Before estimating the term 
$|||\partial_{\mc X} \nabla u|||_{L^1_T(L^p)}$ appearing in \eqref{est:Q_1}, 
it is convenient to establish a bound for the quantity $\d_{\mc X}P$.

Let $X$ be any vector field of the family $\mc X(t)$.
Directly from the definition of the pressure function $P$, see \eqref{DefinitionP}, we can write
\begin{align*}
 \d_XP\,&=\,\Big(\pi_0'(\rho)\,+\,\pi_1'(\rho)\,\theta\Big)\,\d_X\rho\,+\,\pi_1(\rho)\,\d_X\theta \\
 &=\,\Big(\pi_0'(\rho)\,+\,\pi_1'(\rho)\,\theta\Big)\,\div\big(\rho\,X\big)\,-\,\Big(\pi_0'(\rho)\,+\,\pi_1'(\rho)\,\theta\Big)\,\rho\,\div X
\,+\,\pi_1(\rho)\,\d_X\theta\,,
\end{align*}
where we used that $\d_X\rho\,=\,\div\big(\rho\,X\big)\,-\,\rho\,\div X$, since $\rho$ is discontinuous.
Using the estimates of Lemma \ref{LemmaRegularitytheta} for $\theta$ and $\nabla\theta$, together with the fact that
$s_0>s_1$ and the fact that $\rho\in L^\infty_T(L^\infty)$, we then obtain, for any $T\in[0,T_2]$, 
the bound
\begin{align*}
\left\|\d_XP\right\|_{L^{s_1}_T(L^p)}\,&\lesssim\,T^{1/s_1\,-\,1/s_0}\,\Big(1\,+\,\mc N(T)\Big)\,\left\|\div\big(\rho\,X\big)\right\|_{L^\infty_T(L^p)} \\
&\qquad +\,T^{1/s_1\,-\,1/s_0}\,\Big(1\,+\,\mc N(T)\Big)\,\left\|\nabla X\right\|_{L^\infty_T(L^p)}\,+\,\left\|X\right\|_{L^\infty_T(L^\infty)}\,\mc N(T)\,.
\end{align*}

Next, we compute the equation satisfied by $\Vert \divergence (\rho X)\Vert_{L^p}$. 
Repeating the manipulations of \cite{DFP} (see Paragraph 4A2 therein), we arrive at
\[
\partial_t \divergence \big(\rho\, X\big)\, +\, u\cdot \nabla \divergence \big( \rho\, X\big)\, =\,-\,\div u \; \divergence \big( \rho\, X\big)\,. 
\]
Performing a standard $L^p$ estimate for this transport equation and using definition \eqref{def:T_2} of the time $T_2$, we conclude that
\begin{equation}\label{EstimateForDivergenceRhoX}
\forall\,t\in[0,T_2]\,,\qquad\qquad
\Vert \divergence \big(\rho(t)\,X(t)\big)\Vert_{L^p}\,\lesssim\,\Vert \divergence \big(\rho_0\,X_0\big)\Vert_{L^p}\,. 
\end{equation}

Inserting this bound into the previous estimate for $\d_XP$, we find the following inequality:
\begin{align*}
\forall\,T\in[0,T_2]\,,\quad 
\left\|\d_XP\right\|_{L^{s_1}_T(L^p)}\,&\lesssim\,\left\|X\right\|_{L^\infty_T(L^\infty)}\,\mc N(T) \\
&\qquad \,+\,T^{1/s_1\,-\,1/s_0}\,\Big(1\,+\,\mc N(T)\Big)\,\Big(1\,+\,\left\|\nabla X\right\|_{L^\infty_T(L^p)}\Big)\,.
\end{align*}
Passing to the $\sup$ for all the vector fields $X\in\mc X(t)$ and recalling the definition of $\mc Q$ given in \eqref{def:Q}, in turn we get 
\begin{align}
\label{est:d_XP}
\forall\,T\in[0,T_2]\,,\quad 
|||\d_{\mc X}P|||_{L^{s_1}_T(L^p)}\,&\lesssim\,\mc Q(T)\,\mc N(T)
\,+\,T^{\g_s}\,\Big(1\,+\,\mc N(T)\Big)\,\Big(1\,+\,\mc Q(T)\Big)\,,
\end{align}
where we have defined
\[
 \g_s\,:=\,\frac{1}{s_1}\,-\,\frac{1}{s_0}\,=\,\frac{1}{2s}\,\left(1\,-\,\frac{d^2}{p^2}\right)\,>\,0\,.
\]

\subsubsection{Tangential regularity of $\nabla u$}
It is now time to establish an estimate for the $L^1_T(L^p)$ norm of $\d_{\mc X}\nabla u$ appearing in \eqref{est:Q_1}.

Again, let $X$ be any vector field of the family $\mc X(t)$.
By writing $u\,=\,w+v$, we infer that
\begin{equation} \label{eq:d_XDu}
\partial_X \nabla u\, =\, \partial_X \nabla w\, +\, \partial_X\nabla v\,,\qquad\qquad \mbox{ with }\qquad \nabla v\,=\,-\,\nabla^2\left( \Id- \Delta\right)^{-1}P\,, 
\end{equation}
Observe that we have, for any $T\in[0,T_2]$, the bound
\begin{equation*} 
\Vert \partial_X \nabla w \Vert_{L^1_T(L^p)}\,\leq\,T^{1-1/r}\,\Vert X\Vert_{L^\infty_T(L^\infty)}\,\Vert \nabla^2 w\Vert_{L^r_T(L^p)}\,\lesssim\,
T^{1-1/r}\,\left\|X\right\|_{L^\infty_T(L^\infty)}\,\mc N(T)\,,
\end{equation*}
whence we infer that
\begin{equation} \label{est:d_XDw}
\forall\,T\in[0,T_2]\,,\qquad |||\d_{\mc X}\nabla w|||_{L^1_T(L^p)}\,\lesssim\,T^{1-1/r}\,\mc Q(T)\,\mc N(T)\,.
\end{equation}

The second term on the right-hand side of \eqref{eq:d_XDu} is more complex to bound.
We profit from a tangential regularity estimate established in \cite{DFP} (see Propositions 2.12 and 2.13 therein).

\begin{prop}
	\label{PropositionEstimationForTangentialRegularity}
Fix $p \in \,]d, +\infty[\,$ and an integer $m\ge d-1$. Take a non-degenerate family $\mathcal X:= (X_\lambda)_{1\le \lambda\le m}$ of $m$ vector fields belonging
to $\mathbb L^{\infty, p}$. Let $g \in L^\infty\left(\R^d \right)$ be such that $g \in \mathbb L^p_{\mathcal X}$.

Then, for any $\eta>0$, one has the following properties:
\[
\nabla^2\left(\eta\Id-\Delta\right)^{-1} g\, \in\, L^\infty\qquad\qquad \mbox{ and }\qquad\qquad
\nabla^2\left(\eta\Id-\Delta\right)^{-1} g\,\in\,\mbb L^p_{\mc X}\,.
\]
Morover, there exists $C>0$ such that, after defining
\[
 A(\mc X)\,:=\,1+ \frac{|||\mathcal X|||^{4d-5}_{L^\infty}|||\nabla \mathcal X|||_{L^p}}{\big( I(\mathcal X)\big)^{4d-4}}
\qquad \mbox{ and }\qquad
B(\mc X)\,:=\,\frac{|||\mathcal X|||_{L^\infty}^{4d-5}}{\big( I(\mathcal X)\big)^{4d-4}}\,,
\]
one has the following estimates:
\begin{align*}
\left\|\nabla^2\left(\eta\Id-\Delta\right)^{-1} g\right\|_{L^\infty}\,&\lesssim\,
A(\mc X)\,\Vert g \Vert_{L^\infty}\, +\, B(\mc X)\, |||\partial_{\mathcal X}g |||_{L^p}\,, \\
|||\d_{\mc X}\nabla^2\left(\eta\Id-\Delta\right)^{-1} g |||_{L^p}\,&\lesssim\,
|||\nabla \mathcal X|||_{L^p}\,A(\mc X)\,\Vert g \Vert_{L^\infty} \\
&\qquad
+\,B(\mc X)\,|||\mathcal X|||_{L^\infty}\,|||\nabla \mathcal X|||_{L^p}\,\Vert g \Vert_{L^\infty}\,+\,A(\mc X)\,|||\partial_{\mathcal X}g |||_{L^p}\,.
\end{align*}
\end{prop}

We now use the previous proposition with $g=P$ and $\mc X=\mc X(t)$ the given non-degenerate family of vector fields.
We notice that
\begin{equation} \label{est:A-B}
 A\big(\mc X(t)\big)\,\lesssim\,1\,+\,\mc Q^{8d-8}(t)\qquad \mbox{ and }\qquad
 B\big(\mc X(t)\big)\,\lesssim\,\mc Q^{8d-9}(t)
\end{equation}
Thus, the last estimate of Proposition \ref{PropositionEstimationForTangentialRegularity} implies that, for any time $T\in[0,T_2]$, one has
\begin{align*}
|||\d_{\mc X}\nabla v |||_{L^1_T(L^p)}\,&\lesssim\,\left(1\,+\,\mc Q^{8d-7}(T)\right)\,\left\|P\right\|_{L^1_T(L^\infty)}\,+\,
\left(1\,+\,\mc Q^{8d-8}(T)\right)\,|||\d_{\mc X}P|||_{L^1_T(L^p)} \\
&\lesssim\,T^{1-1/s_0}\,\left(1\,+\,\mc Q^{8d-7}(T)\right)\,\Big(1\,+\,\mc N(T)\Big) \\
&\qquad\qquad\qquad\qquad\qquad\qquad\qquad
\,+\,T^{1-1/s_1}\,\left(1\,+\,\mc Q^{8d-7}(T)\right)\,\mc N(T)\,, 
\end{align*}
where, in passing from the first inequality to the second one, we have also used the bounds from Lemma \ref{LemmaRegularityv} and \eqref{est:d_XP}.

At this point, we can gather this latter estimate and \eqref{est:d_XDw} together: owing to the decomposition \eqref{eq:d_XDu}, we find
\begin{align}
\label{est:d_XDu}
\forall\,T\in[0,T_2]\,,\qquad
|||\d_{\mc X}\nabla u |||_{L^1_T(L^p)}\,&\lesssim\,T^{1-1/r}\,\left(1\,+\,\mc Q^{8d-7}(T)\right)\,\Big(1\,+\,\mc N(T)\Big)
\end{align}
Here above, we have also used the fact that $s_0>s_1>s>r$ (as it follows from assumptions \eqref{HypForS} and Remark \ref{r:s_0-s_1}),
so that $1-1/s_0\,>\,1-1/s_1\,>\,1-1/r$, together with the condition $T_2\leq 1$ we have implicitly taken above,
see the comments below \eqref{def:T_2}.

\subsubsection{Bounding $\nabla v$ in $L^{s_1}_T(L^\infty)$}
As mentioned above, in order to close the bootstrap argument a key step is to prove \eqref{TesisP}.
This bound follows from Proposition \ref{PropositionEstimationForTangentialRegularity}:
using the first estimate appearing therein and \eqref{est:A-B}, we deduce, for any $T\in[0,T_2]$, the inequality
\begin{align*}
\left\|\nabla v\right\|_{L^{s_1}_T(L^\infty)}\,\lesssim\,
T^{\g_s}\,\left(1\,+\,\mc Q^{8d-8}(T)\right)\,\left\|P\right\|_{L^{s_0}_T(L^\infty)}\,+\,
\mc Q^{8d-9}(T)\,|||\d_{\mc X}P|||_{L^{s_1}_T(L^p)}\,,
\end{align*}
where the exponent $\g_s>0$ is the same appearing in \eqref{est:d_XP}.
Using Lemma \ref{LemmaRegularityv} and estimate \eqref{est:d_XP} again, we finally arrive at the following inequality:
\begin{align}
\label{est:Dv-inf}
\forall\,T\in[0,T_2]\,,\qquad
\left\|\nabla v\right\|_{L^{s_1}_T(L^\infty)}\,&\lesssim\,
T^{\g_s}\,\left(1\,+\,\mc Q^{8d-8}(T)\right)\,\Big(1\,+\,\mc N(T)\Big) \\
\nonumber
&\qquad\qquad\qquad\qquad\qquad\qquad
\,+\,\mc Q^{8d-8}(T)\,\mc N(T)\,.
\end{align}

\subsection{Closing the bootstrap argument} \label{ss:close}
Thanks to the estimates established in Subsections \ref{SubsectionAPrioriEstimatesForDensity}, \ref{ss:apriori-temp}, \ref{SubsectionMaximalRegularityForVelocity}
and \ref{SecTanReg}, we can now close the bootstrap argument and establish a uniform bound for the quantities
$\mc N(T)$ and $\mc Q(T)$, defined respectively in \eqref{DefinitionOfN} and \eqref{def:Q}, in some small time interval $[0,T^*]$.

Recall that the series of inequalities we have established so far hold true on the time interval $[0,T_2]$, with $T_2>0$ defined in \eqref{def:T_2}.
In this time interval, let us define
\[
 \mc E(T)\,:=\,\mc N(T)\,+\,\mc Q(T)\qquad \mbox{ and }\qquad \mc E_0\,:=\,\mc N_0\,+\,\mc Q_0\,.
\]
Thus, our goal is to find a uniform bound for $\mc E$ on the (possibly) smaller time interval $[0,T^*]$.

First of all, we insert inequality \eqref{est:d_XDu} into \eqref{est:Q_1} to obtain
\begin{equation} \label{est:Q}
\forall\,T\in[0,T_2]\,,\qquad
\mc Q(T)\,\lesssim\,\mc Q_0\,+\,T^{1-1/r}\,\Big(1\,+\,\mc E(T)\Big)^{8d-6}\,.
\end{equation}
Next, from \eqref{est:Dv-inf} we see that
\begin{align*}
 \left\|\nabla v\right\|_{L^{\s}_T(L^\infty)}\,&\lesssim\,T^{1/\s\,-\,1/s_1}\,\left(1\,+\,\mc Q^{8d-8}(T)\right)\,\Big(1\,+\,\mc N(T)\Big) \\
 &\lesssim\,T^{1/\s\,-\,1/s_1}\,\Big(1\,+\,\mc E(T)\Big)^{8d-7}\,,
\end{align*}
where we recall that $1/\s\,-\,1/s_1\,>\,0$, according to \eqref{ConditionSigma} and Remark \ref{r:s_0-s_1}.
Summing up this latter inequality with estimates \eqref{EstZeta2}, \eqref{EstLp_Zeta}, \eqref{est:temp_L^p}, \eqref{EstimationForTheta}
and \eqref{FEstimationForw}, and keeping Remark \ref{r:norm-temp} in mind, we finally find
\begin{align} \label{est:N}
\forall\,T\in[0,T_2]\,,\qquad
\mc N(T)\,\lesssim\,\mc N_0\,+\,T^\g\,\Big(1\,+\,\mc E(T)\Big)^{8d-7}
\end{align}
where we have used the condition $T_2\leq1$ and the fact that, for $d\geq2$, one has $8d-7\geq 4$, and where and we have set
\[
 \g\,:=\,\min\Big\{\frac{1}{\s}-\frac{1}{s_1}\,,\;1-\frac{1}{r}\,,\;\g_{r,s}\,,\;\g_{\rm temp}\,,\;\g_{\rm vel}\Big\}\,.
\]
At this point, putting inequalities \eqref{est:Q} and \eqref{est:N} together, up to redefine the value of $\g>0$ and for a suitable exponent $k\geq8d-7$,
we find the bound
\begin{align} \label{est:E}
\forall\,T\in[0,T_2]\,,\qquad
\mc E(T)\,\lesssim\,\mc E_0\,+\,T^\g\,\Big(1\,+\,\mc E(T)\Big)^k
\end{align}
As already claimed above, based on estimate \eqref{est:E}, it is a routine matter to find a suitable time $T^*>0$, only depending
on $\mc E_0$, and a universal constant $K_0>0$ such that
\begin{align} \label{est:E-final}
\sup_{T\in[0,T^*]}\mc E(T)\,\leq\,K_0\,\mc E_0\,.
\end{align}

\section{Proof of existence} \label{s:existence}

In this section, we give the details of the proof of existence of a solution, as claimed in Theorem \ref{th:well-p}. The same argument trivially
implies the claim of Proposition \ref{PropositionEstimatesForModifiedSolution}, as actually the proofs of the two statements go in parallel.

We divide our discussion into two cases: in the first one, treated in Subsection \ref{ss:ex-energy}, we assume the initial data to have finite energy;
the general case will instead be considered in Subsection \ref{ss:ex-general}.
Recall that, in Theorem \ref{th:well-p}, we have fixed the value of the reference density $\oline\rho$ to be equal to $1$, for notational
simplicity.

\subsection{The case of finite energy initial data} \label{ss:ex-energy}

In this subsection, we are going to suppose that, besides the hypotheses fixed in Theorem \ref{th:well-p}, the initial datum
$\big(\rho_0,u_0,\theta_0\big)$ enjoys the extra assumption
\begin{equation} \label{hyp:energy}
\rho_0\,-\,1\,,\; u_0\,,\; \theta_0\quad \in\; L^2(\R^d)\,.
\end{equation}
According to \eqref{definitione}, we define the quantities
\[
 e_0\,:=\,c_v\,\theta_0\,+\,h(\rho_0)\qquad \mbox{ and }\qquad \oline e\,:=\,h(1)\,.
\]
Observe that, up to adding a suitable constant when solving the differential equation \eqref{RelationBetweenPih}, we can assume without any loss of
generality that
\[
 e_0\,\geq\,e_*\,>\,0\,,
\]
for some positive constant $e_*$. In addition, it follows from \eqref{hyp:energy} that
\[
 e_0\,-\,\oline e\,\in\,L^2(\R^d)\,.
\]

Next, we regularise the initial datum \tsl{e.g.} by use of the low frequency cut-off operators $S_m$ introduced in Subsection \ref{ss:LP}. We thus define,
for $m\in\N$, the functions
\[
 \rho_{0,m}\,:=\,S_m\rho_0\,,\qquad u_{0,m}\,:=\,S_mu_0\,,\qquad \theta_{0,m}\,:=\,S_m\theta_0\,.
\]
Of course, we make a similar regularisation procedure also for each element of the non-degenerate family of vector fields $\mc X_0$,
thus producing a regularised family $\mc X_{0,m}$.

By arguing as in the proof of Lemma 5.1 of \cite{Hoff_1997}, we can construct a sequence $\big(\rho_m,u_m,\theta_m\big)_{m\in\N}$ of smooth local solutions
to system \eqref{eq:NSF}, each one defined on some small time interval $[0,T_m]$ and having $H^\infty$ regularity in space.
Since these solutions decay fast enough at spatial infinity, we can rigorously perform the \tsl{a priori} estimates from Subsections \ref{ss:static},
\ref{SubsectionAPrioriEstimatesForDensity}, \ref{ss:apriori-temp}, \ref{SubsectionMaximalRegularityForVelocity}, \ref{SecTanReg} and \ref{ss:close}
for them. More precisely, introducing the quantity $\nu\,:=\,\lambda\,+\,2\,\mu$, we define
\begin{align*}
&\z_m\,:=\,\rho_m\,-\,1 \\
& \qquad \mbox{ and }\qquad w_m\,:=\,u_m\,-\,v_m\,,\quad \mbox{ with } \ v_m\,:=\,-\,\frac{1}{\nu}\,\nabla(\Id-\Delta)^{-1}P(\rho_m,\theta_m)\,.
\end{align*}
We also set, for any time $T>0$, the quantities $\mc N_m(T)$ and $\mc Q_m(T)$ to be, respectively, the functions $\mc N$, defined in \eqref{DefinitionOfN},
and $\mc Q$, defined in \eqref{def:Q}, but computed on the smooth $m$-th solution:
\begin{align*}
	\mathcal N_m(T)\,&:=\,\Vert w_m \Vert_{L_T^\infty \left(\dot B^{2-2/r}_{p,r}\right)\cap L_T^{r_0}\left( L^{p_0}\right)}\,+\,
	\Vert \nabla w_m \Vert_{L_T^{r_1}\left( L^{p}\right)}\,+\,  \Vert \left( \partial_t w_m\,,\, \nabla^2 w_m\right) \Vert_{L_T^{r}\left( L^{p}\right)} \\ 
 	&\qquad\qquad\qquad\qquad +\,  \Vert  \zeta_m \Vert_{L_T^{\infty}\left( L^{p}\cap L^\infty\right)} \, +\,
 \Vert \nabla^2 \left( \Id-\Delta\right)^{-1}P(\rho_m,\theta_m)\Vert_{L_T^{\sigma}\left( L^\infty\right)}  \\
&\qquad\qquad\qquad\qquad\qquad +\,\Vert \theta_m \Vert_{L_T^{\infty}\left( B^{2-2/s}_{p/2,s}\cap L^{p}\right)} \, +\,
 \left\Vert \left( \partial_t \theta_m\,,\, \nabla^2\theta_m\right)\right\Vert_{L_T^s\left( L^{p/2}\right)}\,, \\
 \mc Q_m(T)\,&:=\,|||\mc X_m|||_{L^\infty_T(L^\infty)}\,+\,|||\nabla\mc X_m |||_{L^\infty_T(L^p)}\,+\,\sup_{t\in[0,T]}\frac{1}{I\big(\mc X_m(t)\big)}\,.
\end{align*}
For later use, we also set
\begin{align*}
& \mc N_{0,m}\,:=\,
 \Vert w_{0,m} \Vert_{\dot B^{2-2/r}_{p,r}}\, +\,  \Vert  \zeta_{0,m} \Vert_{L^{p}\cap L^\infty}\, +\,  \Vert \theta_{0,m} \Vert_{B^{2-2/s}_{p/2,s}\cap L^{p}} \\
&\qquad\qquad\qquad\qquad \mbox{ and }\qquad\qquad
\mc Q_{0,m}\,:=\,|||\mc X_{0,m}|||_{\mbb L^{\infty,p}}\,+\,\frac{1}{I(\mc X_{0,m})}\,.
\end{align*}
Observe that, by standard properties of the low frequency cut-off operators, one has
\[
 \mc N_{0,m}\,\lesssim\,
 \mc N_0\,:=\,\Vert w_{0} \Vert_{\dot B^{2-2/r}_{p,r}}\, +\,  \Vert  \zeta_{0} \Vert_{L^{p}\cap L^\infty}\, +\,  \Vert \theta_{0} \Vert_{B^{2-2/s}_{p/2,s}\cap L^{p}}\,,
\]
for an implicit multiplicative constant which is independent of $m\in\N$. It goes without saying that a similar bound holds true also
for the quantity $\mc Q_{0,m}$.

Then, repeating the computations of Section \ref{s:a-priori}, we deduce that there exists a time $T>0$ such that,
after defining $\mc E_m(T)\,:=\,\mc N_m(T)+\mc Q_m(T)$ and $\mc E_{0,m}\,:=\,\mc N_{0,m}+\mc Q_{0,m}$, one has
\begin{equation} \label{est:approx_m}
\inf_{m\in\N}T_m\geq T\qquad \mbox{ and }\qquad
 \mc E_m(T)\,\lesssim\,\mc E_{0,m}\,\lesssim\,\mc E_0\,:=\,\mc N_0\,+\,\mc Q_0\,,
\end{equation}
for a new implicit constant $C>0$, independent of the approximation parameter $m\in\N$.
The above inequality implies weak-$*$ compactness of the family $\big(\rho_m,w_m,\theta_m\big)$,
so convergence (up to the extraction of a suitable subsequence) to some target profile $\big(\rho,w,\theta\big)$ such that
the following properties are satisfied:
\begin{itemize}
 \item[(a)] $\z\,:=\,\rho-1\;\in\, L_T^{\infty}\left( L^{p}\cap L^\infty\right)$ verifies  $\Vert \zeta \Vert_{L_T^\infty\left( L^\infty\right)}\,\leq\,4\,\varepsilon$;
\item[(b)] $\theta \in L^\infty_T\left( B^{2-2/s}_{p/2,s}\cap L^{p}\right)$, with $\partial_t \theta, \nabla^2\theta\; \in L_T^s\left( L^{p/2}\right)$;
\item[(c)] $w \in L^\infty_T \left(\dot B^{2-2/r}_{p,r}\right)\cap L_T^{r_0}\left( L^{p_0}\right)$, with $\nabla w \in L_T^{r_1}\left( L^{p}\right)$ and
$\partial_t w, \nabla^2 w\;\in L_T^{r}\left( L^{p}\right)$.
\end{itemize}
Similar weak-$*$ convergences can be stated also for the vector fields $X_m$ of the family $\mc X_m$, giving rise to the family $\mc X$ of vector fields
in $\mbb L^{\infty,p}$.

At this point, strong convergence of $w_m$ and $\theta_m$ can be obtained by a standard argument, based on the control of their time derivatives
(to be obtained through their respective equations)
and an application of the Ascoli-Arzel\`a theorem in suitable spaces.
By the same token, using in addition the corresponding equation \eqref{EquationPressure} for the pressure $P_m:=P(\rho_m,\theta_m)$
and the space compactness of $v_m$, one also gets strong convergence of $v_m$ towards some profile $v$. So far, we only know that
$v\,=\,-\,\frac{1}{\nu}\,\nabla(\Id-\Delta)^{-1}\Phi$, for a suitable potential function $\Phi$ belonging to
$L^{s_0}_T(L^\infty)\cap L^\infty_T(L^p)$, which is the weak-$*$ limit of the pressure functions $P_m$ (keep in mind Lemma \ref{LemmaRegularityv},
which yields uniform bounds for $\big(P_m\big)_{m\in\N}$).

It is now time to use the uniform bounds coming from striated regularity: from them, we get the boundedness of the sequence $\big(u_m\big)_m$
in the space $L^{\wtilde m}_T(W^{1,\infty})$, with $\wtilde m\,:=\,\min\big\{r_{0,\infty}\,,\,m_1\big\}$ (see Corollary \ref{CorRegularityu}).
Thus, setting $u\,:=\,w+v$, we infer the weak-$*$ convergence of
$u_m\,\stackrel{*}{\rightharpoonup}\,u$ in \tsl{e.g.} $L^{\wtilde m}_T(W^{1,\infty})$. In particular, one has
$\nabla u\,\in\,L^{\wtilde m}_T\big(L^\infty\big)$. With this property, together with the compactness of $w_m$ and $v_m$,
we can pass to the limit in equation \eqref{EquationForX} at level $m$ and get the correcponding equation for the limit
vector fields of the family $\mc X(t)$, with $t\in[0,T]$.

Next, let us focus on the density functions $\rho_m$ for a while. One can repeat the argument of Lions \cite{Lions} (see also \cite{HPZ})
to deduce the strong convergence (up to a further extraction, not relabelled here) of $\big(\rho_m\big)_m$ to the target $\rho$ in the space
$\mc C\big([0,T];L^q_{\loc}\big)$, for any $q\in[p,+\infty[\,$.
The consequence of this latter property is two-fold. First of all, we infer the convergence of the pressure functions
$P(\rho_m,\theta_m)\,\longrightarrow\,P(\rho,\theta)$ a.e. on $[0,T]\times\R^d$, in the limit $m\to+\infty$, which yields
in particular the identification of the limit velocity field
\[
 v\,=\,-\,\frac{1}{\nu}\,\nabla(\Id-\Delta)^{-1}P(\rho,\theta)\,.
\]
In addition, this compactness result, together with the above mentioned ones for $w_m$, $\theta_m$ and $v_m$, allows us to pass to the limit
in the weak formulation of the equations, implying thus that the triplet $\big(\rho,u,\theta\big)$, with $u\,=\,w+v$, is indeed a weak solution
to system \eqref{eq:NSF} related to the initial datum $\big(\rho_0,u_0,\theta_0\big)$.

Finally, the claimed time continuity properties of the solution follows from standard results about solutions to continuity and parabolic equations
in $L^q$ spaces, combined with classical maximal regularity results to establish the time continuity with values in suitable homogeneous Besov spaces.
Tangential regularity properties for the target density $\rho$ with respect to the family $\mc X(t)$ can also be obtained by standard arguments.
We omit the precise details here.

This proves the existence statement in the case the initial datum has finite energy, namely under the assumption that
the conditions in \eqref{hyp:energy} are satisfied.

Before going on, let us point out the following important fact.
\begin{rmk} \label{r:L^p}
The previous proof, supplemented with minor modifications,  also yields a similar existence result whenever we replace assumption \eqref{hyp:energy}
with the following conditions:
\begin{equation} \label{hyp:L^p}
\rho_0\,-\,1\,,\; \theta_0\quad \mbox{ belong to }\ L^2(\R^d)\qquad\quad \mbox{ and }\qquad\quad w_0\in L^p(\R^d)\,.
\end{equation}
Indeed, in this case it is enough to argue as follows (see a similar approach in \cite{F-V}, for instance).
Take a smooth compactly supported function $\eta\,\in\,\mc C^\infty_0(\R^d)$ such that
\[
0\,\leq\,\eta\,\leq 1\,,\qquad\qquad \eta\equiv1 \; \mbox{ for } \; |x|\leq1\qquad \mbox{ and }\qquad \eta\equiv0 \; \mbox{ for } \; |x|\geq2\,,
\]
For any $k\in\N\setminus\{0\}$, set
\[
 \eta_k(x)\,:=\,\eta\left(\frac{x}{k}\right)\qquad \mbox{ and }\qquad w^k_0\,:=\,\eta_k\,w_0\,.
\]
Then, since $p>d\geq 2$, it is easy to see that $\big(w^k_0\big)_{k\geq1}$ is a bounded sequence in $L^2(\R^d)$.
Therefore, after defining $u_0^k\,:=\,w^k_0\,+\,v_0$, with $v_0\,:=\,-\,\frac{1}{\nu}\,\nabla(\Id-\Delta)^{-1}P(\rho_0,\theta_0)$, we have that
also $\big(u^k_0\big)_{k\geq1}$ is a bounded sequence in $L^2(\R^d)$, whence the claimed result follows from the above proof, combined with a compactness argument
for taking the limit $k\to+\infty$.
\end{rmk}

\subsection{The general case} \label{ss:ex-general}

We now treat the general case, in which we do not assume hypothesis \eqref{hyp:energy}, nor \eqref{hyp:L^p}, to hold. In this case, we need
another approximation procedure for the velocity fields, instead of the one proposed in Remark \ref{r:L^p}.

To begin with, we take the function $\eta\,\in\,\mc C^\infty_0(\R^d)$ fixed in Remark \ref{r:L^p}
and we define, for $k\in\N\setminus\{0\}$, the rescaled functions $\eta_k$ accordingly:
\[
 \eta_k(x)\,:=\,\eta\left(\frac{x}{k}\right)\,.
\]

Recall that $\zeta_0\,:=\,\rho_0-1$. Then, for any $k\in\N\setminus\{0\}$ we set
\[
 \z_0^k\,:=\,\eta_k\,\z_0\qquad\qquad \mbox{ and }\qquad\qquad 
 \theta^k_0\,:=\,\eta_k\,\theta_0\,.
\]
Let us focus on the density functions for a while. Owing to assumption \eqref{th_HypothesisForZetaPropositionEstimatesForModifiedSolution},
we see that $\z_0^k$ is a bounded function, with $\|\z_0^k\|_{L^\infty}\leq\veps$, having compact support; thus $\z_0^k\,\in \,L^2$.
Then, we can define the approximate initial densities
$\rho_0^k$ as
\[
 \rho_0^k\,:=\,1\,+\,\z_0^k\,.
\]
It follows from the previous arguments that
\begin{equation} \label{est:rho_0^k-B}
\left\|\rho^k_0-1\right\|_{L^p\cap L^\infty}\,\leq\,\left\|\rho_0-1\right\|_{L^p\cap L^\infty}\,.
\end{equation}
Finally, we point out that, by definition of $\eta_k$, one has $\rho_0^k\equiv\rho_0$ on the ball
$B(0,k)$ of center $0$ and radius $k$. 
Thus, we infer the convergence property 
\[ 
 \rho_0^k\,\longrightarrow\,\rho_0\qquad\qquad \mbox{ in }\quad L^\infty_{\rm loc}(\R^d)\,,
\] 
in the limit $k\to+\infty$.
By the same token, we also infer that, for any $k\in\N\setminus\{0\}$ fixed, one has $\theta_0^k\in L^2$, together with
the convergence property, when $k\to+\infty$,
\[ 
 \theta_0^k\,\longrightarrow\,\theta_0\qquad\qquad \mbox{ in }\quad L^p_{\rm loc}(\R^d)\,.
\] 

Next, we propose a suitable approximation of the velocity fields. Consider $w_0\,:=\,u_0\,-\,v_0$, where, as usual, we have set
$v_0\,:=\,-\,\frac{1}{\nu}\,\nabla(\Id-\Delta)^{-1}P(\rho_0,\theta_0)$.
Recall that $w_0\,\in\,\dot B^{2-2/r}_{p,r}$ by assumption.
Then, for any $k\in\N\setminus\{0\}$, we define
\[
 w_0^k\,:=\,\left(\Id-\dot S_{-k}\right)\,w_0\,,
\]
where the operators $\dot S_j$ have been introduced in \eqref{eq:S_j-dot}
We have the following easy result, which however we will prove, for the sake of completeness, at the end of this part.

\begin{lemma} \label{l:approx-w}
For any $k\in\N\setminus\{0\}$, one has the property $w_0^k\in L^p$.
Moreover, the sequence $\big(w_0^k\big)_{k\geq1}$ is uniformly bounded in the space $\dot B^{2-2/r}_{p,r}$. In addition,
in the limit $k\to+\infty$, one has the convergence $w^k_0\,\longrightarrow\,w_0$ in the sense of $\mc S'$.
\end{lemma}

With Lemma \ref{l:approx-w} at hand, we can conclude our proof. To this scope, for any integer $k\geq1$, define
the approximate initial velocity field
\[
 u_0^k\,:=\,w^k_0\,+\,v_0^k\,,\qquad\qquad \mbox{ with }\qquad v_0^k\,:=\,-\,\frac{1}{\nu}\,\nabla(\Id-\Delta)^{-1}P(\rho^k_0,\theta_0^k).
\]
Then, for any $k\in\N\setminus\{0\}$ fixed, from the considerations above we infer that the triplet
$\big(\rho_0^k,u^k_0,\theta_0^k\big)$ satisfies the conditions in \eqref{hyp:L^p}, of course not uniformly with respect to $k$.
Then, repeating the arguments of Subsection \ref{ss:ex-energy} and Remark \ref{r:L^p}, we can construct a sequence
$\big(\rho^k,u^k,\theta^k\big)_{k\geq1}$  of local in time solutions to system \eqref{eq:NSF}, related to those initial data.
By following the computations performed in Section \ref{s:a-priori}, it is not hard to find a uniform (in $k\in\N\setminus\{0\}$) time interval $[0,T]$
such that the sequence $\big(\rho^k,w^k,\theta^k\big)_{k\geq1}$ satisfies the properties (a)-(b)-(c) of Subsection \ref{ss:ex-energy},
\emph{uniformly} with respect to $k\in\N\setminus\{0\}$.
We stress the fact that the time $T>0$ does not depend on the $L^p$ norm of $w_0^k$, so it is uniform for the whole sequence of approximate solutions;
in this respect, recall the discussion before \eqref{est:E-final}.
Finally, it is a matter of repeating the compactness arguments of the previous subsection, combined with the convergence of the initial datum,
to prove convergence of the sequence (up to a suitable extraction) to a solution to equations \eqref{eq:NSF} verifying the properties claimed in Theorem \ref{th:well-p}.

Therefore, in order to complete the proof of the existence part, we only have to show the proof of Lemma \ref{l:approx-w}, which we now tackle.

\begin{proof}[Proof of Lemma \ref{l:approx-w}]
First of all, by definition, it is easy to see that $w_0^k\in \mc S'_h$ for any integer $k\geq1$. In addition, for tempered distributions $f$ in 
$\mc S'_h$, one has the equality
\begin{equation} \label{eq:dot-S}
 \dot S_jf\,=\,\sum_{\ell\leq j-1}\dot \Delta_\ell f
\end{equation}
for any $j\in\Z$.
This immediately implies that
$ \left\|w^k_0\right\|_{\dot B^s_{p,r}}\,\leq\,\left\|w_0\right\|_{\dot B^s_{p,r}}$,
where, for simplicity, we have denoted $s\,=\,2\,-\,2/r$. Thus, $\big(w_0^k\big)_{k\geq1}$ is bounded in $\dot B^s_{p,r}$, as claimed.

Next, by definition, we also see that
\[
 w_0^k\,=\,\sum_{\ell\geq -k}\dot\Delta_\ell w_0\,,
\]
from which we deduce, for a constant $C_k$ which exploses when $k\to+\infty$, that
\begin{align*}
 \left\|w^k_0\right\|_{L^p}\,&\leq\,\sum_{\ell\geq -k}\left\|\dot\Delta_\ell w_0\right\|_{L^p}\,
\leq\,\left\|w_0\right\|_{\dot B^s_{p,r}}\,\left(\sum_{\ell\geq -k} 2^{-\,s\,\frac{r}{r-1}\,\ell}\right)^{1-1/r}\,\leq\,
C_k\,\left\|w_0\right\|_{\dot B^s_{p,r}}\,.
\end{align*}
Therefore, for any $k\in \N\setminus\{0\}$ fixed, one has $w_0^k\in L^p$.

It remains us to show the convergence $w^k_0\,\longrightarrow\,w_0$ in the sense of $\mc S'$, when $k\to+\infty$. By a duality argument, it is enough
to show that, for any $\phi\in \mc S(\R^d)$, one has
\begin{equation} \label{conv:to-show}
\left(\Id\,-\,\dot S_{-k}\right)\phi\,\longrightarrow\,\phi\,,\qquad\qquad \mbox{ that is }\qquad
\dot S_{-k}\phi\,\longrightarrow\,0\,,
\end{equation}
in some suitable topology. 
Observe that, if $f\in L^p$, then $f\in \mc S'_h$. Thus $\dot S_{-k}\phi$ belongs to $\mc S'_h$ and, by definition \eqref{def:S_h},
one has $\dot S_{-k}\phi\,\longrightarrow\,0$ for $k\to+\infty$ in the sense of $\mc S'$. Since $\big(\dot S_{-k}\phi\big)_{k\geq1}$
is a bounded sequence in \tsl{e.g.} $L^p$ for $1<p<+\infty$, by uniqueness of the weak limit we immediately deduce that the convergence \eqref{conv:to-show}
holds true in the weak topology of $L^p$. This in turn implies the sought convergnce of $w^k_0$ towards $w_0$ in the sense of $\mc S'$.

The lemma is now completely proven.
\end{proof}

\section{Uniqueness} \label{s:uniqueness}

This section deals with the proof of uniqueness of solutions claimed in Theorem \ref{th:well-p}.
In order to deal with the low regularity of the density, we follow a nowadays standard approach \cite{Hoff_2006, D-2014, Chik-D, DFP}
and reformulate system \eqref{eq:NSF} in Lagrangian coordinates. This will be done in Subsection \ref{ss:lagrange}, where we essentially follow the presentation
of \cite{DFP}, adapted to the heat-conducting case.

Like in \cite{DFP}, the approach adopted here is different than the one used for \tsl{a priori} estimates. Skew-symmetry of the differential
operators at stake (in particular, of the viscous term) is crucial for the computation; this fact prompts us to perform stablity estimates in the energy space
$L^2$, instead of general $L^p$ spaces. These stability estimates are the object of Subsection \ref{ss:stab-est}. 
The core of the uniqueness argument is contained in Proposition \ref{p:G-integr}, whose proof will be the object of the final Subsection \ref{ss:proof-prop-u}.
It is precisley there where we will need to resort to the additional assumptions imposed in the last part of Theorem \ref{th:well-p}.

\subsection{Lagragian formulation} \label{ss:lagrange}

Let $\big(\rho,u,\theta\big)$ be a solution to system \eqref{eq:NSF} over $[0,T]\times\R^d$, for some time $T>0$,
satisfying the regularity conditions stated in the existence part of Theorem \ref{th:well-p}.
We immediately fix a time $T_\bull\in \,]0, T]$ such that 
\begin{equation}\label{ConditionSmallnessu}
\int_0^{T_\bull} \Vert \nabla u (t)\Vert_{L^\infty}\,\dt\, \le\, k_0\,,	
\end{equation}
for a suitably small constant $k_0>0$, whose precise value will be determined later on (at the same time, this will
determine the size of the time $T_\bull$).

We are going to reformulate system \eqref{eq:NSF} in Lagrangian coordinates in the time interval $[0,T_\bull]$ and prove uniqueness in that interval.
By a standard argument, local uniqueness implies global uniqueness, whence the conclusion of Theorem \ref{th:well-p}.

\subsubsection{Introduction to Lagrangian coordinates}
To begin with, we collect here some generalities about Lagrangian coordinates. 
We refer to \tsl{e.g.} \cite{D-2014, D-2016} for more details.

Given $u\in L^1\big([0,T_\bull];W^{1,\infty}(\R^d)\big)$, 
we define the \emph{flow} map $\psi_u$ associated to it as the solution to the integral equation  \eqref{eq:flow}.
Lagrangian coordianates are a system of coordinates based on the use of $\psi_u$. More precisely,
given any function $f$ defined over $[0,T_\bull]\times\R^d$, we can express if in terms of the Lagrangian coordinates $(t, y)$ as 
\[
\overline f(t,y)\,:=\, f\big(t, \psi_u(t,y)\big)\,. 	
\]
In particular, by writing the velocity $u$ ifself in Lagrangian coordinates, we can re-define the flux $\psi_u$ by the formula
\begin{equation}\label{DefinitionOfFlow}
\psi_u(t,y)\,=\, y\, +\, \int_0^t \overline u\big(\t,y)\,\dd\t\,. 	
\end{equation}

Let us now introduce the following convenient notation: we set
\begin{equation}\label{JuAu}
J_u\,:=\,\det\left( D\psi_u \right)\,, \qquad A_u\, :=\, \left( D \psi_u\right)^{-1}\,, \qquad \text{ so that }\quad \adj D \psi_u\, =\, J_u\, A_u\,,
\end{equation}
where $\adj D\psi_u$ denotes the adjugate matrix of $D\psi_u$.
Observe that, for any funtion $f$, one has
\[
\overline {D_xf}\,=\, D_y \oline f\cdot A_u\,.
\]
Lemma A.2 of \cite{D-2014} provides us with the following alternative differentiation rules, for passing from Eulerian to Lagrangian coordinates.
 
\begin{lemma} \label{LemmaA2Danchin}
For any $\mc C^1$ function $K$ and any $\mc C^1$ vector-field $H$ defined over $\R^d$, one has
\begin{align*}
\oline{\nabla_xK}\,&=\,J_u^{-1}\,\div_y\bigl(\adj D\psi_u\,\oline K\bigr)\,, \\
\oline{\div_xH}\,&=\,J_u^{-1}\,\div_y\bigl(\adj D\psi_u\,\oline H\bigr)\,.
\end{align*}
\end{lemma}

Let us now recall some basic estimations of the quantities above, which will be useful in what follows.
We start by recalling Lemma 4.3 of \cite{DFP}. Notice that, here, the time variable which may possibly appear can be treated as a parameter.

\begin{lemma} \label{l:equiv-L^p}
For any function $f= f(x)$, denote $f$ in lagrangian coordinates as $\overline f(t, y)\,=\, f\big(\psi_u(t, y)\big)$.
Then, for any $p \in [1, +\infty]$, one has that
		\[
		\left\|\overline f\right\|_{L^p_y}\,\leq\, \Vert J_u^{-1} \Vert_{L^\infty_y}^{1/p}\,\Vert f \Vert_{L^p_x} \qquad \mbox{ and } \qquad 
    \Vert f \Vert_{L^p_x}\le \Vert J_u\Vert_{L^\infty_y}^{1/p}\,\left\| \overline f \right\|_{L^p_y}\,. 	
		\]
	\end{lemma}

Next, we recall some estimates for the various quantities $D\psi_u$, $J_u$ and $A_u$. We refer to Appendix A of \cite{D-2014} for a proof
(with easy adaptations to the $L^\infty$ framework).

\begin{lemma}\label{LemmaJacobian0}
Let $u$ be a velocity field with $\nabla u\,\in\,L^1\bigl([0,T_\bull];L^\infty(\R^d)\bigr)$, and let $\psi_u$ be its flow, defined by \eqref{eq:flow}
or equivalently by \eqref{DefinitionOfFlow}.
Suppose that condition \eqref{ConditionSmallnessu} is fulfilled with some $k_0<1$. Define $J_u$ and $A_u$ according to \eqref{JuAu}.

Then there exists a multiplicative constant $C=C(k_0)>0$, just depending on $k_0$ and such that $C(k_0)\longrightarrow0$ for $k_0\to0$,
such that one has, for all times $t\in[0,T_\bull]$, the following estimates (where $C$ is kept implicit):
\begin{align*}
\|D\psi_u(t)\,-\,\Id\|_{L^\infty}\,&+\,\left\|\Id\,-\,A_u(t)\right\|_{L^\infty} \\  
&\quad +\,\left\|J_u^{\pm1}(t)\,-\,1\right\|_{L^\infty} 
\,+\,\left\|\Id\,-\,\adj D\psi_u(t)\right\|_{L^\infty}\,\lesssim\,\left\|Du\right\|_{L^1_{T_\bull}(L^\infty)}\,.
\end{align*}
\end{lemma}

Before going on, we point out that a direct application of the Piola identity \cite{Kup-S} yields
\begin{equation} \label{eq:Piola}
 \divergence \left(\adj D\psi_u\right)\,=\,0\,,
\end{equation}
where, according to our convention in this paper, the divergence operator of a matrix-valued function is taken column-by-column.
The above relation will be crucially exploited in the stability estimates for the temperature equation.

\subsubsection{The Navier-Stokes-Fourier system in Lagrangian coordinates} \label{sss:NSF-Lagr}
Here we recast system \eqref{eq:NSF} in Lagrangian coordinates.

First, we use Lemma \ref{LemmaA2Danchin} to get, for any function $f=f(t,x)$, the following formula:
\begin{equation}\label{ecuacionalgo3}
\overline{\partial_t f\, +\, \divergence\left( f\,u\right)} \,= \,J_u^{-1}\,\d_t\! \left( J_u\, \overline f \right)\,. 	
\end{equation}
Then, the mass balance equation in \eqref{eq:NSF} can be written as
$\partial_t \left( J_u\, \overline \rho \right)\,=\,0$, 	
which immediately implies the equality
\begin{equation}\label{DensityInLC}
J_u\, \overline \rho\,=\, \rho_0\,. 	
\end{equation}

Next, we write the operator $\mathcal L$ in Lagrangian coordinates as 
\begin{equation} \label{def:L-Lagr}
 \wtilde{\mc L}\,f\,:=\,-\,\mu\,\div_y\Big(\adj(D\psi_u)\;^t\!A_u\,\nabla_y f\Big)\,-\,
(\lambda+\mu)\,\div_y\Big(\adj(D\psi_u)\;(\,^t\!A_u:\nabla_y f)\Big)\,.
\end{equation}
Hence, taking into account \eqref{ecuacionalgo3} and \eqref{DensityInLC}, the equation for the velocity field $u$ can be expressed in Lagrangian coordinates as 
\begin{equation}\label{EqLagrangeU}
\rho_0\, \partial_t \overline u\, +\, \wtilde{\mathcal L}\,\overline u \,=\, -\, \divergence_y\left( \adj D\psi_u\, P\left( J_u^{-1}\,\rho_0, \overline \theta \right)\right)\,.
\end{equation}

Finally, for the equation for the temperature, using again \eqref{ecuacionalgo3} and \eqref{DensityInLC}, we have
\begin{align*}
&\rho_0\, \partial_t \overline \theta\, +\, \overline \theta \, \pi_1\left( J_u^{-1}\,\rho_0\right)\,\divergence_y\left( \adj D\psi_u\; \overline u \right)
\\
&\qquad\qquad\qquad -\,\kappa\, \divergence_y \left( \adj D\psi_u\; ^t\!A_u\, \nabla_y \overline \theta \right)
\,= \,J_u\,\overline {\lambda \vert \divergence u \vert^2 + \mu \vert \nabla u\vert^2 + \mu \nabla u : \nabla u}. 
\end{align*}
To deal with the right-hand side terms, we perform the following computations.
\begin{itemize}
\item We know that $\overline{\left( \divergence u \right)^2}\,=\, \left( \overline {\divergence u}\right)^2$; thus, from Lemma \ref{LemmaA2Danchin}
we infer
\[
\overline{\left( \divergence u \right)^2}\, =\, J_{u}^{-2}\,\big(\divergence_y \left( \adj D \psi_u \, \overline u \right)\big)^2\,.
\]
\item Since $\overline {D_x u}\, =\, D_y \overline u\, A_u$, we immediately obtain 
\[
\oline{\vert \nabla u\vert^2}\, =\,  \vert D_y \overline u\, A_u \vert^2\,.
\]
\item Finally, by definition we have that $\nabla u:\nabla u\,=\,\tr\left((\nabla u)^2\right)$. So, by use of the formula
$\overline {D_x u}\, =\, D_y\overline u \,A_u$ again, we find
\[
\overline {\nabla u : \nabla u }\,=\,\tr\left( \left(D_y \overline u\, A_u\right)^2 \right)\,.
\]
\end{itemize}
We conclude that the equation for the temperature can be writen as 
\begin{align}
\label{EqLagrTemp}
&\rho_0\, \partial_t \overline \theta\, +\, \overline \theta \, \pi_1\left( J_u^{-1}\,\rho_0\right)\,\divergence_y\left( \adj D\psi_u\; \overline u \right)\, -\,
\kappa\, \divergence_y \left( \adj D\psi_u\; ^t\!A_u\, \nabla_y \overline \theta \right) \\
\nonumber
& \qquad\qquad =\, \lambda\, J_{u}^{-1}\,\big(\divergence_y \left( \adj D \psi_u \, \overline u \right)\big)^2\, +\,  \mu\, J_u\,\vert D_y \overline u\, A_u \vert^2\,+\,
\mu\, J_u\,\tr\left( \left(D_y \overline u\, A_u\right)^2 \right)\,.
\end{align}

\subsection{Difference of solutions in Lagrangian coordinates} \label{ss:stab-est}
This part is the core of the uniqueness proof. 
We assume to have two sets of initial data
\[
\Big( \rho_{0,j}\,,\, u_{0,j}\,,\, \theta_{0,j}\Big)\,, \qquad \mbox{ where }\ j\,=\,1\,,\, 2\,.
\]
Assume that they satisfy the assumptions of Theorem \ref{th:well-p}. For simplicity,
let us assume that the initial densities are equal, that is
\[
\rho_{0,1}\,=\,\rho_{0,2}\,,
\]
whereas the initial velocity and temperature fields may differ.
We then obtain two sets of solutions for system \eqref{eq:NSF}, namely $\big( \rho_1, u_1, \theta_1\big)$ and $\big( \rho_2, u_2, \theta_2\big)$,
defined on a same time interval $[0,T_\bull]$. Recall that those solutions possess the regularity properties stated in the existence part of
Theorem \ref{th:well-p} on $[0,T_\bull]\times\R^d$.
Recall that $T_\bull$ has been chosen so small that condition \eqref{ConditionSmallnessu} holds true, for a suitable constant $k_0<1$, to be fixed later on.

For convenience, from now on we adopt the following notation: for $j=1,2$, we set
\begin{equation}\label{Notation1}
\psi_j(t,y)\,=\,\psi_{u_j}\,, \qquad J_j\,=\,J_{u_j} 
\qquad \text{ and } \qquad 
A_j\,=\,A_{u_j}\,. 
\end{equation}
In addition, let us define 
\begin{equation}\label{Notation2}
\delta {\overline u}\,:=\, \overline u_1 - \overline u_2 \qquad \mbox{ and }\qquad
\delta {\overline \theta}\,:=\, \overline \theta_1 - \overline \theta_2\,. 
\end{equation}
Similarly, for any quantity $f_j$ related to the solutions $\big(\oline u_j,\oline\theta_j\big)$, with $j=1,2$, we set
$\delta f\,:=\,f_1-f_2$. So, for instance, we have $\delta {\adj}\,:=\, \adj D\psi_1-\adj D\psi_2$ and,
after setting $\wtilde{\mc L}_j$ to be the operator $\widetilde{\mathcal L}$ associated with
the velocity field $u_j$ (recall \eqref{def:L-Lagr} above),
$\delta {\widetilde {\mathcal L}}\,:=\,\widetilde{ \mathcal L}_1 -\widetilde{ \mathcal L}_2$.

Then, using \eqref{EqLagrangeU} and \eqref{EqLagrTemp}, we can derive an equation for $\delta{\overline u}$ and $\delta{\overline \theta}$.
First of all, we notice that the equation for $\de\oline u$ can be written as
\begin{align}\label{EqDeltau}
\rho_0\, \partial_t \delta{\overline u}\, +\, \mathcal L \,\delta{\overline u}\,&=\, 
\,\left(\mathcal L\,-\,\widetilde {\mathcal L}_1\right)\,\delta{ \overline u} \,-\, \delta{\widetilde{\mathcal L}}\,\overline u_2\,-\,
\divergence \Big( \delta{\adj}\, P\left( J_{1}^{-1}\rho_0, \overline \theta_1\right)\Big) \\
\nonumber
&\qquad\qquad -\,\divergence \bigg( \adj D\psi_2\,\Big( P\left( J_{1}^{-1}\rho_0, \overline \theta_1\right)\,-\,
P\left( J_{2}^{-1}\rho_0, \overline \theta_2\right)\Big)\bigg)\,. 
\end{align}
Here, the operator $\mc L$ is the ``true'' Lam\'e operator defined in \eqref{DefinitionL}. Indeed, Lemma \ref{LemmaJacobian0} guarantees that,
by taking $k_0$ small enough in \eqref{ConditionSmallnessu}, both $\wtilde{\mc L}_1$ and  $\wtilde{\mc L}_2$ remain close to it.

Let us now turn our attention to the equation for $\de\oline\theta$. By use of \eqref{EqLagrTemp}, straightforward computations
give us the following relation:
\begin{align}
\label{EqDeltatheta}
&\hspace{-0.2cm}
\rho_0\, \partial_t\de \overline \theta\,-\,\k\,\Delta\de\oline\theta \\
\nonumber
&=\,\k\,\div\Big( \left(\adj D\psi_1\; ^t\!A_1\,-\,\Id\right)\, \nabla \de\overline \theta \Big)\,-\,
\k\,\div\Big( \left(\de\adj\;^t\!A_1\,+\,\adj D\psi_2\; ^t\!\de A\right)\, \nabla \overline \theta_2 \Big) \\
\nonumber
&\quad -\,\de\oline\theta\,\pi_1\left( J_1^{-1}\,\rho_0\right)\,\divergence\left( \adj D\psi_1\; \overline u_1 \right)\,-\,
\oline\theta_2\,\Big(\pi_1\left( J_1^{-1}\,\rho_0\right)\,-\,\pi_1\left( J_2^{-1}\,\rho_0\right)\Big)\,\divergence\left( \adj D\psi_1\; \overline u_1 \right) \\
\nonumber
&\quad -\,\oline\theta_2\,\pi_1\left( J_2^{-1}\,\rho_0\right)\,\divergence\Big( \de\adj\; \overline u_1 \Big)
-\,\oline\theta_2\,\pi_1\left( J_2^{-1}\,\rho_0\right)\,\divergence\Big( \adj D\psi_2\; \de\overline u \Big) \\
\nonumber
&\quad
+\,\lambda\, \de J^{-1}\,\big(\divergence \left( \adj D \psi_1 \, \overline u_1 \right)\big)^2 \\
\nonumber
&\quad
+\,\lambda\, J_{2}^{-1}\,\divergence\Big(\de\adj \, \overline u_1\,+\,\adj D \psi_2 \,\de\overline u \Big)\;
\div\Big(\adj D \psi_1 \, \overline u_1\,+\,\adj D \psi_2 \, \overline u_2\Big) \\
\nonumber
&\quad
+\,\mu\,\de J\,\vert D \overline u_1\, A_1 \vert^2\,+\,
\mu\,J_2\,\sum_{j,k=1}^d\left(D\de\oline u\,A_1\,+\,D\oline u_2\,\de A\right)_{j,k}\;\left(D \overline u_1\, A_1\,+\,D \overline u_2\, A_2\right)_{j,k} \\
\nonumber
&\quad
+\,\mu\,\de J\,\tr\left( \left(D \overline u_1\, A_1\right)^2 \right)\,+\,
\mu\,J_2\,\tr\left( \left(D\de \overline u\, A_1\,+\,D\oline u_2\,\de A\right)\cdot\Big(D\oline u_1\,A_1\,+\,D\oline u_2\,A_2\Big) \right)\,.
\end{align}

The strategy now consists in performing an energy estimate on equations \eqref{EqDeltau} and \eqref{EqDeltatheta}.

\subsubsection{Energy estimate for equation \eqref{EqDeltau}}

First, we do an energy estimates on $\delta {\overline u}$. For this, we multiply equation \eqref{EqDeltau} by $\delta{\overline u}$ and integrate over $\R^d$,
and get
\begin{align}
\label{est:du_prelim}
&\dfrac12  \, \frac{\dd}{\dt}\int \rho_0 \vert \delta{ \overline u} \vert^2\,\dd y\, +\,
\int \mathcal L \delta{ \overline u }\cdot \delta {\overline u }\,\dd y \\
\nonumber
& \qquad =\,\int \left(\mathcal L\,-\,\widetilde {\mathcal L}_1\right)\,\delta{ \overline u}\cdot \delta{ \overline u}\, \dd y
\,+\,\int\delta{\widetilde{\mathcal L}}\,\overline u_2\cdot\de\oline u\,\dd y\,-\,
\int \divergence \Big( \delta{\adj}\, P\left( J_{1}^{-1}\rho_0, \overline \theta_1\right)\Big)\cdot\de\oline u\,\dd y \\
\nonumber
&\qquad\qquad\qquad\qquad\qquad\qquad -\,\int
\divergence \bigg( \adj D\psi_2\,\Big( P\left( J_{1}^{-1}\rho_0, \overline \theta_1\right)\,-\,P\left( J_{2}^{-1}\rho_0, \overline \theta_2\right)\Big)\bigg)
\cdot \de\oline u\, \dd y\,. 
\end{align}

Keeping definition \eqref{DefinitionL} in mind, after an integration by parts we easily see that
\begin{equation} \label{eq:Lu-lagrange}
 \int \mathcal L \delta{ \overline u }\cdot \delta {\overline u }\,\dd y\,=\,
\mu\,\int\left|\nabla \de\oline u\right|^2\,\dd y\,+\,(\mu+\lambda)\,\int\left|\div \de\oline u\right|^2\,\dd y\,.
\end{equation}
In order to control the right-hand side of \eqref{est:du_prelim}, instead, we resort to Lemma 4.4 of \cite{DFP},
which we recall here below for the reader's convenience (see also Lemma A.4 of \cite{D-2014} for details).

\begin{lemma} \label{l:est_delta_A-J}
Assume that condition \eqref{ConditionSmallnessu} is fulfilled by $u_1$ and $u_2$, for some $k_0\in\,]0,1]$.

Then, there exists a constant $C>0$, proportional to $k_0$ or positive powers of it, such that the following estimates hold true,
for all $t\in[0,T_\bull]$ and all $p\in[1,+\infty]$:
\begin{align*}
\left\|\de\adj(t)\right\|_{L^p}\,+\,\left\|\de A\right\|_{L^p}\,+\,
\left\|\de J^{\pm1}(t)\right\|_{L^p}\,&\leq\,C\,\int_0^t\left\|\nabla\de u(\tau)\right\|_{L^p}\,\dd\tau\,.
\end{align*}
\end{lemma}

Now, let us come back to the bounds for the right-hand side of \eqref{est:du_prelim}.
Firstly, employing Lemmas \ref{LemmaJacobian0} and \ref{l:est_delta_A-J} and making use of \eqref{ConditionSmallnessu},
we deduce the following bounds, for implicit multiplicative constants which may depend on the time $T_\bull$ and on $k_0$,
but not degenerating for $k_0\to0$:
\begin{align}
\label{est:du_1}
\left|\int \left(\mathcal L\,-\,\widetilde {\mathcal L}_1\right)\,\delta{ \overline u}\cdot \delta{ \overline u}\, \dd y\right|\,&\lesssim\,
k_0\,\left\|\nabla\de\oline{u}\right\|^2_{L^2}\,, \\
\label{est:du_2}
\left|\int\delta{\widetilde{\mathcal L}}\,\overline u_2\cdot\de\oline u\,\dd y\right|\, &\lesssim
\left(\int^t_0\|\nabla\de\oline u\|_{L^2}\,\dd\tau\right)\,\left\|\nabla\oline{u}_2\right\|_{L^\infty}\,\left\|\nabla\de\oline{u}\right\|_{L^2}\,.
\end{align}

Next, using Lemmas \ref{LemmaJacobian0} and \ref{l:est_delta_A-J} and condition \eqref{ConditionSmallnessu} again, we find
\begin{align*}
\left|\int \divergence \Big( \delta{\adj}\, P\left( J_{1}^{-1}\rho_0, \overline \theta_1\right)\Big)\cdot\de\oline u\,\dd y\right|\, &\lesssim
\left(\int^t_0\|\nabla\de\oline u\|_{L^2}\,\dd\tau\right)\,\left\|P\left( J_{1}^{-1}\rho_0, \overline \theta_1\right)\right\|_{L^\infty}\,
\left\|\nabla\de\oline{u}\right\|_{L^2}\,.
\end{align*}
Observe that $P\left( J_{1}^{-1}\rho_0, \overline \theta_1\right)\,=\,P(\oline\rho_1,\oline\theta_1)$. Thus, by definition \eqref{DefinitionP} of the pressure
function, making use of Lemmas \ref{l:equiv-L^p} and \ref{LemmaJacobian0} and of bound \eqref{ConditionSmallnessu}, we find
\[
 \left\|P\left( J_{1}^{-1}\rho_0, \overline \theta_1\right)\right\|_{L^\infty}\,\lesssim\,1\,+\,\left\|\oline\theta_1\right\|_{L^\infty}\,,
\]
for an implicit multiplicative constant with may depend on $k_0$ and on the $L^\infty$ norm of $\oline\rho_1$ over $[0,T_\bull]\times\R^d$.
Thus, in turn we get
\begin{align}
 \label{est:du_3}
 \left|\int \divergence \Big( \delta{\adj}\, P\left( J_{1}^{-1}\rho_0, \overline \theta_1\right)\Big)\cdot\de\oline u\,\dd y\right|\, &\lesssim
\left(\int^t_0\|\nabla\de\oline u\|_{L^2}\,\dd\tau\right)\,\Big(1+\left\|\oline\theta_1\right\|_{L^\infty}\Big)\,
\left\|\nabla\de\oline{u}\right\|_{L^2}\,.
\end{align}

Concerning the last term appearing in \eqref{est:du_prelim}, we can proceed similarly, but with some more care: to begin with,
we can bound it as
\begin{align*}
&\left|\int\div \bigg( \adj D\psi_2\,\Big( P\left( J_{1}^{-1}\rho_0, \overline \theta_1\right)\,-\,P\left( J_{2}^{-1}\rho_0, \overline \theta_2\right)\Big)\bigg)
\cdot \de\oline u\, \dd y\right| \\
&\qquad\qquad\qquad
\lesssim\,\left\|\nabla\de\oline{u}\right\|_{L^2}\,\Big(\left\|\d_\rho P\right\|_{L^\infty}\,\left\|\rho_0\right\|_{L^\infty}\,\left\|\de J^{-1}\right\|_{L^2}\,+\,
\left\|\d_\theta P\right\|_{L^\infty}\,\left\|\de\oline\theta\right\|_{L^2}\Big) \\
&\qquad\qquad\qquad 
\lesssim\,\left(\int^t_0\|\nabla\de\oline u\|_{L^2}\,\dd\tau\right)\,\left\|\d_\rho P\right\|_{L^\infty}\,\left\|\nabla\de\oline{u}\right\|_{L^2}\,+\,
\left\|\d_\theta P\right\|_{L^\infty}\,\left\|\de\oline\theta\right\|_{L^2}\,\left\|\nabla\de\oline{u}\right\|_{L^2}\,,
\end{align*}
where, of course, we agree that the $L^\infty$ norm of $\d_\rho P$ and $\d_\theta P$ must be taken over the support of the two sets of solutions
$\big(\oline\rho_1,\oline\theta_1\big)$ and $\big(\oline\rho_2,\oline\theta_2\big)$. 
We observe that, owing to \eqref{DefinitionP}, we have
\[
 \d_\rho P(a,b)\,=\,\pi_0'(a)\,+\,b\,\pi_1'(a)\qquad \mbox{ and }\qquad
 \d_\theta P(a,b)\,=\,\pi_1(a)\,.
\]
Keep in mind that $\oline\rho_1$ and $\oline\rho_2$ are uniformly bounded on $[0,T_\bull]\times \R^d$, owing to the regularity properties
stated in Theorem \ref{th:well-p} for their Eulerian counterparts, combined with 
Lemmas \ref{l:equiv-L^p} and \ref{LemmaJacobian0}, and also with \eqref{ConditionSmallnessu}. It follows that
\[
 \left\|\d_\rho P\right\|_{L^\infty}\,\lesssim\,1\,+\,\left\|\oline\theta_1\right\|_{L^\infty}\,+\,\left\|\oline\theta_2\right\|_{L^\infty}
\qquad \mbox{ and }\qquad \left\|\d_\theta P\right\|_{L^\infty}\,\lesssim\,1\,,
\]
whence in turn we find
\begin{align}
\label{est:du_4}
&\left|\int\div \bigg( \adj D\psi_2\,\Big( P\left( J_{1}^{-1}\rho_0, \overline \theta_1\right)\,-\,P\left( J_{2}^{-1}\rho_0, \overline \theta_2\right)\Big)\bigg)
\cdot \de\oline u\, \dd y\right| \\
\nonumber
&\qquad\qquad
\lesssim\,\left(\int^t_0\|\nabla\de\oline u\|_{L^2}\,\dd\tau\right)\,\Big(1\,+\,\left\|\oline\theta_1\right\|_{L^\infty}\,+\,\left\|\oline\theta_2\right\|_{L^\infty}\Big)\,
\left\|\nabla\de\oline{u}\right\|_{L^2}\,+\,\left\|\de\oline\theta\right\|_{L^2}\,\left\|\nabla\de\oline{u}\right\|_{L^2}\,.
\end{align}

Collecting relation \eqref{eq:Lu-lagrange} and estimates \eqref{est:du_1}, \eqref{est:du_2}, \eqref{est:du_3} and \eqref{est:du_4},
and plugging them into \eqref{est:du_prelim}, in the end we obtain
\begin{align*}
&\dfrac12\, \frac{\dd}{\dt}\int \rho_0 \vert \delta{ \overline u} \vert^2\,\dd y\, +\,
\mu\,\int\left|\nabla \de\oline u\right|^2\,\dd y\,+\,(\mu+\lambda)\,\int\left|\div \de\oline u\right|^2\,\dd y \\
& \qquad \lesssim\,k_0\,\left\|\nabla\de\oline{u}\right\|^2_{L^2}\,+\,
\left\|\de\oline\theta\right\|_{L^2}\,\left\|\nabla\de\oline{u}\right\|_{L^2} \\
&\qquad\qquad\qquad\qquad\qquad +\,\left(\int^t_0\|\nabla\de\oline u\|_{L^2}\,\dd\tau\right)\,
\Big(1+\left\|\big(\oline\theta_1,\oline\theta_2,\nabla\oline{u}_2\big)\right\|_{L^\infty}\Big)\,\left\|\nabla\de\oline{u}\right\|_{L^2}\,,
\end{align*}
where, from now on, we adopt the convenient notation $\|(a,b)\|_{\mc B}\,=\,\|a\|_{\mc B}+\|b\|_{\mc B}$.
By a direct application of the Young inequality, 
and taking $k_0$ in \eqref{ConditionSmallnessu} small enough compared to $\mu>0$, we thus get
\begin{align*}
&\frac{\dd}{\dt}\left\|\sqrt{\rho_0}\, \delta{ \overline u}\right\|^2_{L^2}\, +\,
\mu\,\left\|\nabla \de\oline u\right\|^2_{L^2}\, \lesssim\,
\left\|\de\oline\theta\right\|^2_{L^2}\,+\,
\left(\int^t_0\|\nabla\de\oline u\|_{L^2}\,\dd\tau\right)^2\,
\Big(1+\left\|\big(\oline\theta_1,\oline\theta_2,\nabla\oline{u}_2\big)\right\|_{L^\infty}\Big)^2\,.
\end{align*}
for a new implicit multiplicative constant, depending on $k_0$, $\mu$ and $T_\bull$.
By a simple application of the Jensen inequality, we finally conclude that
\begin{align}
\label{est:du}
&\frac{\dd}{\dt}\left\|\sqrt{\rho_0}\, \delta{ \overline u}\right\|^2_{L^2}\, +\,
\mu\,\left\|\nabla \de\oline u\right\|^2_{L^2} \\
\nonumber
& \qquad\qquad\qquad\qquad \lesssim\,
\left\|\de\oline\theta\right\|^2_{L^2}\,+\,
t\,\Big(1+\left\|\big(\oline\theta_1,\oline\theta_2,\nabla\oline{u}_2\big)\right\|^2_{L^\infty}\Big)\,\int^t_0\|\nabla\de\oline u\|^2_{L^2}\,\dd\tau\,\,.
\end{align}

\subsubsection{Energy estimate for equation \eqref{EqDeltatheta}}

Now, we estimate $\delta \overline \theta$ by taking the $L^2$ scalar product of equation \eqref{EqDeltatheta} with $\de\oline\theta$ itself.
It goes without saying that
\begin{align}
\label{est:dT_0}
\int\Big(\rho_0\, \partial_t\de \overline \theta\,-\,\k\,\Delta\de\oline\theta\Big)\,\de\oline\theta\,\dd y\,=\,
\frac{1}{2}\,\frac{\dd}{\dt}\int\rho_0\,\left|\de\oline\theta\right|^2\,\dd y\,+\,\k\,\int\left|\nabla\de\oline\theta\right|^2\,\dd y\,.
\end{align}

Next, we are going to bound the integrals related to the right-hand side terms of \eqref{EqDeltatheta} one by one. To begin with,
we see that, repeating \tsl{mutatis mutandis} the computations leading to \eqref{est:du_1} and \eqref{est:du_2}, we find
\begin{align}
\label{est:dT_1}
&\left|\k\int\div\Big( \left(\adj D\psi_1\; ^t\!A_1\,-\,\Id\right)\, \nabla \de\overline \theta \Big)\,\de\oline\theta \,\dd y\right|\,\lesssim\,
\k\,k_0\,\left\|\nabla\de\oline\theta\right\|_{L^2}^2\,, \\
\label{est:dT_2}
&\hspace{-0.9cm}
\left|\k\int\,\div\Big( \left(\de\adj\;^t\!A_1\,+\,\adj D\psi_2\; ^t\!\de A\right)\, \nabla \overline \theta_2 \Big)\de\oline\theta \,\dd y\right| \\
\nonumber
&\qquad\qquad\qquad\qquad\qquad \lesssim\,
\k\,\left(\int^t_0\|\nabla\de\oline u\|_{L^2}\,\dd\tau\right)\,\left\|\nabla\oline{\theta}_2\right\|_{L^\infty}\,\left\|\nabla\de\oline{\theta}\right\|_{L^2} \\
\nonumber
&\qquad\qquad\qquad\qquad\qquad \lesssim\,\k\,k_0\,\left\|\nabla\de\oline\theta\right\|_{L^2}^2\,+\,
\k\,t\,\left\|\nabla\oline{\theta}_2\right\|^2_{L^\infty}\,\int^t_0\|\nabla\de\oline u\|^2_{L^2}\,\dd\tau\,,
\end{align}
where, in the last step, we have argued precisely as done above to obtain inequality \eqref{est:du}.

Also the estimate of the third term on the right-hand side of \eqref{EqDeltatheta} is easy:
making use of the uniform boundedness of the densities $\oline\rho_{1,2}$ over $[0,T_\bull]\times\R^d$, together with Lemma \ref{LemmaJacobian0}, condition \eqref{ConditionSmallnessu} and the Piola identity \eqref{eq:Piola}, we get
\begin{align}
\label{est:dT_3}
 \left|\int \de\oline\theta\,\pi_1\left( J_1^{-1}\,\rho_0\right)\,\divergence\left( \adj D\psi_1\; \overline u_1 \right)\,\de\oline\theta\,\dd y\right|\,&\lesssim\,
 \left\|\nabla\oline u_1\right\|_{L^\infty}\,\left\|\de\oline\theta\right\|_{L^2}^2\,.
\end{align}

For the other terms appearing on the right-hand side of \eqref{EqDeltatheta}, instead, we need to proceed carefully.
As a matter of fact, those terms are more delicate to treat, essentially because of
the scarse time integrability of certain norms of $\nabla\oline u_j$ and $\oline\theta_j$, for $j=1,2$.
It is precisely here that the additional assumptions imposed in the uniqueness part of Theorem \ref{th:well-p} come into play.
However, in order isolate the difficulty and show the necessity of those assumptions, let us treat the above mentioned terms by means of a precise interpolation argument
and push the existence hypotheses up to their boundaries.
The strategy is the following. All those terms depend (either directly or indirectly, see Lemma \ref{l:est_delta_A-J})
on the quantity $\nabla\de\oline u$, which we are forced to control in $L^2$.
At the same time, the term $\de\oline\theta$ coming from the $L^2$ scalar product is
multiplied (roughly speaking) by another factor, depending on either $\nabla\oline u_j$ or $\oline\theta_j$.
In order to gain time integrability, the idea is to use H\"older inequality and control the $L^2$ norm of the product not in $L^2\times L^\infty$,
but rather in $L^{q_1}\times L^{q_2}$, for $q_{1,2}\in\,]2,+\infty[\,$ such that $1/q_1\,+\,1/q_2\,=\,1/2$. 
By writing $\nabla\oline u_j\,=\,\nabla \oline w_j\,+\,\nabla\oline v_j$ and
taking into account the properties of the solutions (and also Lemma \ref{l:equiv-L^p}), it is convenient to bound
$\nabla\oline v_j$ and $\oline\theta_j$ in $L^p$. This consideration prompts us to fix $q\in\,]2,+\infty[\,$ as
\begin{equation} \label{def:q}
 \frac{1}{q}\,:=\,\frac{1}{2}\,-\,\frac{1}{p}\,.
\end{equation}
The same argument applies also for the term $\nabla\oline w_j$. However,
since we only have $\nabla^2w\in L^r_{T_\bull}(L^p)$ with $r\in\,]1,2[\,$ small, direct computations show that the related term in the final Gr\"onwall
argument would not belong to $L^1_{T_\bull}$. How to overcome this issue will be explained in detail in Subsection \ref{ss:proof-prop-u}.

With the above considerations in mind, we can now estimate the remaining terms arising from the right-hand side of equation \eqref{EqDeltatheta}.
Here below, we will make a systematic use of the uniform boundedness of the densities $\oline\rho_1$ and $\oline\rho_2$ over $[0,T_\bull]\times\R^d$,
as well as of Lemma \ref{LemmaJacobian0}, condition \eqref{ConditionSmallnessu} and the Piola identity \eqref{eq:Piola}.
Then, we proceed to bound the pressure-dependent terms in the following way:
\begin{align*}
&\left|\int\oline\theta_2\,\Big(\pi_1\left( J_1^{-1}\,\rho_0\right)\,-\,\pi_1\left( J_2^{-1}\,\rho_0\right)\Big)\,
\divergence\left( \adj D\psi_1\; \overline u_1 \right)\,\de\oline\theta\,\dd y\right| \\
&\qquad\qquad\qquad\qquad\qquad\qquad\qquad\qquad\qquad
\lesssim\,\left\|\oline\theta_2\right\|_{L^p}\,\left\|\nabla\oline u_1\right\|_{L^\infty}\,\left\|\de\oline\theta\right\|_{L^q}\,\left\|\de J^{-1}\right\|_{L^2}\,, \\
&\quad
\left|\int\oline\theta_2\,\pi_1\left( J_2^{-1}\,\rho_0\right)\,\divergence\Big( \de\adj\; \overline u_1 \Big)\,\de\oline\theta\,\dd y\right|\,
\lesssim\,\left\|\oline\theta_2\right\|_{L^p}\,\left\|\nabla\oline u_1\right\|_{L^\infty}\,\left\|\de\oline\theta\right\|_{L^q}\,\left\|\de\adj\right\|_{L^2}\, \\
&\hspace{-0.2cm}\left|\int\oline\theta_2\,\pi_1\left( J_2^{-1}\,\rho_0\right)\,\divergence\Big( \adj D\psi_2\; \de\overline u \Big)\,\de\oline\theta\,\dd y\right|\,
\lesssim\,\left\|\oline\theta_2\right\|_{L^p}\,\left\|\de\oline\theta\right\|_{L^q}\,\left\|\nabla\de\oline u\right\|_{L^2}\,.
\end{align*}
Here, $q\in\,]2,+\infty[\,$ is defined as in \eqref{def:q}. 
By the same token, we can derive the following estimates, where, for the sake of better clarity, we split the terms $\oline v_j$ and $\oline w_j$:
\begin{align*}
&\left|\int\de J^{-1}\,\big(\divergence \left( \adj D \psi_1 \, \overline u_1 \right)\big)^2\,\,\de\oline\theta\,\dd y\right| \\
&\qquad \lesssim\,
\left\|\nabla\oline v_1\right\|_{L^p}\,\left\|\nabla\oline u_1\right\|_{L^\infty}\,\left\|\de\oline\theta\right\|_{L^q}\,\left\|\de J^{-1}\right\|_{L^2}\,+\,
\left\|\nabla\oline w_1\right\|_{L^{p}}\,\left\|\nabla\oline u_1\right\|_{L^\infty}\,\left\|\de\oline\theta\right\|_{L^{q}}\,\left\|\de J^{-1}\right\|_{L^2}
\,, \\
&\left|\int J_{2}^{-1}\,\divergence\Big(\de\adj \, \overline u_1\Big)\,\sum_{j=1,2}\div\Big(\adj D \psi_j \, \overline u_j\Big)\,\,\de\oline\theta\,\dd y\right| \\
&\qquad
\lesssim\,\left\|\nabla\oline v_1\right\|_{L^p}\,\left\|\big(\nabla\oline u_1,\nabla\oline u_2\big)\right\|_{L^\infty}\,\left\|\de\oline\theta\right\|_{L^q}\,
\left\|\de\adj\right\|_{L^2} \\
&\qquad\qquad\qquad\qquad\qquad\qquad\qquad\qquad+\,
\left\|\nabla\oline w_1\right\|_{L^{p}}\,\left\|\big(\nabla\oline u_1,\nabla\oline u_2\big)\right\|_{L^\infty}\,
\left\|\de\oline\theta\right\|_{L^{q}}\,\left\|\de J^{-1}\right\|_{L^2}\,, \\
&\left|\int J_{2}^{-1}\,\divergence\Big(\adj D \psi_2 \,\de\overline u \Big)\,\sum_{j=1,2}\div\Big(\adj D \psi_j \, \overline u_j\Big)\,\de\oline\theta\,\dd y
\right| \\
&\qquad
\lesssim\,
\left\|\big(\nabla\oline v_1,\nabla\oline v_2\big)\right\|_{L^p}\,\left\|\de\oline\theta\right\|_{L^q}\,\left\|\nabla\de\oline u\right\|_{L^2}
\,+\,
\left\|\big(\nabla\oline w_1,\nabla\oline w_2\big)\right\|_{L^{p}}\,\left\|\de\oline\theta\right\|_{L^{q}}\,\left\|\nabla\de\oline u\right\|_{L^2}
\,, \\
&
\left|\int\de J\,\vert D \overline u_1\, A_1 \vert^2\,\de\oline\theta\,\dd y\right| \\
&\qquad \lesssim\,
\left\|\nabla\oline v_1\right\|_{L^p}\,\left\|\nabla\oline u_1\right\|_{L^\infty}\,\left\|\de\oline\theta\right\|_{L^q}\,\left\|\de J^{-1}\right\|_{L^2}\,+\,
\left\|\nabla\oline w_1\right\|_{L^{p}}\,\left\|\nabla\oline u_1\right\|_{L^\infty}\,\left\|\de\oline\theta\right\|_{L^{q}}\,\left\|\de J^{-1}\right\|_{L^2}\,, \\
&\left|\int J_2\,\sum_{j,k=1}^d\left(D\de\oline u\,A_1\right)_{j,k}\;\left(D \overline u_1\, A_1\,+\,D \overline u_2\, A_2\right)_{j,k}\,\de\oline\theta\,\dd y\right| \\
&\qquad \lesssim\,
\left\|\big(\nabla\oline v_1,\nabla\oline v_2\big)\right\|_{L^p}\,\left\|\de\oline\theta\right\|_{L^q}\,\left\|\nabla\de\oline u\right\|_{L^2}
\,+\,
\left\|\big(\nabla\oline w_1,\nabla\oline w_2\big)\right\|_{L^{p}}\,\left\|\de\oline\theta\right\|_{L^{q}}\,\left\|\nabla\de\oline u\right\|_{L^2}
\,,
\\
&\left|\int J_2\,\sum_{j,k=1}^d\left(D\oline u_2\,\de A\right)_{j,k}\;\left(D \overline u_1\, A_1\,+\,D \overline u_2\, A_2\right)_{j,k}\,\de\oline\theta\,\dd y\right| \\
&\qquad\lesssim\,
\left\|\nabla\oline v_2\right\|_{L^p}\,\left\|\big(\nabla\oline u_1,\nabla\oline u_2\big)\right\|_{L^\infty}\,\left\|\de\oline\theta\right\|_{L^q}\,
\left\|\de A\right\|_{L^2}  \\
&\qquad\qquad\qquad\qquad\qquad\qquad\qquad\qquad+\,
\left\|\nabla\oline w_2\right\|_{L^{p}}\,\left\|\big(\nabla\oline u_1,\nabla\oline u_2\big)\right\|_{L^\infty}\,\left\|\de\oline\theta\right\|_{L^{q}}\,
\left\|\de A\right\|_{L^2}\,.
\end{align*}
Similarly, we can bound the last terms, where the trace operator appears:
\begin{align*}
\left|\int \de J\,\tr\left( \left(D \overline u_1\, A_1\right)^2 \right)\,\de\oline\theta\,\dd y\right|\,&\lesssim\,
\left\|\nabla\oline v_1\right\|_{L^p}\,\left\|\nabla\oline u_1\right\|_{L^\infty}\,\left\|\de\oline\theta\right\|_{L^q}\,\left\|\de J^{-1}\right\|_{L^2} \\
&\qquad +\,
\left\|\nabla\oline w_1\right\|_{L^{p}}\,\left\|\nabla\oline u_1\right\|_{L^\infty}\,\left\|\de\oline\theta\right\|_{L^{q}}\,\left\|\de J^{-1}\right\|_{L^2}
\,, \\
\left|\int J_2\,\tr\left( D\de \overline u\, A_1\cdot\sum_{j=1,2}D\oline u_j\,A_j \right)\,\de\oline\theta\,\dd y\right|\,&\lesssim\,
\left\|\big(\nabla\oline v_1,\nabla\oline v_2\big)\right\|_{L^p}\,\left\|\de\oline\theta\right\|_{L^q}\,\left\|\nabla\de\oline u\right\|_{L^2} \\
&\qquad +\,
\left\|\big(\nabla\oline w_1,\nabla\oline w_2\big)\right\|_{L^{p}}\,\left\|\de\oline\theta\right\|_{L^{q}}\,\left\|\nabla\de\oline u\right\|_{L^2}\,, \\
%
\left|\int J_2\tr\left( D\oline u_2\de A\cdot\sum_{j=1,2}D\oline u_jA_j \right)\de\oline\theta\dd y\right|\,&\lesssim\,
\left\|\nabla\oline v_2\right\|_{L^p}\,\left\|\big(\nabla\oline u_1,\nabla\oline u_2\big)\right\|_{L^\infty}\,\left\|\de\oline\theta\right\|_{L^q}\,
\left\|\de A\right\|_{L^2} \\
&\qquad +
\left\|\nabla\oline w_2\right\|_{L^{p}}\left\|\big(\nabla\oline u_1,\nabla\oline u_2\big)\right\|_{L^\infty}
\left\|\de\oline\theta\right\|_{L^{q}}\left\|\de A\right\|_{L^2}.
\end{align*}

Gathering all these inequalities together,
we find that their sum (which we call $\mf S$ for convenience of notation) may be controlled by the following quantity:
\begin{align*}
\mf S\,&\lesssim\,
\left\|\big(\oline\theta_2,\nabla\oline v_1,\nabla\oline v_2\big)\right\|_{L^p}\,\left\|\de\oline\theta\right\|_{L^q}\,\left\|\nabla\de \oline u\right\|_{L^2}
\,+\,\left\|\big(\nabla\oline w_1,\nabla\oline w_2\big)\right\|_{L^{p}}\,\left\|\de\oline\theta\right\|_{L^{q}}\,
\left\|\nabla\de \oline u\right\|_{L^2} \\
&\quad +\,\left\|\big(\oline\theta_2,\nabla\oline v_1,\nabla\oline v_2\big)\right\|_{L^p}\,\left\|\big(\nabla\oline u_1,\nabla\oline u_2\big)\right\|_{L^\infty}\,
\left\|\de\oline\theta\right\|_{L^q}\,\Big(\left\|\de J^{-1}\right\|_{L^2}+\left\|\de\adj\right\|_{L^2}+\left\|\de A\right\|_{L^2}\Big) \\
&\qquad\qquad\qquad
+\, \left\|\big(\nabla\oline w_1,\nabla\oline w_2\big)\right\|_{L^{p}}\,\left\|\big(\nabla\oline u_1,\nabla\oline u_2\big)\right\|_{L^\infty}\,
\left\|\de\oline\theta\right\|_{L^{q}}\,\Big(\left\|\de J^{-1}\right\|_{L^2}+\left\|\de A\right\|_{L^2}\Big)\,,
\end{align*}
and an application of Lemma \ref{l:est_delta_A-J} yields
\begin{align*}
\mf S\,&\lesssim\,
\left\|\big(\oline\theta_2,\nabla\oline v_1,\nabla\oline v_2,\nabla\oline w_1,\nabla\oline w_2\big)\right\|_{L^p}\,
\left\|\de\oline\theta\right\|_{L^q}\,\left\|\nabla\de \oline u\right\|_{L^2} \\
&\qquad +\,\left\|\big(\oline\theta_2,\nabla\oline v_1,\nabla\oline v_2,\nabla\oline w_1,\nabla\oline w_2\big)\right\|_{L^p}\,
\left\|\big(\nabla\oline u_1,\nabla\oline u_2\big)\right\|_{L^\infty}\,
\left\|\de\oline\theta\right\|_{L^q}\,\left(\int^t_0\|\nabla\de\oline u\|_{L^2}\,\dd\tau\right)\,. 
\end{align*}
Now, similarly to what has been done in \eqref{est:dT_2}, we make use of the Young and Jensen inequalities
to get, for any $h>0$ as small as we want (to be chosen later) and for a multiplicative constant $C=C(h)>0$ also depending on the value of $h$,
the following bound:
\begin{align} \label{est:dT_4}
\mf S\,&\leq\,
C\,\left\|\big(\oline\theta_2,\nabla\oline v_1,\nabla\oline v_2,\nabla\oline w_1,\nabla\oline w_2\big)\right\|^2_{L^p}\,\left\|\de\oline\theta\right\|^2_{L^q} \\
\nonumber
&\qquad\qquad\qquad\qquad\qquad
+\,h\,\left\|\nabla\de \oline u\right\|^2_{L^2}\,+\,
t\,\left\|\big(\nabla\oline u_1,\nabla\oline u_2\big)\right\|^2_{L^\infty}\,\int^t_0\|\nabla\de\oline u\|^2_{L^2}\,\dd\tau\,.
\end{align}
Next, we make an interpolation on the $L^q$ norm of $\de\oline\theta$. Since $q>2$, an application of
Lemma \ref{LemmaUsefulInequalityGeneralVersion} yields
$ 
 \left\|\de\oline\theta\right\|_{L^q}\,\leq\,\left\|\de\oline\theta\right\|_{L^2}^\alpha\,\left\|\nabla\de\oline\theta\right\|_{L^2}^{1-\alpha}\,,
$ 
where, keeping in mind definition \eqref{def:q}, we have
\begin{equation} \label{def:alpha}
 \alpha\,=\,1\,-\,\frac{d}{p}\,. 
\end{equation}
By use of those interpolation inequalities, combined with an application of the Young inequality, from \eqref{est:dT_4} we find, for any $h>0$, the bound
\begin{align} \label{est:dT_5}
\mf S\,&\leq\,
C\,\left\|\big(\oline\theta_2,\nabla\oline v_1,\nabla\oline v_2,\nabla\oline w_1,\nabla\oline w_2\big)\right\|^{2/\alpha}_{L^p}\,\left\|\de\oline\theta\right\|^2_{L^2} \\
\nonumber
&\qquad\qquad
+\,h\,\left\|\nabla\de \oline u\right\|^2_{L^2}\,+\,h\,\left\|\nabla\de \oline \theta\right\|^2_{L^2}\,+\,
t\,\left\|\big(\nabla\oline u_1,\nabla\oline u_2\big)\right\|^2_{L^\infty}\,\int^t_0\|\nabla\de\oline u\|^2_{L^2}\,\dd\tau\,,
\end{align}
for a new multiplicative constant $C=C(h)>0$, blowing up for $h\to0^+$.

Finally, putting relations \eqref{est:dT_0}, \eqref{est:dT_1}, \eqref{est:dT_2}, \eqref{est:dT_3} and \eqref{est:dT_4} together, we deduce an energy estimate
for the quantity $\de\oline\theta$:
\begin{align}
\label{est:dT_final}
& \frac{1}{2}\,\frac{\dd}{\dt}\int\rho_0\,\left|\de\oline\theta\right|^2\,\dd y\,+\,\k\,\int\left|\nabla\de\oline\theta\right|^2\,\dd y \\
\nonumber
&\qquad \leq\,
\Big(C_0\,\k\,k_0\,+\,h\Big)\,\left\|\nabla\de\oline\theta\right\|_{L^2}^2\,+\,h\,\left\|\nabla\de \oline u\right\|^2_{L^2} \\
\nonumber
&\qquad\qquad +\,
C_0\,\k\,t\,\left\|\big(\nabla\oline{\theta}_2,\nabla\oline u_1,\nabla\oline u_2\big)\right\|^2_{L^\infty}\,\int^t_0\|\nabla\de\oline u\|^2_{L^2}\,\dd\tau \\
\nonumber
&\qquad\qquad\qquad +\,
C_h\,\left(\left\|\nabla\oline u_1\right\|_{L^\infty}\,+\,
\left\|\big(\oline\theta_2,\nabla\oline v_1,\nabla\oline v_2,\nabla\oline w_1,\nabla\oline w_2\big)\right\|^{2/\alpha}_{L^p}\right)\,
\left\|\de\oline\theta\right\|_{L^2}^2\,,
\end{align}
where we tacitly agree on the following notation: the constant $C_0>0$ depends only on the parameters of the problem, but not on $h$, while
$C_h$ depends also on the value of $h$.

\subsubsection{The final energy estimate} \label{sss:final-en}
At this point, we put inequalities \eqref{est:du} and \eqref{est:dT_final} together.
By taking $k_0$ small enough (that is, $T_\bull>0$ small enough) and $h>0$ small enough, we infer
\begin{align*}
& \frac{\dd}{\dt}\left\|\sqrt{\rho_0}\, \delta{ \overline u}\right\|^2_{L^2}\, +\,\frac{\dd}{\dt}\left\|\sqrt{\rho_0}\, \delta{ \overline \theta}\right\|^2_{L^2}\, +\,
\left\|\nabla \de\oline u\right\|^2_{L^2}\,+\,\left\|\nabla \de\oline \theta\right\|^2_{L^2} \\
\nonumber
& \qquad\qquad\qquad \lesssim\,
t\,\Big(1\,+\,\left\|\big(\oline\theta_1,\oline\theta_2,\nabla\oline\theta_1,\nabla\oline{\theta}_2,\nabla\oline u_1,\nabla\oline u_2\big)\right\|^2_{L^\infty}\Big)\,
\int^t_0\|\nabla\de\oline u\|^2_{L^2}\,\dd\tau \\
&\qquad\qquad\qquad\qquad\qquad\quad +\,
\left(1\,+\,\left\|\nabla\oline u_1\right\|_{L^\infty}\,+\,
\left\|\big(\oline\theta_2,\nabla\oline v_1,\nabla\oline v_2,\nabla\oline w_1,\nabla\oline w_2\big)\right\|^{2/\alpha}_{L^p}\right)
\,\left\|\de\oline\theta\right\|_{L^2}^2\,,
\end{align*}
for a universal multiplicative constant, which now depends only on $\mu$ and $\k$  (and on the value of $h$, which is now fixed once for all).

We now define the energy of the difference of the solutions as
\[
 \mc D(t)\,:=\,\left\|\sqrt{\rho_0} \, \delta{ \overline u}(t) \right\|^2_{L^2} +\,\left\|\sqrt{\rho_0}\, \delta{ \overline \theta}(t) \right\|^2_{L^2}\,+\,
\int^t_0\left\|\nabla \delta{ \overline u}(\t) \right\|^2_{L^2}\,\dd\t\,+\,\int^t_0\left\|\nabla \de\oline\theta(\t)\right|^2\,\dd\t
\]
for any $t\in [0,T_\bull]$, with the convention that
\[
 \mc D(0)\,\equiv\,\mc D_0\,:=\,\left\|\sqrt{\rho_0}\, \delta{ \overline u}_0 \right\|^2_{L^2} +\,\left\|\sqrt{\rho_0}\, \delta{ \overline \theta}_0 \right\|^2_{L^2}\,.
\]

From the previous inequality, using the fact that $\rho_0\,=\,1\,+\,\z_0\,\geq\,1\,-\,\veps_0\,>\,0$, an application of the Gr\"onwall lemma yields,
for any time $t\in[0,T_\bull]$, the estimate
\begin{align} \label{est:Gronwall}
\mc D(t)\,\leq\,\mc D_0\,\exp\left(C\,\int^t_0\Big(\Gamma_1(\t)\,+\,\Gamma_2(\t)\Big)\,\dd\t\right)\,,
\end{align}
where, for $j=1,2$, we have defined the function
\begin{align*}
 \Gamma_j(t)\,&:=\,t\,\Big(1+\left\|\big(\oline\theta_j,\nabla\oline\theta_j,\nabla\oline u_j\big)(t)\right\|^2_{L^\infty}\Big)\, +\,
 \left(1+\left\|\nabla\oline u_j(t)\right\|_{L^\infty}+\left\|\big(\oline\theta_j,\nabla\oline v_j,\nabla\oline w_j\big)(t)\right\|^{2/\alpha}_{L^p}\right)\,.
\end{align*}
Recall that the exponent $\alpha\in\,]0,1[\,$ has been defined in \eqref{def:alpha} as $\alpha\,=\,1\,-\,d/p$.

Thanks to estimate \eqref{est:Gronwall} and Lemma \ref{l:equiv-L^p},
it is clear that the uniqueness statement boils down to proving the following statement.
As it will appear clear from its proof, this result makes crucial use of the additional
assumptions imposed in the last part of Theorem \ref{th:well-p}.

\begin{prop} \label{p:G-integr}
Let the general assumptions of Theorem \ref{th:well-p} be in force.
Reinforce the lower bound for  $1/r$ in \eqref{HypForPs} with the assumption $1/r>7/8$ and suppose $p\,\in\;\big]2d\,,\,4d\big]$.
Let the triplet $\big(\rho_0,u_0,\theta_0\big)$ satisfy
the conditions stated in the existence part, plus the properties $w_0\,\in\, B^{2-2/R}_{p,R}$
and $\theta_0\,\in\,B^{2-2/S}_{p/2,S}$, with $R$ and $S$ defined as in the uniqueness part.
Let $\big(\rho,u,\theta\big)$ a solution related to this initial datum, defined on $[0,T]\times\R^d$ (for some $T>0$)
and satisfying the properties listed in the existence part Theorem \ref{th:well-p}.

Then, if we define $\alpha\in\,]0,1[\,$ as in \eqref{def:alpha} and, for $t\in[0,T]$, the function
\[
\Gamma(t)\,:=\,t\,\Big(1\,+\,\left\|\big(\theta,\nabla\theta,\nabla u\big)(t)\right\|^2_{L^\infty}\Big)\,+\,
 \left(1\,+\,\left\|\nabla u(t)\right\|_{L^\infty}\,+\,
\left\|\big(\theta,\nabla v,\nabla w\big)(t)\right\|^{2/\alpha}_{L^p}\right)\,,
\]
one has the property $\Gamma\,\in\,L^1\big([0,T]\big)$.
\end{prop}

The proof of Proposition \ref{p:G-integr} is the matter of the next subsection.

\subsection{Proof of Proposition \ref{p:G-integr}} \label{ss:proof-prop-u}

Here, we show that the property $\Gamma\in L^1\big([0,T]\big)$.
We start by treating the terms which are already under control within the assumptions of the existence theory, 
supplemented with the additional condition $p>2d$.

\begin{lemma} \label{l:u-by-exist}
Let the hypotheses of the existence part of Theorem \ref{th:well-p} be in force. Assume moreover the condition $p>2d$.
Let $\alpha\in\,]0,1[\,$ be defined as in \eqref{def:alpha}.

Then, one has the properties $\left\|\nabla u\right\|_{L^\infty}\,\in\,L^1\big([0,T]\big)$ and
$\left\|\theta\right\|_{L^\infty}\,\in\,L^2\big([0,T]\big)$. In addition,
one has
$\left\|\big(\theta,\nabla v\big)\right\|_{L^p}\,\in\,L^{2/\alpha}\big([0,T]\big)$
and $\left\|\nabla v\right\|_{L^\infty}\,\in\,L^2\big([0,T]\big)$.
\end{lemma}

\begin{proof}
First of all, we notice that $\nabla u$ belongs to $L^1_T(L^\infty)$, as a consequence of Corollary \ref{CorRegularityu}.
Moreover, $\theta$ and $\nabla v$ both belong to $L^\infty_T(L^p)$, owing to estimate \eqref{est:temp_L^p} and Lemma \ref{LemmaRegularityv}.

Next, we use that $\theta\in L^\infty_T(L^{p/2})$ and $\Delta\theta\in L^s(L^{p/2})$ to get, by interpolation,
$\nabla\theta\in L^{2s}_T(L^{p/2})$; since $p>2d$ and $2s>2$, by Sobolev embeddings we immediately deduce that
$\theta\in L^{2s}_T(L^\infty)\,\hookrightarrow\,\theta\in L^2_T(L^\infty)$.

Finally, let us focus on the term $\nabla v\,=\,-\,\nabla^2\big(\Id-\Delta\big)^{-1}P(\rho,\theta)$. Keep in mind definition \eqref{DefinitionP} of the pressure function,
with $\left\|\pi_0(\rho)\right\|_{L^p\cap L^\infty}\,\lesssim\,\|\z\|_{L^p\cap L^\infty}$ (recall that $\z=\rho-1$)
under the smallness property \eqref{TesisLInftyZeta}, satisfied by the given solution.
Thus, on the one hand we have
\[
 \left\|\nabla^2\big(\Id-\Delta\big)^{-1}\pi_0(\rho)\right\|_{L^\infty}\,\in\,L^\infty\big([0,T]\big)\,,
\]
as a consequence of the first estimate in Proposition \ref{PropositionEstimationForTangentialRegularity}.
On the other hand, using again the condition $p>2d$, the above proved property $\theta\in L^{2s}_T(L^\infty)$ and Lemma \ref{LemmaRegularitytheta}, the same
estimate from Proposition \ref{PropositionEstimationForTangentialRegularity} (applied with $p/2$ this time) implies that
\[
 \left\|\nabla^2\big(\Id-\Delta\big)^{-1}\Big(\pi_1(\rho)\,\theta\Big)\right\|_{L^\infty}\,\in\,L^{2s}\big([0,T]\big)\,\hookrightarrow\,L^2\big([0,T]\big)\,.
\]
We conclude that $\left\|\nabla v\right\|_{L^\infty}\,\in\,L^2\big([0,T]\big)$, as claimed.
This ends the proof of the lemma.
\end{proof}

With these properties at hand, it remains us to prove that $\nabla\theta$ and $\nabla w$ both belong to $L^2_T(L^\infty)$ and that
$\nabla w$ belongs to $L^{2/\alpha}_T(L^p)$.
Before showing those properties, we need some preparation. First of all,
we establish a low frequency control on the effective viscous flux $w$.

\begin{lemma} \label{l:w-L^p}
Let the assumptions of the existence part of Theorem \ref{th:well-p} be in force,
together with the condition $p>2d$.  
If moreover $w_0\in L^p(\R^d)$, then one has
\[
w\,\in\,L^{\infty}\big([0,T];L^{p}(\R^d)\big)\,.
\]
In particular, one also deduces that $u\,\in\,L^{\infty}\big([0,T];L^{p}(\R^d)\big)$.
\end{lemma}

\begin{proof}
 We observe that the second equation in system \eqref{FNSwithTemperatureModified}, can be written as
\[
 \d_tw\,-\,\mu\,\Delta w\,=\,-\,\z\,\d_tw\,+\,\rho\,F\,+\,\big(\mu\,+\,\lambda\big)\,\nabla\div w\,,
\]
where we have also used the specific form \eqref{DefinitionL} of the operator $\mc L$. Thus, by assumption on the solution triplet $\big(\rho,u,\theta\big)$,
each component of $w$ solves a heat equation with right-hand side belonging to $L^r_T(L^p)$, so in particular to $L^1_T(L^p)$.
Since the initial datum $w_0$ belongs to $L^p$ by, by standard parabolic theory (see the computations leading to \tsl{a priori} estimates in
\cite{F-FD-MM}, for instance) one infers that every component of $w$ belongs to $L^\infty_T(L^p)$, as claimed.

Finally, using that $v\in L^\infty_T(L^p)$ owing to Lemma \ref{LemmaRegularityv}, we see that also $u$ belongs to the same space.
This concludes the proof of the lemma.
\end{proof}

The next results is a first step in the improvement of the time integrability of the temperature function.
Remark that the index $\wtilde r$ appearing in this statement satisfies $\wtilde r< S$, so that
$B^{2-2/S}_{p/2,S}\,\hookrightarrow\,B^{2-2/\wtilde r}_{p/2,\wtilde r}$ and the assumptions are then consistent with the ones of Proposition \ref{p:G-integr}.

\begin{lemma} \label{l:integr-D^2-theta}
Let the assumptions of the existence part of Theorem \ref{th:well-p} be in force, together with the conditions $p>2d$ and $w_0\in L^p(\R^d)$.
Assume in addition that $\theta_0\in B^{2-2/\wtilde r}_{p/2,\wtilde r}$, with $\wtilde r$ defined as
$\wtilde r\,:=\,r_1/2$.
Then one has the properties
\[
\d_t\theta\,,\,\Delta \theta\;\in\,L^{\wtilde r}\big([0,T];L^{p/2}(\R^d)\big)\,.
\]
In particular, one also has
\[
 \nabla\theta\,\in\,L^{2\wtilde r}\big([0,T];L^{p/2}(\R^d)\big) 
\qquad \mbox{ and }\qquad \theta\,\in\,L^{2\wtilde r(p-d)/p}\big([0,T];L^{\infty}(\R^d)\big)\,.
\]
\end{lemma}

\begin{proof}
We argue as in the proof of Lemma \ref{lem:RHS_MR_theta}.
By using that $\nabla u\in L^{r_1}_T(L^p)$ (keep in mind Corollary \ref{CorRegularityu}),
since $\wtilde r= r_1/2$, we immediately infer that $\mc B(\nabla u)$ belongs to $L^{\wtilde r}_T(L^{p/2})$, as claimed.
Next, combining the same property $\nabla u\in L^{r_1}_T(L^p)$ with $\theta\in L^\infty_T(L^p)$ and $\rho\in L^\infty([0,T]\times\R^d)$, we get
$\pi_1(\rho)\,\theta\,\div u\in L^{r_1}_T(L^{p/2})$, so in particular it belongs to $L^{\wtilde r}_T(L^{p/2})$.
The control of the term $\rho\,u\cdot\nabla\theta$, is slightly more involved: let us present it in detail.

First of all, we write $u\,=\,v+w$. Since $v\in L^\infty_T(W^{1,p})\,\hookrightarrow\,L^\infty_T(L^\infty)$ owing to Lemma \ref{LemmaRegularityv},
from the fact that $\nabla\theta$ belongs to $L^{2s}_T(L^{p/2})$ we find that
\[
 \rho\,v\cdot\nabla\theta\,\in\,L^{2s}\big([0,T];L^{p/2}(\R^d)\big)\,.
\]
For the corresponding term involving $w$, we first use the interpolation inequality
\[
 \left\|\nabla\theta\right\|_{L^p}\,\lesssim\,\left\|\theta\right\|_{L^p}^{1-\beta}\,\left\|\nabla^2\theta\right\|_{L^{p/2}}^{\beta}\,,\qquad\qquad
\mbox{ with }\qquad \beta\,=\,\frac{p}{2p-d}\,,
\]
together with the properties $\theta\in L^\infty_T(L^p)$ and $\Delta\theta\in L^s_T(L^{p/2})$, to get that
 $\nabla\theta$ belongs to $L^{s/\beta}_T(L^p)$. 
Thus, since $w\in L^\infty_T(L^p)$ by Lemma \ref{l:w-L^p} and $\rho\in L^\infty\big([0,T]\times\R^d\big)$, we obtain that
\[
 \rho\,w\cdot\nabla\theta\,\in\,L^{s/\beta}\big([0,T];L^{p/2}(\R^d)\big)\,.
\]
At this point, we observe that, thanks to the assumption $p>2d$, one has
\begin{equation} \label{est:parameters}
 \frac{\beta}{s}\,\leq\,\frac{1}{\wtilde r}\,.
\end{equation}
As a matter of fact, by the explicit definition of $\beta$ and $\wtilde r$ and by the fact that $1/s<1/r$, we see that, for inequality \eqref{est:parameters}
to be true, it is enough that
\[
 \frac{p}{2p-d}\,\frac{1}{r}\,\leq\,\frac{2}{r}\,-\,1\,,\qquad \mbox{ that is }\qquad
1\,\leq\,\frac{1}{r}\,\left(2\,-\,\frac{p}{2p-d}\right)\,.
\]
Using the lower bound for $1/r$ appearing in \eqref{HypForPs}, we see that this latter inequality is automatically satisfied whenever 
\[
1\,\leq\, \left(1\,-\,\frac{d}{2p}\right)\,\left(2\,-\,\frac{p}{2p-d}\right)\,=\,\frac{3p-2d}{2p}\,,
\]
which can be easily verified to be equivalent to $p\geq 2d$. Now, this requirement is fulfilled for free as, by assumption, we have $p>2d$.
This in turn implies that \eqref{est:parameters} holds true, hence in particular we get the sought property
$\rho\,w\cdot\nabla\theta\,\in\,L^{\wtilde r}\big([0,T];L^{p/2}(\R^d)\big)$.

All in all, we have proved that the terms
$\mc B(\nabla u)$, $\pi_1(\rho)\,\theta\,\div u$ and $\rho\,u\cdot\nabla\theta$ all belong to $L^{\wtilde r}\big([0,T];L^{p/2}(\R^d)\big)$.
Next, we observe that the assumption $\theta_0\in B^{2-2/\wtilde r}_{p/2,\wtilde r}$ in particular implies that $\theta_0$ also belongs
to the homogeneous space $\dot B^{2-2/\wtilde r}_{p/2,\wtilde r}$. Owing to Proposition \ref{PropositionBahouri}, this is equivalent to having
 $e^{t\Delta}\,\theta_0\;\in\,L^{\wtilde r}\big([0,T];L^{p/2}(\R^d)\big)$. 
Thus, by maximal regularity estimates of Lemma \ref{LemmaMaximalRegularitySecondDerivative}, we finally deduce that
$\Delta\theta$ and $\d_t\theta$ both belong to the same space $L^{\wtilde r}_T(L^{p/2})$, as claimed in the statement.

To conclude the proof of the lemma, we now go back to the property $\theta\in L^\infty_T(L^{p/2})$: this, combined with
$\Delta\theta\in L^{\wtilde r}_T(L^{p/2})$, yields $\nabla\theta\in L^{2\wtilde r}_T(L^{p/2})$ by interpolation.
In addition, by interpolating again this last property with the $L^\infty_T(L^p)$ bound for $\theta$, we easily get the
sought time integrability of $\theta$ in $L^\infty$.
The lemma is thus finally proven.
\end{proof}

The next result is the crucial step in the proof of Proposition \ref{p:G-integr}. It deals with improved integrability for the effective velocity $w$
and its derivatives. Similarly as before, we notice that the index $\wtilde r$ also satisfies $\wtilde r< R$, so one has
$B^{2-2/R}_{p,R}\,\hookrightarrow\,B^{2-2/\wtilde r}_{p,\wtilde r}$; thus the assumptions below are consistent with the ones of Proposition \ref{p:G-integr}.

\begin{lemma} \label{l:integr-time}
Let the hypotheses of the existence part of Theorem \ref{th:well-p} be in force, together with the condition $p>2d$.
After setting $\wtilde r\,:=\,r_1/2$, assume in addition that $\theta_0\in B^{2-2/\wtilde r}_{p/2,\wtilde r}$ and that $w_0\in B^{2-2/\wtilde r}_{p,\wtilde r}$.
Then one has
\[
\d_tw\,,\,\nabla^2w\;\in\,L^{\wtilde r}\big([0,T];L^p(\R^d)\big)\,.
\]
In addition, define $\wtilde r_{0,\infty}$, $\wtilde{r}_1$ and $\wtilde r_{1,\infty}$ by the formulas
\[
 \frac{1}{\wtilde r_1}\,:=\,\frac{1}{\wtilde r}\,-\,\frac{1}{2}\,,\qquad 
\frac{1}{\wtilde r_{1,\infty}}\,:=\,\frac{1}{\wtilde r}\,-\,\frac{1}{2}\,+\,\frac{d}{2p} 
 \qquad \mbox{ and }\qquad
\frac{1}{\wtilde r_{0,\infty}}\,:=\,\frac{1}{\wtilde r}\,-\,1\,+\,\frac{d}{2p}\,.
\]
Then one has the property $\nabla w\,\in\,L^{\wtilde r_1}\big([0,T];L^p(\R^d)\big)\,\cap\,L^{\wtilde r_{1,\infty}}\big([0,T];L^\infty(\R^d)\big)$.
together with $w\,\in\,L^{\wtilde r_{0,\infty}}\big([0,T];L^\infty(\R^d)\big)$.
\end{lemma}

\begin{proof}
We begin the proof by observing that, as it was the case for the temperature equation, the assumption on the initial datum $w_0$ is made
precisely to have $e^{t\Delta}\,w_0\;\in\,L^{\wtilde r}_T(L^{p})$.
Thus, by splitting the time derivative term as $\rho\,\d_tw\,=\,\d_tw\,+\,\z\,\d_tw$ and seeing the second term as a small perturbation
owing to \eqref{EstZeta2}, by maximal regularity
we only have to verify that the right-hand side of \eqref{EquationForw} belongs to $L^{\wtilde r}_T(L^p)$. Since $\rho$ is uniformly bounded
over $[0,T]\times\R^d$, it is enough to do so for the vector field $F$, defined in \eqref{DefinitionF}. We will essentially repeat the arguments
carried out in Subsection \ref{SubsectionMaximalRegularityForVelocity}.

First of all, thanks to Lemma \ref{l:integr-D^2-theta}, we imediately see that the term $\nabla(\Id-\Delta)^{-1}\left(\frac{\pi_1(\rho)}{\rho}\,G\right)$
belongs to the sought space $L^{\wtilde r}_T(L^p)$. The next term is $\frac{1}{\rho}\,\nabla(\Id-\Delta)^{-1}P\,=\,-\,\frac{1}{\rho}\,v$,
which is even easier to treat, because it belongs to $L^\infty_T(L^p)$.
We then proceed to bound the term
$\nabla(\Id-\Delta)^{-1}\div(P\,u)$: by arguing exactly as in \eqref{MaximalRegularityForWSecondEstimation}, we see that it
belongs to $L^{r_{0,\infty}}_T(L^p)$, and it follows directly from the definitions that $r_{0,\infty}\,>\,\wtilde r=r_1/2$.

We now focus on the term $\nabla(\Id-\Delta)^{-1}\Big(\big(g_0(\rho)+\theta g_1(\rho)\big)\,\div u\Big)$. For this, we use
Corollary \ref{ClaimEstimationForOperatorPGeneralp} twice, as done before \eqref{MaximalRegularityForWThirdEstimation}, with the additional observation
that $\div u$ in fact belongs to $L^{r_1}_T(L^p)$, with of course $r_1>\wtilde r=r_1/2$.
Finally, the estimate of the non-linear term $u\cdot\nabla u$ does not change either: we use the properties $u\in L^{r_{0,\infty}}_T(L^\infty)$
and $\nabla u\in L^{r_1}_T(L^p)$, together with the fact that $1/r_{0,\infty}\,+\,1/r_1\,\leq\,1/\wtilde r$.

By maximal regularity, we then deduce $\d_tw,\nabla^2w\,\in L^{\wtilde r}_T(L^p)$, which yields the first conclusion of the Lemma.
At this point, the bound for the first order derivative is an easy consequence. As a matter of fact, from Lemma \ref{l:w-L^p}
and the fact that $\nabla w\in L^{r_1}_T(L^p)$, we infer that $w\in W^{1,\wtilde r}_T(L^p)\cap L^{\wtilde r}_T(W^{2,p})$.
It is then a direct consequence of the functional inequality (2.14) of \cite{B-F-P} (itself derived from the analogous estimate on bounded domains \cite{Sol04}) that
$\nabla w\in L^{\wtilde r_1}_T(L^p)$. Moreover, by interpolating this last property with $\nabla^2w\in L^{\wtilde r}_T(L^p)$,
exactly as done in Lemma \ref{LemmaRegularityw}, we finally get also the claimed $L^{\wtilde r_{1,\infty}}_T(L^\infty)$
bound on $\nabla w$.

Concerning $w$, its $L^{\wtilde r_{0,\infty}}_T(L^\infty)$ bound is a consequence of the Gaussian bheaviour of the fundamental solution to the heat equation;
we refer to formula (2.13) in \cite{B-F-P} for details (see also \cite{Aron, Nash_par}).
The lemma is now completely proven.
\end{proof}

With the above result at hand, we can easily conclude the time integrability of the terms involving $\nabla w$, which appear in the definition
of the function $\Gamma$.

\begin{cor} \label{c:integr-G}
Under the assumptions of Lemma \ref{l:integr-time}, one has
$\nabla w\,\in\,L^{2/\alpha}\big([0,T];L^p(\R^d)\big)\,\cap\,L^{2}\big([0,T];L^\infty(\R^d)\big)$, with $\alpha\in\,]0,1[\,$ defined in \eqref{def:alpha}.
\end{cor}

\begin{proof}
 It is a direct consequence of our definitions and of conditions \eqref{HypForPs} that
 \[
\frac{1}{\wtilde r_1}\,=\,\frac{1}{\wtilde r}\,-\,\frac{1}{2}\,=\,\frac{2}{r}\,-\,\frac{3}{2}\,\leq\,\frac{\alpha}{2}\qquad \mbox{ and }\qquad
\frac{1}{\wtilde r_{1,\infty}}\,=\,\frac{1}{\wtilde r}\,-\,\frac{1}{2}\,+\,\frac{d}{2p}\,\leq\,\frac{1}{2}\,.
 \]
This proves the corollary.
\end{proof}

As a last step, we prove the right time integrability of the term depending on the $L^\infty$ norm of $\nabla\theta$. Doing this, however, requires
an intricate argument, where we need to iterate the above procedure (the one starting from Lemma \ref{l:integr-D^2-theta}) twice.

\begin{lemma} \label{l:D-theta_inf}
Under the assumptions of Proposition \ref{p:G-integr}, one has the property
\[
\nabla \theta\;\in\,L^{2}\big([0,T];L^\infty(\R^d)\big)\,.
\]
\end{lemma}

\begin{proof}
Our starting point consists in going back to the proof of Lemma \ref{l:integr-D^2-theta}. Using the properties of Lemma \ref{l:integr-time} together with
Lemma \ref{LemmaRegularityv}, we see that $\nabla u\in L^{\wtilde r_1}_T(L^p)$, whence both
$\mc B(\nabla u)$ and $\pi_1(\rho)\,\theta\,\div u$ belong to the space $L^{\wtilde r_1/2}_T(L^{p/2})$. Next, we write
$u=v+w$ and notice that $u\in L^{\wtilde r_{0,\infty}}_T(L^\infty)$, thanks to Lemmas \ref{LemmaRegularityv} and \ref{l:integr-time}.
Observe that $\wtilde r_1/2\leq 2\wtilde r$ whenever $1/r\geq 5/6$, which is always true under our assumptions
(as $1/r>7/8$).
Then, using Lemma \ref{l:integr-D^2-theta} and further remarking that $1/\wtilde r_{0,\infty}\,+\,1/\wtilde r_1\,\leq\,2/\wtilde r_1$,
we deduce that $\rho u\cdot\nabla\theta\in L^{\wtilde r_1/2}_T(L^{p/2})$.

At this point, set ${\wtilde r}_2\,:=\,\wtilde r_1/2$. Observe that ${\wtilde r}_2=R> \wtilde r$.
From the computations above, we deduce that $\theta$ satisfies a heat equation with right-hand side belonging to $L^{\wtilde r_2}_T(L^{p/2})$.
As $S>R$ by assumption, by Besov embeddings we find that the initial datum
$\theta_0$ belongs to the space $B^{2-2/\wtilde r_2}_{p/2,\wtilde r_2}$, so also to its homogeneous counterpart $\dot B^{2-2/\wtilde r_2}_{p/2,\wtilde r_2}$.
Hence, by maximal regularity estimates we deduce
\begin{equation} \label{est:interm-theta}
 \d_t\theta\,,\Delta\theta\;\in\,L^{\wtilde r_2}\big([0,T];L^{p/2}\big)\qquad \mbox{ and } \qquad
\nabla\theta\;\in\,L^{2\wtilde r_2}\big([0,T];L^{p/2}\big)\,,
\end{equation}
the latter property following by interpolation of the former with the $L^\infty_T(L^{p/2})$ bound on $\theta$.

Next, we repeat exactly the same argument of the proof of Lemma \ref{l:integr-time}, to get that all the terms appearing in the definition of $F$
in \eqref{DefinitionF} belong to $L^{\wtilde r_2}_T(L^p)$. The only term which presents some difficulty is the one involving $\div (P\,u)$, which we
write as
\[
 \nabla(\Id-\Delta)^{-1}\div(P\,u)\,=\,\nabla(\Id-\Delta)^{-1}\div\big(\pi_0(\rho)\,u\big)\,+\,\nabla(\Id-\Delta)^{-1}\div\big(\pi_1(\rho)\,\theta\,u\big)\,.
\]
For the former term appearing on the right-hand side we use that $\rho$ is bounded over $[0,T]\times\R^d$, together with the bound $u\in L^\infty_T(L^p)$.
For the second term, instead, we use that $\theta\in L^\infty_T(L^p)$ and that $u\in L^{\wtilde r_{0,\infty}}_T(L^\infty)$,
as it follows from Lemmas \ref{l:integr-time} and \ref{LemmaRegularityv}, with (as already observed)
$1/\wtilde r_{0,\infty}<1/\wtilde r_{0,\infty}\,+\,1/\wtilde r_1\,\leq\,2/\wtilde r_1=1/\wtilde r_2$.
Combining these bounds with the property $w_0\in B^{2-2/\wtilde r_2}_{p,r_2}$, we immediately infer that
\[
 \d_tw\,,\,\nabla^2w\;\in\,L^{\wtilde r_2}\big([0,T];L^p(\R^d)\big)\,.
\]
Owing to these properties, we can conclude as done in the proof of Lemma \ref{l:integr-time} that
\[
\nabla w\,\in\,L^{\wtilde s_1}\big([0,T];L^p(\R^d)\big)
\qquad \mbox{ and }\qquad
w\,\in\,L^{\wtilde s_{0,\infty}}\big([0,T];L^\infty(\R^d)\big)\,,
\]
where we have defined the two integrability indices $\wtilde s_1$ and $\wtilde s_{0,\infty}$ by the usual formulas
\[
 \frac{1}{\wtilde s_1}\,:=\,\frac{1}{\wtilde r_2}\,-\,\frac{1}{2}\qquad \mbox{ and }\qquad
 \frac{1}{\wtilde s_{0,\infty}}\,:=\,\frac{1}{\wtilde r_2}\,-\,1\,+\,\frac{d}{2p}\,.
\]

As a last step, we come back again to the heat equation for $\theta$. Observe that $\nabla u$ enjoys the same time integrability properties as $\nabla w$
with values in $L^p$, thanks to Lemma \ref{LemmaRegularityv}.
Therefore, we can repeat \tsl{mutatis mutandis} the same computations above, to bound the terms $\mc B(\nabla u)$, $\pi_1(\rho)\,\theta\,\div u$
and $\rho\,u\cdot\nabla\theta$ in view of the new properties satisfied by $w$ and its derivatives and of the second condition appearing in \eqref{est:interm-theta}.
At this point, we remark that $\wtilde s_1/2\,=\,S$ by definition; furthermore, we have
\[
\frac{1}{\wtilde s_{0,\infty}}\,+\,\frac{1}{2\,\wtilde r_2}\,\leq\,\frac{2}{\wtilde{s}_1}\,=\,\frac{1}{S}\,.
\]
Indeed, by using the definitions, we see that this inequality is equivalent to have the lower bound $1/r\,\geq\, 3/4\, -\,d/(4p) $,
which is always satisfied under our assumptions, because $ 3/4\, -\,d/(4p) < 7/8< 1/r$.
Therefore, using also the assumption $\theta_0\in B^{2-2/S}_{p/2,S}$ on the initial datum, by maximal regularity again we find
\[
 \d_t\theta\,,\Delta\theta\;\in\,L^{S}\big([0,T];L^{p/2}(\R^d)\big)\qquad \mbox{ and } \qquad
\nabla\theta\;\in\,L^{2S}\big([0,T];L^{p/2}(\R^d)\big)\,.
\]
Now, using that $p>2d$, a control on the $L^\infty$ norm of $\nabla\theta$ comes from an interpolation of the above two quantities. In any case, we can always affirm that
\[
 \nabla\theta\,\in\,L^{S}\big([0,T];L^{\infty}(\R^d)\big)\,.
\]
To finish the proof, we check that having $1/S\,=\,2/\wtilde s_1\,\leq\,1/2$ reduces to the condition
$1/r\,\leq\,15/16$, which holds certainly true, thanks to \eqref{HypForPs}, whenever $p\leq 4d$.
\end{proof}


\appendix

\section{Equivalence of the two formulations} \label{app:equiv}

In this section, we carry out the computations leading from the original formulation \eqref{FNS} of the full compressible Navier-Stokes system to the
Navier-Stokes-Fourier system \eqref{eq:NSF} in the unkowns $\big(\rho,u,\theta\big)$.
We notice that, owing to the regularity properties satisfied by $u$ and $\theta$, the next computations are fully justified in our functional framework.

To begin with, we compute the equation for the kinetic energy. For doing this in a rigorous way, we test the momentum equation, \tsl{i.e.} the second equation
in \eqref{FNS}, against $\vphi\,u$, where $\vphi\in\mc C^\infty_0\big(\R_+\times\R^d\big)$ is a classical (scalar) test function.
The terms contained in $\mc Lu\cdot (\vphi u)$ can be rewritten in the following form, where the equal sign must be understood in the
weak sense:
\begin{align*}
-\,\mu\,\Delta u\cdot (\vphi\,u)\,&=\,\vphi\,\Big(-\,\mu\,\div\big(\nabla u\cdot u\big)\,+\,\mu\,|\nabla u|^2\Big) \\
-\,\lambda\,\nabla\div u\cdot (\vphi\,u)\,&=\,\vphi\,\Big(-\,\lambda\,\div\big(\div u\;u\big)\,+\,\lambda\,|\div u|^2\Big) \\
-\,\mu\,\nabla\div u\cdot (\vphi\,u)\,&=\,-\,\mu\sum_{j,k=1}^d\d_j\d_ku^k\,\vphi\,u^j\,=\,
-\,\mu\sum_{j,k=1}^d\d_k\Big(\d_ju^k\,u^j\Big)\,\vphi\,+\,\mu\,\sum_{j,k=1}^d\d_ju^k\,\d_ku^j\,\vphi \\
&=\,\vphi\,\Big(-\,\mu\,\div\big(Du\cdot u\big)\,+\,\mu\,\nabla u:\nabla u\Big)\,.
\end{align*}
In light of these computations, using also the mass conservation equation, we thus obtain
\begin{align}
	&\partial_t \left( \rho\, \frac{\vert u \vert^2}{2}\right)\,+\, \divergence \left( \rho\,\frac{\vert u \vert^2}2\, u \right)\,
	+\,\div\big(P\,u\big)\,-\, P\,\div u\, =\,\div\big(\tau\cdot u\big)\,-\,\mc B(\nabla u)\,, 
\label{KineticEnergy}
\end{align}
where we have resorted to the notation introduced in \eqref{def:B_u}.

By subtracting the above relation \eqref{KineticEnergy} from the third equation in \eqref{FNS}, we get
\begin{align} \label{eq:e}
\d_t(\rho\,e)\,+\,\div\big(\rho\,e\,u\big)\,+\,P\,\div u\,-\,\k\,\Delta\theta\,=\,\mc B(\nabla u)
\end{align}
To conclude, we recall that, by definition \eqref{definitione}, one has $e\,=\,\theta\,+\,h(\rho)$. Since $\rho$ solves the first equation
of \eqref{FNS} in the weak sense and $u$ belongs to $L^{m_1}\big([0,T];W^{1,\infty}(\R^d)\big)$, with $m_1>1$, we can apply \tsl{e.g.} Theorem 10.29 and Lemma 10.13
of \cite{F-N_book} to see that $\rho$ is also a normalised solution to this equation. Hence,
after defining $H(\rho)\,=\,\rho\,h(\rho)$, one has
\[
 \d_t\big( H(\rho)\big)\,+\,\div\big( H(\rho)\,u\big)\,+\,\big(\rho\,H'(\rho)\,-\,H(\rho)\big)\,\div u\,=\,0\,.
\]
Taking the diffence of \eqref{eq:e} by the above relation, and noticing that
\[
 \rho\,H'(\rho)\,-\,H(\rho)\,=\,\rho^2\,h'(\rho)\,+\,\rho\,h(\rho)\,-\,\rho\,h(\rho)\,=\,\pi_0(\rho)
\]
thanks to \eqref{RelationBetweenPih}, we finally get equation \eqref{eq:NSF}, as sought.

To conclude this part, we point out that, owing to the regularity properties of $u$, the above computations are reversible. Thus, disposing of a solution
to system \eqref{eq:NSF} at the level of regularity of Theorem \ref{th:well-p} provides us with a weak solution to the original
full Navier-Stokes system \eqref{FNS}.
This shows that, within the functional setting considered within this paper, the two formulations are equivalent.


\end{document}